\documentclass[a4paper,11pt,twoside]{article}
\usepackage{pst-fill,pst-grad}  
\usepackage{textcomp} 
\usepackage[english]{babel}      
\usepackage[utf8]{inputenc}    
\usepackage[titletoc]{appendix} 
\usepackage{titlesec}    
\usepackage{graphicx}        
\usepackage{amsmath} 
\usepackage{float} 
\usepackage{fancyhdr}  
\usepackage[matrix,arrow,curve]{xy}
\usepackage{pstricks} 
\usepackage{amsmath,amsfonts,verbatim,afterpage,theorem,euscript,mathrsfs,amssymb}
\usepackage{amsfonts} 
\usepackage{amssymb}
\usepackage{array}
\usepackage{dsfont}
\usepackage[colorlinks=true,linkcolor=blue,citecolor=red]{hyperref}
\usepackage{authblk} 
\usepackage{color} 
\newtheorem{Definition}{Definition}[section] 
 
\newtheorem{Proposition}{Proposition}[section]

\newtheorem{Lemme}{Lemma}[section]
\newtheorem{Theoreme}{Theorem}[section]

\newtheorem{Corollaire}{Corollary}[section]
\newtheorem{Remarque}{Remark}

\def \fe{\vec{f}} 

\def \vu{\vec{u}}
\def \vv{\vec{v}}
\def \vw{\vec{w}}
\def \P{\mathbb{P}}
\def \U{\vec{U}}

\def \R{\mathbb{R}}

\def \Rt{\mathbb{R}^3}

\def \finpv{\hfill $\blacksquare$}
\def \pv{{\bf{Proof.}}~} 

\def \ds{\displaystyle}

\def \a{{\bf a}}
\def \b{{\bf b}}
\def \c{{\bf c}}
\def \d{{\bf d}}

\title{\bf Asymptotic behavior  of a generalized Navier-Stokes-alpha  model and  applications to related models}  

\author[1]{ Oscar Jarr\'in\footnote{corresponding author: oscar.jarrin@udla.edu.ec}}  

\affil[1]{\scriptsize Escuela de Ciencias Físicas y Matemáticas, Universidad de Las Américas, Vía a Nayón, C.P.170124, Quito, Ecuador.} 
\begin{document} 
\maketitle
\begin{abstract}
We consider a generalized alpha-type model in  the whole three-dimensional space and driven by a stationary (time-independent) external force. This model   contains as  particular cases some relevant equations of the fluid dynamics, among them the Navier-Stokes-Bardina's model, the critical alpha-model, the fractional and the classical Navier-Stokes equations with an additional drag/friction term. First,  we study  the existence and in some cases the uniqueness of finite energy solutions. Then, we use a general framework  to study their long time behavior with respect to the \emph{weak} and the \emph{strong} topology of the phase space. When the uniqueness of solutions is known, we prove the existence of a \emph{strong} global attractor. Moreover, we  proof the existence of a \emph{weak} global attractor in the case when  the uniqueness of solutions is unknown. 

The weak/global attractor contains a particular kind of solutions to our model, so-called the stationary solutions. In all generality we construct these solutions,  and we study their uniqueness, \emph{orbital and  asymptotic stability} in the case when some physical constants in our model are large enough. As a bi-product, we show that   in some  cases  the weak/global attractor reduces down to the unique stationary solution. 
 \medskip 

\end{abstract} 

\section{Introduction} 
The study of the fluid dynamics provides us several evolution models of great importance,  among them, the well-known Navier-Stokes equations and some related equations  so-called  the  alpha-type models. The alpha-type models   have been developed in the mathematical literature as physically relevant  approximations of the Navier-Stokes equations, for instance,  the Bardina's model \cite{Bardina1,Bardina} the viscous Camassa-Holm model \cite{Chen}, the Leray-alpha model \cite{CheskidovHolm} and the Clark-alpha model \cite{Cao}. 

\medskip 

 Numerical solutions of the Navier–Stokes equations for problems of physically and  engineering  relevance is not possible at present  as the mathematical theory for uniqueness and  regularity of Leray's solutions  is one of the most challenging open questions \cite{Dubois,Pope}. Thus,  the alpha-type  models are regularized versions of the classical Navier-Stokes equations for which the global well-posedness of finite energy solutions can be solved.  Moreover, a deep comprehension of their  long time behavior is one of the \emph{key questions} to a better understanding of these
 models. The long time behavior of solutions has been extensively  studied when considering the alpha-type model with \emph{spatial periodic conditions} in the torus $\Omega = [0,L]^3$ (with $0<L$). See, for instance,   \cite{Cao,Cao1,Chai,CheskidovHolm,Chen,Chen2,Olsen-Titi} and the references therein.

\medskip

In this paper, we are interested in studying the long time behavior of solutions of the  following fractional alpha-type model which we shall consider on the \emph{whole space} $\Rt$. For the parameters $0<\alpha$, $0\leq \beta$, $0<\gamma$, $0<\delta$ and $0<\nu$,  this equation writes down as follows:

\begin{equation}\label{Equation}
\begin{cases}\vspace{1mm}
\partial_t \vu + \nu  (-\Delta)^{\frac{\alpha}{2}}\,  \vu +   (I_d - \delta^2\Delta)^{- \frac{\beta}{2}}\, \P  \left(  \text{div} (\vu \otimes \vu) \right)  = \fe -\gamma \vu, \qquad \text{div}(\vu)=0,\\
\vu(0,\cdot)=\vu_0.
\end{cases} 
\end{equation}
Here, $\vu : [0,+\infty[\times \Rt \to \Rt$ is the  velocity of the fluid. The  second equation $\text{div}(\vu)=0$ describes the  fluid's incompressibility. Moreover,  $0<\nu$ is the viscosity parameter,  the function $\fe : \Rt \to \Rt$ is the external force acting on the system, which is assumed as a time independent and divergence free vector field, while the function $\vu_0 : \Rt \to \Rt$ denotes the initial velocity field at the time $t=0$.   On the other hand,  the operator  $\P$ stands for the Leray's projector given by $\P(\vec{\varphi})=\vec{\varphi}- \vec{\nabla} \frac{1}{\Delta}(\text{div}(\vec{\varphi}))$.   We have applied the Leray's projector to this equation as the pressure term does not play any substantial role in our study.

\medskip

The main features of  equation (\ref{Equation}) are, on the one hand,  the fractional derivative  operators in both the linear diffusion term and the  nonlinear transport term and, on the other hand,  the damping term on the right-hand side. In what follows, we shall briefly explain the mathematically and physically relevance of these  terms. 

\medskip 

From the physical and the experimental point of view, the fractional diffusion term $(-\Delta)^{\frac{\alpha}{2}} \vu $ and the fractional transport term $\big( I_d - \delta^2 \Delta \big)^{-\frac{\beta}{2}} \left( \text{div}(\vu \otimes \vu) \right)$ have been successfully employed to model anomalous reaction-diffusion process in porous media models \cite{Meerschaert,Meerschaert2} and in computational turbulence models \cite{Pope}.  In these last models,  the operator  $\nu(-\Delta)^{\frac{\alpha}{2}}$ is used to  characterize anomalous viscous diffusion effects in turbulent fluids which are driven by the parameters $\alpha$ and $\nu$. On the other hand, the operator $\big( I_d - \delta^2 \Delta \big)^{-\frac{\beta}{2}}$, also known the Bessel potential \cite{Grafakos}, acts as  filtering-averaging operator: the parameter $\delta$ allow us to obtain an accurate model describing the large-scale motion of the fluid while filtering or averaging the fluid motion at small scales smaller than $\delta$. We refer to \cite{Pope} for more details. 

\medskip

From the mathematical point of view, the parameter $\alpha$ measures  the dissipative  degree, while the parameter $\beta$ affects the strength of the nonlinear effects. In this sense, one of the main mathematical interest of  equation (\ref{Equation}) is the study of how these parameters work together to yield a sharp description of the long time dynamics of solutions. 
 
 \medskip 

Equation (\ref{Equation}) was inspired by the following fractional   alpha-like model  introduced by E. Olsen and E.S. Titi in \cite{Olsen-Titi}: 
\begin{equation}\label{Model-Rel}
\partial_t \vu + \nu  (-\Delta)^{\frac{\alpha}{2}} \vu +  \text{div}   \left(  \left( \big( I_d -  \delta^2\Delta \big)^{-\frac{\beta}{2}}\, \vu\right)  \otimes \vu \right) + \vec{\nabla} p   = \fe.
\end{equation}
Mathematically, this model  has been  studied in the space periodic setting of the box $\Omega = [0,L]^3$.  As noticed in \cite{Olsen-Titi}, a weaker non-linearity and a stronger  dissipation  yield the global well - posedness of finite energy solutions. More precisely, when $5\leq 2\alpha +\beta$ the main result of \cite{Olsen-Titi} shows that  for any $\vu_0 \in H^{\frac{\beta}{2}}(\Omega)$ there exists a unique global in time weak solution $\vu \in L^{\infty}_{loc}([0,+\infty)), H^{\frac{\beta}{2}}(\Omega)) \cap L^{2}_{loc}([0,+\infty), H^{\frac{\alpha+\beta}{2}}(\Omega))$. 

\medskip

The long time behavior of these solutions  was then studied  in    \cite{Chai,TTT}. The main result states the existence of a \emph{strong} global attractor (see Definition \ref{Def-global-attractor}) for the equation (\ref{Model-Rel}). This result essentially bases on two key ideas: on the one hand, uniqueness of finite energy weak solutions allows to define a semigruop $S(t): H^{\frac{\beta}{2}}(\Omega) \to H^{\frac{\beta}{2}}(\Omega)$, where $S(t)\vu_0= \vu(t,\cdot)$ is the \emph{unique solution} of (\ref{Model-Rel}) arising from $\vu_0$.  On the other hand, the energy equality verified by these solutions and the \emph{Poincaré's inequality} yield to the following control in time: 
\begin{equation}\label{Control-Per}
\Vert \vu(t,\cdot) \Vert^{2}_{H^{\frac{\beta}{2}}(\Omega)} \lesssim \frac{L^2}{\nu^2} \, \Vert \fe \Vert^{2}_{H^{\frac{\beta}{2}}(\Omega)} , \qquad  t \to +\infty.
\end{equation}
The semigroup $(S(t))_{t\geq 0}$ and the  estimate above  are the key ingredients to apply some results in the  theory of dynamical systems (see Section \ref{Sec:strong-global-attractor} for more details) to yield the existence of a \emph{strong} global attractor. 
 
\medskip

Getting back to our model (\ref{Equation}), which is posed on the whole space $\Rt$,  to the best of  our knowledge  we are not able to obtain an analogous estimate to (\ref{Control-Per}), see for instance  \cite{Constantin2,FCortezOJarrin,Ilyn,Jarrin}, as the  \emph{Poincaré's inequality} is not longer valid. In this sense,   the damping term $-\gamma \vu$ mathematically acts as compensation of the  lack of the \emph{Poincaré's inequality}.  It is worth mentioning  another  damping terms can be considered to study the long-time behavior of  Navier-Stokes type  equations on the whole space \cite{Chamorro,Liua}. However, we will consider here the damping term $-\gamma \vu$ for its relevant  physical meaning: the parameter $0<\gamma$ is known as  the Rayleigh or Ekman friction coefficient; and the term $-\gamma \vu $  models the bottom friction in  ocean models and is the main energy sink in large scale atmospheric models  \cite{Pedlosky}. 

\medskip

Equation (\ref{Equation}) is also of interest since it contains as  particular case some relevant models. Consequently, our results also hold  for the following equations. See the Section \ref{Sec:Applications} for a more detailed discussion. When  we set $\alpha=\beta=2$,   equation (\ref{Equation}) agrees with the  damped Navier-Stokes-Bardina's model: 
\begin{equation}\label{Bardina-damped}
\partial_t \vu - \nu \Delta\,  \vu +   (I_d -\delta^2\Delta)^{- 1}\, \P  \left(  \text{div} (\vu \otimes \vu) \right)  = \fe -\gamma \vu, \qquad \text{div}(\vu)=0,
\end{equation}
previously studied in \cite{FCortezOJarrin}. In this sense the equation (\ref{Equation}) can be maned a \emph{generalized Navier-Stokes-Bardina's model}.   For the values $\alpha=2$ and $\beta = \frac{1}{2}$ we obtain a damped version of  the critical Leray-alpha model 
\begin{equation}\label{Critical}
\partial_t \vu - \nu \Delta\,  \vu +   (I_d - \delta^2\Delta)^{- \frac{1}{4}}\, \P  \left(  \text{div} (\vu \otimes \vu) \right)  = \fe -\gamma \vu, \qquad \text{div}(\vu)=0,
\end{equation}
for which the global well-posedness problem  was studied in \cite{Ali}, in the space periodic setting (when $\gamma=0$). Thereafter, for $0<\alpha$ and $\beta=0$ we get the following damped version of the fractional Navier-Stokes equations: 
\begin{equation}\label{Fractional-NS}
\partial_t \vu +  \nu(-\Delta)^{\frac{\alpha}{2}}\,  \vu +  \P \left(  \text{div} (\vu \otimes \vu)\right)   = \fe -\gamma \vu, \quad \text{div}(\vu)=0.
\end{equation}
This equation  has recently attired the attention of researchers in the mathematical fluid dynamics  to understand the dissipative effects (given by the   fractional Laplacian operator) in the study of  outstanding open  problems in   the  classical Navier-Stokes equations, for instance,  uniqueness and  regularity issues of Leray's weak solutions \cite{Cholewa,Colombo,Nan}. Finally, in the particular case when $\alpha=2$ and $\beta=0$, the equation (\ref{Equation})  deals with the classical damped Navier-Stokes equations: 
\begin{equation}\label{Damped-NS}
\partial_t \vu -\nu \Delta\,  \vu + \P \left(   \text{div} (\vu \otimes \vu) \right)   = \fe -\gamma \vu, \quad  \text{div}(\vu)=0.
\end{equation} 
In the setting of the whole space $\Rt$, when studying the large time behavior of solutions  this equation is an interesting \emph{ counterpart} of the classical (when $\gamma=0$) Navier-Stokes equations with \emph{space-periodic conditions}.  We refer to \cite{Constantin2} and \cite{Ilyn} for some interesting  previous related works on this equation. 

\medskip 

Once we have introduced the model (\ref{Equation}), we briefly summarize our main results.  In the Section \ref{Sec:Results} below  we make a detailed presentation  and discussion of them. We recall first that  the dissipative effects of the fractional Laplacian operator  are measured by the parameter $0<\alpha$, while the parameter $0\leq \beta$ measures the  regularizing effects of the Bessel potential in the nonlinear  transport term.  Therefore, the quantity $0<\alpha+\beta$ quantifies  the total contribution   of both dissipative and regularizing effects  in the qualitative study of the equation (\ref{Equation}).

\medskip

In all generality,  for $0<\alpha+\beta$ in Theorem \ref{Th1} we  construct global in time finite energy solutions for the equation (\ref{Equation}). These solutions belong to the energy space $(L^{\infty}_{t})_{loc} H^{\frac{\beta}{2}}_{x} \cap (L^{2}_{t})_{loc}H^{\frac{\alpha+\beta}{2}}_{x}$.  Moreover, in Proposition \ref{Prop1} we  study the effects of the damping term $-\gamma \vu$ and we obtain  useful controls in time on these solutions.  Thereafter, the study of their  long time  behavior  is divided in the following cases: 
\begin{enumerate}
\item[$\bullet$] When $\frac{5}{2}\leq \alpha+\beta$,  we  prove the uniqueness of finite energy solutions. Here,  the \emph{critical value} $\frac{5}{2}$ was also pointed out in prior  related works \cite{Ali,Olsen-Titi,TTT}. In this case,  uniqueness allows us to define a strongly continuous  semigroup $(S(t))_{t\geq 0}$ on the space $H^{\frac{\beta}{2}}(\Rt)$: for $\vu_0 \in H^{\frac{\beta}{2}}(\Rt)$ and for $0\leq t$ we have $S(t)\vu_0 = \vu(t,\cdot)$, where $\vu(t,\cdot)$ is the unique solution of the equation (\ref{Equation}) arising from $\vu_0$.  Uniqueness is also one of the key properties to show that  $(S(t))_{t\geq 0}$ is an  asymptotically compact semigroup in the strong topology of the space $H^{\frac{\beta}{2}}(\Rt)$ (see Definition  \ref{def-asymptotically-compact}). These properties   yield the existence of a \emph{strong} global attractor  for the equation (\ref{Equation}). See the first point of  Definition \ref{Def-global-attractor} and the first point of Theorem \ref{Th2} for more details. 

\item[$\bullet$] When $0<\alpha +\beta < \frac{5}{2}$,  uniqueness of finite energy solutions remains an outstanding open problem. In this case, we use a different approach to study their  long time behavior. More precisely, we introduce here the \emph{set} $R(t)\vu_0$  containing all the possibly finite energy solutions of  equation (\ref{Equation}) at the time $0<t$, which  arise from the initial datum $\vu_0$. Then, we are able to prove that  the family $(R(t))_{t\geq 0}$ is uniformly compact in the weak topology of the space $H^{\frac{\beta}{2}}(\Rt)$ (see Definition \ref{Def-Uniform-weak-compact}) and this fact  yields the existence of a \emph{weak} global attractor. See the second point of  Definition \ref{Def-global-attractor} and the second point of Theorem \ref{Th2} respectively.   

\medskip

To the best of our knowledge,  these results in the case $0 < \alpha+\beta <\frac{5}{2}$ have not been studied before
in the existent literature on fractional alpha-type models \cite{Ali,Chai,Holst,Ilyn1,Ilyn2,Olsen-Titi,TTT,Zhao}. On the other hand, it is  worth mentioning the  existence of a strong global attractor in the case when $0<\alpha +\beta < \frac{5}{2}$ is another open problem far from obvious. See \cite{CheskidovFoias} for a discussion in the case of the classical Navier-Stokes equations with periodic conditions.
\end{enumerate}	

The weak/global attractor can be precisely  characterized through the notion of eternal solutions to the equation (\ref{Equation}) (see expressions (\ref{Charac-weak}) and (\ref{Charac-strong}) below for more details). Thus, a simple but  key remark is that (when exist) stationary  solutions $\U$ to the equation (\ref{Equation}) belong to the weak/strong global attractor. Stationary solutions solve the elliptic problem:
\begin{equation}\label{Stationary}
\nu(-\Delta)^{\frac{\alpha}{2}} \U + (I_d-\delta^2\Delta)^{-\frac{\beta}{2}} \P \, \text{div}(\U \otimes \U) = \fe - \gamma \U, \qquad \text{div}(\U)=0.
\end{equation}
For the general case $0<\alpha+\beta$ and for any  (divergence-free) external  force $\fe \in H^{\frac{\beta}{2}}(\Rt)$, in Theorem \ref{Th3}  we  construct these solutions in the natural energy space $H^{\frac{\alpha+\beta}{2}}(\Rt)$. Then, for the range of values $2\leq \alpha+\beta$, in Theorem \ref{Th4}  we find some natural \emph{sufficient  conditions}, only depending on the external force $\fe$ and the parameters $\alpha, \beta, \gamma,\delta,\nu$  in  equation (\ref{Stationary}), which yield, on the one hand, the \emph{orbital stability} of stationary solutions and, on the other hand, a stronger result concerning the uniqueness and  \emph{asymptotic stability} of stationary solutions in the \emph{strong} topology of the space $H^{\frac{\beta}{2}}(\Rt)$. Consequently, in the case of uniqueness and asymptotic stability  we deduce that  weak/strong global attractor reduces down to the singleton $\{ \U \}$. In particular, for the range of values $2\leq \alpha+\beta < \frac{5}{2}$ the weak global attractor becomes  a strong one. 

\medskip

Finally, let mention that  our results are essentially  obtained by  energy methods, which   make them more interesting from a physical point of view since  we only control the natural energy quantities derived from equations (\ref{Equation}) and (\ref{Stationary}).

\medskip

{\bf Organization of the paper:}  In Section \ref{Sec:Results} we introduce some definition and  we present all our results. Section \ref{Sec:Leray-Solutions} is devoted to the study of the main features (existence and time controls) of finite energy solutions to the equation (\ref{Equation}). In Section \ref{Sec:Long-time-behavior} we focus on their long time asymptotic behavior through the notion of the weak/strong global attractor.  Section \ref{Sec:Stationary} is devoted to the aforementioned study  of stationary solutions. Finally, at Appendix \ref{Appendix-Fractal-dim} we derive an upper bound of the fractal dimension of the strong global attractor. 

\section{Definitions and presentation of the results}\label{Sec:Results}
 We have organized this section in three parts.
\subsection{Finite energy solutions}
 From now on, these solutions  shall be called the Leray-type solutions as they share the main properties of the well-known Leray's solutions in the classical Navier-Stokes theory.
\begin{Definition}[Leray-type solution]\label{Def-Leray} Let $0<\alpha$ and $0 \leq \beta$. We shall say that $\vu$ is a Leray-type solution of  equation (\ref{Equation}) if: 
	\begin{enumerate}
		\item  The function $\vu$ belongs to the energy space: $\ds{ L^{\infty}_{loc}\Big([0,+\infty[,H^{\frac{\beta}{2}}(\Rt)\Big) \cap L^{2}_{loc}\Big([0,+\infty[,H^{\frac{\alpha+\beta}{2}}(\Rt)\Big)}$,  and it  verifies the equation (\ref{Equation}) in the distributional sense.
		\item For all $0\leq t$, the following energy inequality holds: 
		\begin{equation*}  
		\begin{split}
		\Vert (I_d-\delta^2 \Delta)^{\frac{\beta}{4}} & \vu(t,\cdot)\Vert^{2}_{L^2}  \leq \, \Vert  (I_d-\delta^2 \Delta)^{\frac{\beta}{4}} \vu_0 \Vert^{2}_{L^2} -  \, 2 \nu \, \int_{0}^{t}  \left\Vert  (-\Delta)^{\frac{\alpha}{4}} (I_d-\delta^2 \Delta)^{\frac{\beta}{4}} \, \vu(s,\cdot)\right\Vert^{2}_{L^2}   ds \\
		& + 2 \int_{0}^{t} \left( (I_d-\delta^2 \Delta)^{\frac{\beta}{4}} \fe, (I_d-\delta^2 \Delta)^{\frac{\beta}{4}} \vu(s,\cdot) \right)_{L^2} \, ds -2 \gamma \int_{0}^{t}\Vert (I_d-\delta^2 \Delta)^{\frac{\beta}{4}} \vu(s,\cdot)\Vert^{2}_{L^2} ds,
		\end{split}
		\end{equation*}
		provided that $\vu_0 \in H^{\frac{\beta}{2}}(\Rt)$ and $\fe \in L^{2}_{loc}\Big([0,+\infty[, H^{\frac{\beta}{2}}(\Rt)\Big)$.
	\end{enumerate}		
\end{Definition} 
Existence and in some cases uniqueness  of Leray-type solutions is a rather standard issue. However, for the completeness of this article, we start by stating the following:
\begin{Theoreme}\label{Th1} Let $0<\alpha$, $0\leq \beta$.  Let $\vu_0 \in H^{\frac{\beta}{2}}(\Rt)$ be a divergence free initial datum. Moreover, let $\fe \in L^{2}_{loc}\Big([0,+\infty[, H^{\frac{\beta}{2}}(\Rt)\Big)$ be a divergence free external force.  Then, there exists $\vu$ a Leray-type solution of the equation (\ref{Equation}) given in Definition \ref{Def-Leray}. Moreover, if    $\frac{5}{2} \leq \alpha +\beta$ then the equation (\ref{Equation}) has a \emph{unique} Leray-type solution. 
\end{Theoreme}	
The value $\frac{5}{2}$ is the critical one obtained in related studies \cite{Cholewa,Lions,Nan} when studying the  uniqueness of Leray-type  solutions to the fractional Navier-Stokes equation (\ref{Fractional-NS}):  by setting $\beta=0$ we obtain the critical value $\frac{5}{2} \leq \alpha$. Moreover,  uniqueness of Leray-type solutions  in the \emph{supercritical range} $0<\alpha+\beta < \frac{5}{2}$ is still an open problem. We refer to Remark $6.11$ in Chapter $1$ of \cite{Lions} for a related discussion in the case of the equation (\ref{Fractional-NS}). 

\medskip

Leray-type solutions also verify useful energy estimates, which will be  of key importance when studying the long time behavior.  All the energy estimates that we shall perform strongly depend on the parameters $\alpha,\beta,\gamma,\delta,\nu$ in equation (\ref{Equation}). We thus set the  following notation  that we shall frequently use throughout this article:
\begin{equation}\label{Notation}
{\bf a}=\min(1,\delta^\beta), \  \  {\bf b}=\max(1,\delta^\beta)  \ \ \mbox{and} \ \  {\bf c}= m_\alpha \min(\gamma,\nu) \ \ \mbox{with} \ \  0<m_\alpha=\inf_{\xi \in \Rt} \frac{1+|\xi|^\alpha}{(1+|\xi|^2)^{\frac{\alpha}{2}}}<+\infty.
\end{equation}
At Section \ref{Sec:Preliminaries} we provide a more detailed explanation of these quantities. 
\begin{Proposition}\label{Prop1}  Within the framework of Theorem \ref{Th1}, the next energy estimates hold: 
\begin{enumerate} 
\item For $0<\gamma$ and for  all $0 \leq t$ we have:  
\[ \ds{\Vert \vu(t,\cdot) \Vert^{2}_{H^{\frac{\beta}{2}}}  \leq e^{- \gamma\,t} \left( \Vert \vu_0 \Vert^{2}_{H^{\frac{\beta}{2}}} + \frac{\b^2 }{\a^2\gamma}\,\int_{0}^{t} e^{\gamma\, s}\, \Vert \fe(s,\cdot)\Vert^{2}_{H^{\frac{\beta}{2}}} ds\right)}. \]
\item  For all $0\leq t $  and $0<T$ we have: 
\begin{equation*}
 \c \, \int_{t}^{t+T}\left\Vert  \vu(s,\cdot)\right\Vert^{2}_{H^{\frac{\alpha+\beta}{2}}}   ds  \leq  \, e^{- \gamma\,t} \left( \Vert \vu_0 \Vert^{2}_{H^{\frac{\beta}{2}}} + \frac{\b^2 }{\a^2\gamma}\,\int_{0}^{t} e^{\gamma\, s}\, \Vert \fe(s,\cdot)\Vert^{2}_{H^{\frac{\beta}{2}}} ds\right) + \frac{\b^2}{\a^2 \c }\, \int_{t}^{t+T} \Vert \fe(s,\cdot)\Vert^{2}_{H^{\frac{\beta}{2}}} ds. 
\end{equation*}
\end{enumerate}	
\end{Proposition}
The first estimate   is a direct consequence of the dissipative effects of the damping term $-\gamma \vu$. In particular, the expression $e^{-\gamma \, t}$ gives us a very good control in time on the quantity $\Vert \vu(t,\cdot) \Vert^{2}_{H^{\frac{\beta}{2}}}$, and this fact shall be well exploited in the continuing section.  
\subsection{Asymptotic behavior   of Leray-type solutions}
We consider  $\fe \in H^{\frac{\beta}{2}}(\Rt)$ (with $0\leq \beta$) a time-independent external force acting on the evolution equation (\ref{Equation}) and we shall  study  the long time behavior of Leray-type solutions.  Before to state our  results, first we need first to precise some  notation and definition. Our first definition concerns the notion of an absorbing set for the evolution equation (\ref{Equation}):
\begin{Definition}[Absorbing set]\label{def-Abosorbing-set}
A  set $\mathcal{B} \subset H^{\frac{\beta}{2}}(\Rt)$ is an absorbing set for the  equation (\ref{Equation}) if for every initial datum $\vu_0 \in H^{\frac{\beta}{2}}(\Rt)$ there exists a time $0<T=T(\vu_0)$ such that for all $T<t$ all the  Leray-type solutions $u(t,x)$ arising from $\vu_0$ verify $\vu(t,\cdot) \in \mathcal{B}$.  
\end{Definition}	

As a direct consequence of the energy estimate given in first point of  Proposition \ref{Prop1} we have the following result: 
\begin{Proposition}\label{Prop2} Let $0<\gamma$  and  $\fe \in H^{\frac{\beta}{2}}(\Rt)$.  We define 
\begin{equation}\label{Aborbing-set}
\ds{\mathcal{B}=\left\{ \vu_0 \in H^{\frac{\beta}{2}}(\Rt): \Vert \vu_0 \Vert^{2}_{H^{\frac{\beta}{2}}} \leq  \frac{2\b^2}{\a^2\gamma^2}\, \Vert \fe \Vert^{2}_{H^{\frac{\beta}{2}}}\right\}}.
\end{equation}	
Then, the set  $\mathcal{B}$ is an absorbing set for the equation (\ref{Equation}) in the sense of Definition \ref{def-Abosorbing-set}.
\end{Proposition}	

As we may observe, the absorbing set is defined by the quantities $\a$, $\b$ given in (\ref{Notation}), the  damping parameter $\gamma$ and the external force $\fe$.  Here, the expression $\ds{\frac{1}{\gamma^2}}$ clearly  shows that this definition  only makes sense when  $0<\gamma$, \emph{i.e.}, in the damped case of the equation (\ref{Equation}). Moreover, by (\ref{Notation}) we precisely have $\ds{\frac{\b^2}{\a^2}=\frac{\max(1,\delta^{2\beta})}{\min(1,\delta^{2\beta})}}$, which shows the explicit dependence of  the parameters $\beta$ and  $\delta$ in the filtering operator $(I_d - \delta^2 \Delta)^{-\frac{\beta}{2}}$.  

\medskip 

The existence of an  absorbing set for the equation (\ref{Equation}) is one of its key features in the study of the long time behavior of  Leray-type solutions. In Definition \ref{def-Abosorbing-set} we observe that all Leray-type solutions  solutions belong to the set  $\mathcal{B}$ when the time is large enough, and  consequently, their  long-time behavior  can be restricted to the set $\mathcal{B}$.

\medskip

In what follows, we explain how the absorbing set $\mathcal{B}$ is the key tool in our study.  The set $\mathcal{B}\subset H^{\frac{\beta}{2}}(\Rt)$ can be  provided of two topologies: the \emph{strong} topology  and the \emph{weak} topology  inhered from the space $H^{\frac{\beta}{2}}(\Rt)$.  Thus, when considering the  \emph{strong} topology, the absorbing set  $\mathcal{B}$   is a topological space  with the topology generated by  the usual \emph{strong} distance: 
\begin{equation}\label{strong-distance}
d_s(\vu_0, \vv_0)= \Vert \vu_0 - \vv_0 \Vert_{H^{\frac{\beta}{2}}}, \quad  \text{for all} \quad \vu_0, \vv_0 \in \mathcal{B}. 
\end{equation} 
On the other hand, the   weak  topology  on the set  $\mathcal{B}$ is generated by the  \emph{weak} distance, which is defined as follows: as $H^{\frac{\beta}{2}}(\Rt)$ is a separable Hilbert space with its usual inner product $( \cdot , \cdot )_{H^{\frac{\beta}{2}}}$,  we denote by $(\vec{e}_{n})_{n \in \mathbb{N}}$ its numerable Hilbertian basis. Then, for all $\vu_0 \in H^{\frac{\beta}{2}}(\Rt)$ we have $\ds{\vu_0 = \sum_{n\in \mathbb{N}} u_{n} \vec{e}_n}$, where  $u_{n} = (\vu_0, \vec{e}_n)_{H^{\frac{\beta}{2}}}$. Thereafter, the  \emph{weak} distance $d_w$ on $\mathcal{B}$ is given by: 
\begin{equation}\label{weak-distance}
d_w(\vu_0, \vv_0) =  \sum_{n \in \mathbb{N}} \frac{1}{2^n} \frac{\vert u_n-v_n\vert}{1+\vert u_n - v_n \vert}, \quad  \text{for all} \quad \vu_0, \vv_0 \in \mathcal{B}. 
\end{equation}

For the sake of simplicity, we shall denote the topological  metric space  $(\mathcal{B}, d_{\bullet})$, where $\bullet$ stands for either $s$ or $w$ in the case of the \emph{strong} or the \emph{weak} distances given in (\ref{strong-distance}) and (\ref{weak-distance}) respectively. In this framework, when $\bullet=s$ 
all the properties of the topological metric space $(\mathcal{B}, d_s)$ will refer as \emph{strong} properties, while  when $\bullet=w$ all the properties of the topological metric space $(\mathcal{B}, d_w)$ will refer as \emph{weak} properties.

\medskip

Our next definition is devoted to the notion of an \emph{strong} and \emph{weak} attracting set for the evolution equation (\ref{Equation}). We recall first that, in the metric space  $(\mathcal{B}, d_{\bullet})$, for $B \subset \mathcal{B}$ and $\vu_0 \in \mathcal{B}$ we define  by
\begin{equation}\label{ditance-attractor}
\ds{d_{\bullet}\big(\vu_0, B\big)= \inf_{\vv_0 \in B} d_{\bullet}(\vu_0, \vv_0)},
\end{equation}
the distance of  the point $\vu_0$ to the set $B$. 
 \begin{Definition}[Attracting set]\label{Attracting-set}  A set $B \subset \mathcal{B}$ is a $\bullet-$ attracting set for the equation (\ref{Equation}) if for all initial datum $\vu_0 \in H^{\frac{\beta}{2}}(\Rt)$ and for all  $0<\varepsilon$ there exists $0<T=T(\vu_0, \varepsilon)$ such that  all the Leray-type solutions arising from $\vu_0$ verify   $d_\bullet \big(\vu(t,\cdot), B\big) < \varepsilon$, for all $T<t$.
\end{Definition}	

Once we have the notion of the $\bullet-$attracting set, we are able to introduce now the $\bullet-$global attractor for the evolution  equation (\ref{Equation}).

\begin{Definition}[Global attractor]\label{Def-global-attractor} A set $\mathcal{A}_{\bullet} \subset \mathcal{B}$ is a $\bullet-$global attractor for the equation (\ref{Equation}) if:  
\begin{enumerate}
\item The set $\mathcal{A}_\bullet$ is $\bullet-$ compact.
\item The set  $\mathcal{A}_\bullet$ is the minimal $\bullet-$attracting set in the sense of Definition \ref{Attracting-set}
\end{enumerate}		
\end{Definition}	 

As mentioned, the notion of  the $\bullet-$global attractor is the key idea in a sharp understanding of the long time behavior of Leray-type solutions. In this definition we focus on the second point to remark that  when the time goes to infinity  the $\bullet-$ global attractor \emph{attires}  the  Leray-type solutions of the equation (\ref{Equation}). More precisely, by  Definition \ref{Attracting-set} we have that  from any initial datum $\vu_0 \in H^{\frac{\beta}{2}}(\Rt)$  all the arising Leray-type solutions are as close to $\mathcal{A}_{\bullet}$ as we want when the time is large enough. 

\medskip

The convergence of Leray-type solutions to the $\bullet-$ global attractor  is measured in terms of the distance $d_{\bullet}$ given in (\ref{ditance-attractor}). Thus,  when $\bullet=s$  these solutions converge to  the \emph{strong} global attractor $\mathcal{A}_s$   in the  \emph{strong} topology of the metric space $(\mathcal{B}, d_s)$, while when $\bullet=w$ these solutions converge to the \emph{weak} global attractor in the \emph{weak} topology of the metric space $(\mathcal{B}, d_w)$. Consequently, the \emph{strong} or the  \emph{weak} featured of the global attractor is  determined by the type of convergence and the compactness in terms of the \emph{strong}  topology or the \emph{weak} topology respectively.   

\medskip  

 When the \emph{strong} global attractor exists there also exists the \emph{weak} global attractor, and  we have the identity  $\mathcal{A}_{w}=\overline{\mathcal{A}_{s}}^{w}= \mathcal{A}_{s}$,  where $\overline{\mathcal{A}_{s}}^{w}$ denotes the cloture of $\mathcal{A}_{s}$ in the \emph{weak} topology of   the metric space $(\mathcal{B}, d_w)$.  The reverse property may not be true. We refer to \cite{CheskidovFoias} for some examples of simpler evolution equations that possess a \emph{weak} global attractor, but not a \emph{strong} global attractor.  

\medskip

Getting back to our evolution equation (\ref{Equation}), in our next result, we study the existence of a global attractor and its \emph{weak} or \emph{strong} featured.  Moreover, we shall give a characterization of the weak or strong  global attractor. For this we  recall that an eternal solution for the generalized Navier-Stokes-Bardina's model is a function  $\vu_e \in L^{\infty}_{loc}(\R, H^{\frac{\beta}{2}}(\Rt))\cap L^{2}_{loc}(\R, H^{\frac{\alpha+\beta}{2}}(\Rt))$, 
which is defined on the whole real line and it solves the equation 
\begin{equation}\label{Equation-Eternal}
 \partial_t \vu_e + \nu  (-\Delta)^{\frac{\alpha}{2}}\,  \vu_e +   (I_d -\delta^2\Delta)^{- \frac{\beta}{2}}\, \P  \left(  \text{div} (\vu_e \otimes \vu_e) \right)  = \fe -\gamma \vu_e, \qquad \text{div}(\vu_e)=0,
\end{equation}
in the weak sense. Moreover, a \emph{bounded} eternal solution is a weak solution of this equation which belongs to the space  $L^{\infty}(\R, H^{\frac{\beta}{2}}(\Rt))\cap L^{2}_{loc}(\R, H^{\frac{\alpha+\beta}{2}}(\Rt))$. 

\begin{Theoreme}\label{Th2} The following statements hold:
\begin{enumerate}
\item  	When $0<\alpha+\beta < \frac{5}{2}$,   there exists a unique \emph{weak} ($\bullet=w$) global attractor  $\mathcal{A}_{w}$ for the equation (\ref{Equation}) in the sense of Definition \ref{Def-global-attractor}.  Moreover, we have the following characterization:
\begin{equation}\label{Charac-weak}
\mathcal{A}_w=\left\{  \vu_e(0,\cdot) \in H^{\frac{\beta}{2}}(\Rt):  \,\, \vu_e \,\, \text{is an eternal solution of}\,\, (\ref{Equation-Eternal})\right\}.
\end{equation}

\item  When  $\frac{5}{2} \leq \alpha+\beta$,  there exists a unique \emph{strong} ($\bullet=s$) global attractor $\mathcal{A}_{s}$ for the equation (\ref{Equation}) in the sense of Definition \ref{Def-global-attractor}. Moreover, we have the characterization:
\begin{equation}\label{Charac-strong}
\mathcal{A}_s=\left\{  \vu_e(0,\cdot) \in H^{\frac{\beta}{2}}(\Rt):  \,\, \vu_e \,\, \text{is an \emph{bounded} eternal solution of}\,\, (\ref{Equation-Eternal})\right\}.
\end{equation}
\end{enumerate}		
\end{Theoreme}

A particular case of eternal solutions are the stationary ones belonging to the space $H^{\frac{\alpha+\beta}{2}}(\Rt)$, moreover,  by the characterizations  (\ref{Charac-weak})  and  (\ref{Charac-strong})  these solutions (when exist) belong to the global attractor $\mathcal{A}_{\bullet}$ with $\bullet=w,s$ respectively. Indeed, we just remark that all  stationary solutions $\U \in H^{\frac{\alpha+\beta}{2}}(\Rt)$ also belong to the space $L^{\infty}(\R, H^{\frac{\beta}{2}}(\Rt)) \cap L^{2}_{loc}(\R, H^{\frac{\alpha+\beta}{2}})$ (since they do not depend on the time)  and we have  $\U(0, \cdot)=\U$.

\medskip

It is thus interesting to study deeper relationships between the global attractor and stationary solutions. For this, we start by proving the existence of these latter. 
 \begin{Theoreme}\label{Th3} Let $0<\alpha$ and $0\leq \beta$. Let $\fe \in H^{\frac{\beta}{2}}(\Rt)$ be the divergence-free external force.  There exists at least  $\U \in H^{\frac{\alpha+\beta}{2}}(\Rt)$ a solution to the equation (\ref{Stationary}), which verifies the following energy  estimate  $\ds{\Vert \U \Vert_{H^{\frac{\alpha+\beta}{2}}} \leq  \frac{\b}{\a \c}\, \Vert \fe \Vert_{H^{\frac{\beta}{2}}}}$.  
\end{Theoreme}

The proof of this result is based on the Scheafer's fixed point argument, which allows us to prove the existence of solutions associated with  \emph{any} external force.  In this sense,  this is a general result for the elliptic equation (\ref{Stationary}), which is also of independent interest in the particular models (\ref{Bardina-damped}), (\ref{Critical}), (\ref{Fractional-NS}) and (\ref{Damped-NS}).  On the other hand, the uniqueness issue for stationary solutions (in the general case of any external force) seems to be more delicate. In fact, for the time-dependent equation (\ref{Equation}), uniqueness of Leray-type solutions (when $\frac{5}{2}\leq \alpha+\beta$) is obtained by energy estimates and the Gr\"onwall inequality. However, these arguments are not longer valid for equation (\ref{Stationary}) due to the lack of the temporal variable. 

\medskip

For a range of values of the parameters $\alpha$ and $\beta$, and under some \emph{sufficient conditions} depending on the external force and the parameters in equation (\ref{Equation}), we can give a more precise result on  stationary solutions. Precisely,  for a numerical constant $C>0$,   the parameters $\a$ and $\b$ defined in (\ref{Notation})  and the external force $\fe \in H^{\frac{\beta}{2}}(\Rt)$ we introduce the expression 
\[C \frac{\b}{\a^{\frac{3}{2}}} \| \fe \|_{H^{\frac{\beta}{2}}}. \]
Moreover, recall that by (\ref{Notation}) we have $\c=m_\alpha\min(\gamma,\nu)$.  In the next result, we prove that when the damping parameter $\gamma$ and the viscosity parameter $\nu$ are large enough in the sense that
\begin{equation}\label{Assumption1}
C \frac{\b}{\a^{\frac{3}{2}}} \| \fe \|_{H^{\frac{\beta}{2}}} \leq 2  \c^{\frac{3}{2}}, 
\end{equation} 
then all the stationary solutions verifying the energy estimate in Theorem \ref{Th3} are \emph{orbitally stable}, \emph{i.e.}, for all $0<\varepsilon$ we can find a quantity $0<\eta=\eta(\varepsilon)$ such that for all initial data verifying $\ds{\|\vu_0 - \U \|_{H^{\frac{\beta}{2}}} \leq \eta}$, all the arising solutions $\vu(t,\cdot)$ to equation (\ref{Equation}) satisfy $\ds{ \sup_{t>0} \| \vu(t,\cdot)-\U \|_{H^{\frac{\beta}{2}}} \leq \varepsilon}$. 

\medskip

On the other hand, when we assume the stronger control
\begin{equation}\label{Assumption2}
C \frac{\b}{\a^{\frac{3}{2}}} \| \fe \|_{H^{\frac{\beta}{2}}} \leq  \c^{\frac{3}{2}},
\end{equation}
the stationary solution $\U$ obtained in Theorem \ref{Th3} is the unique one, and it is \emph{asymptotically stable}:  for any initial data $\vu_0 \in H^{\frac{\beta}{2}}(\Rt)$ all the arising solutions $\vu(t,\cdot)$ satisfy $\ds{\lim_{t\to +\infty}\| \vu(t,\cdot)- \U \|_{H^{\frac{\beta}{2}}}=0}$. 

\begin{Theoreme}\label{Th4}  Let  $\frac{3}{2}<\alpha+\frac{\beta}{2}$  and $2 \leq \alpha + \beta$. Let $\U \in H^{\frac{\alpha+\beta}{2}}(\Rt)$ be a  solution to the stationary equation (\ref{Stationary}), which  verifies the energy  estimate  $\Vert \U \Vert_{H^{\frac{\alpha+\beta}{2}}} \leq  \frac{\b}{\a \c}\, \Vert \fe \Vert_{H^{\frac{\beta}{2}}}$.
\begin{enumerate}
	\item If (\ref{Assumption1}) holds then $\U$ is  orbitally stable. 
	
	\item If (\ref{Assumption2}) holds then   for any initial datum $\vu_0 \in H^{\frac{\beta}{2}}(\Rt)$ all the arising Leray-type solutions $\vu$ to the evolutionary equation (\ref{Equation}) (constructed in Theorem \ref{Th1}) exponentially converge to $\U$: 
	\begin{equation}\label{Inequality}
	\Vert \vu(t,\cdot) - \U \Vert^{2}_{H^{\frac{\beta}{2}}} \leq \Vert \vu_0 - \U \Vert^{2}_{H^{\frac{\beta}{2}}}\, e^{- \gamma\, t}, \quad 0 \leq t. 
	\end{equation}
	In particular, the stationary problem (\ref{Stationary}) has a unique solution verifying the energy  estimate above  and it is asymptotically stable. 
\end{enumerate}		 
\end{Theoreme}	

The proof is  essentially based on some acute energy estimates, where the conditions $\frac{3}{2}<\alpha+\frac{\beta}{2}$  and $2 \leq \alpha + \beta < \frac{5}{2}$ are required to handle the nonlinear term in equation (\ref{Equation}). Remark that these conditions are less restrictive than $\frac{5}{2}\leq \alpha+\beta$, which provides a regular enough framework to  handle  the nonlinear term in a easier way.  Moreover, these conditions are not too restrictive since they include the physically relevant models (\ref{Bardina-damped}), (\ref{Critical}), (\ref{Fractional-NS}) and (\ref{Damped-NS}). 

\medskip

On the other hand, a more interesting  featured about this result are the conditions (\ref{Assumption1}) and (\ref{Assumption2}), where  we may observe  the effects of the parameters $\gamma$ and $\nu$ in the  long time behavior of solutions to the equation (\ref{Equation}). It is worth emphasizing these kind of results on stationary solutions in the complementary case: $\gamma \lesssim  \Vert \fe \Vert_{H^{\frac{\beta}{2}}}$ or $\nu \lesssim \Vert \fe \Vert_{H^{\frac{\beta}{2}}}$,  are far from obvious and it shall be a matter of further investigations. 

\medskip

As a direct consequence of the asymptotic stability of stationary solutions we obtain the following:
\begin{Corollaire}\label{Corollary} Under the same hypothesis of the second point in Theorem \ref{Th4}, the global attractor $\mathcal{A}_{\bullet}$ obtained in Theorem \ref{Th2} verifies $\mathcal{A}_{\bullet}=\{ \U\}$. In particular, when $2\leq \alpha + \beta < \frac{5}{2}$ the weak global attractor becomes a strong one and we have $\mathcal{A}_{w}=\mathcal{A}_s=\{ \U \}$. 
\end{Corollaire}	

As mentioned, all these results hold for the particular models (\ref{Bardina-damped}), (\ref{Critical}), (\ref{Fractional-NS}) and (\ref{Damped-NS}). We summarize them in the following graphic. 

\medskip

\begin{tabular}{ll}
	\hspace{-1cm} 	
	\begin{minipage}{9cm} 
	In the region $(\alpha, \beta) \in ]0,+\infty[\times [0,+\infty[$  we graphically  represent  our results on the existence of a strong global attractor $\mathcal{A}_s$ and a weak global attractor $\mathcal{A}_w$ of the equation (\ref{Equation})  and their main related models: the Bardina's model (\ref{Bardina-damped}) is represented at the point $(\alpha, \beta)=(2,2)$, the critical Leray-alpha model (\ref{Critical}) is represented at the point $(\alpha, \beta)=(2, 1/2)$, the fractional Navier-Stokes equations (\ref{Fractional-NS}) is represented in the horizontal axis $(\alpha, 0)$;  and the classical Navier-Stokes equation (\ref{Damped-NS}) represented  at the  point $(2,0)$. The  red region  represents the conditions $\frac{3}{2}<\alpha+\beta$ and $2 \leq \alpha+\beta < \frac{5}{2}$, where the weak global attractor becomes a strong global attractor provided that (\ref{Assumption2}) holds.
	\end{minipage}\hspace{0.5cm}	
	\begin{minipage}{7.5cm}
		\begin{center}
			\includegraphics[scale=0.37]{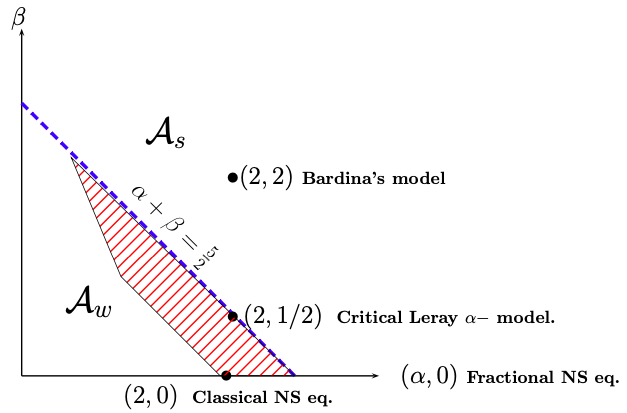} 	
		\end{center}
		{\footnotesize Figure 1: The  weak and the strong  global attractor in  the regions $0<\alpha+\beta<\frac{5}{2}$ and $\frac{5}{2}<\alpha+\beta$ respectively.}
	\end{minipage}	
\end{tabular}

\medskip

\medskip

To close this section, let mention that  it is also interesting to give a sharp estimate of the size of the global attractor.  At Appendix \ref{Appendix-Fractal-dim}, when $1\leq \alpha$ and $2 \leq \beta$  we are able to derive un upper bound for the \emph{fractal dimension} of the strong global attractor $\mathcal{A}_s$. 

\subsection{The case $\beta=0$: applications to the fractional and classical Navier-Stokes equations}\label{Sec:Applications}
It is worth making a brief discussion on  how the  results obtained for the general equation (\ref{Equation}) read down for the fractional and classical damped Navier-Stokes equations.

\medskip

First, recall that by Definition \ref{Def-Leray} for $0<\alpha$ we shall  say that $\vu$ is a Leray-type solution of the equation (\ref{Fractional-NS}) if $\vu \in L^{\infty}_{loc}\Big([0,+\infty[,L^2(\Rt)\Big) \cap L^{2}_{loc}\Big([0,+\infty[, \dot{H}^{\frac{\alpha}{2}} (\Rt)\Big)$ and it verifies  equation (\ref{Fractional-NS}) in the distributional sense. Moreover, for $\vu_0 \in L^2(\Rt)$, $\fe \in L^{2}(\Rt)$ and  for all $0\leq t$ the following energy inequality holds:	
\begin{equation*}  
\begin{split}
\Vert \vu(t,\cdot)\Vert^{2}_{L^2}  \leq \,\, \Vert \vu_0 \Vert^{2}_{L^2} -  \, 2 \nu \, \int_{0}^{t}  \left\Vert  (-\Delta)^{\frac{\alpha}{4}} \, \vu(s,\cdot)\right\Vert^{2}_{L^2}   ds + 2 \int_{0}^{t} \left( \fe, \vu(s,\cdot) \right)_{L^2} \, ds -2 \gamma \int_{0}^{t}\Vert \vu(s,\cdot)\Vert^{2}_{L^2} ds.
\end{split}
\end{equation*}

As a direct application  of Theorem \ref{Th1} (with $\beta=0$) we know that   for all $0<\alpha$  there exists $\vu$ a Leray-type solution of the equation (\ref{Fractional-NS}). Moreover, if $\frac{5}{2}\leq \alpha$ this equation has a unique Leray-type solution.

\medskip

By  the enegy estimates given in Proposition \ref{Prop1}, the equation (\ref{Fractional-NS}) has the following aborbing set given in Definition \ref{def-Abosorbing-set} (with $\beta=0$):
\begin{equation}\label{Aborbing-set-frac-NS}
\ds{\mathcal{B}_1=\left\{ \vu_0 \in L^2(\Rt): \Vert \vu_0 \Vert^{2}_{L^2} \leq  \frac{2}{\gamma^2}\, \Vert \fe \Vert^{2}_{L^2}\right\}}.
\end{equation}	
As before, the absorbing set is defined by the damping parameter $0<\gamma$ and the external force acting on the equation (\ref{Fractional-NS}), which always assumes  a time independent function.  Moreover, in this case we have $\mathcal{B}_1 \subset L^2(\Rt)$ and consequently all  the \emph{weak} or the \emph{strong} properties refer to the weak or the strong topology of this space. Precisely, when $0<\alpha < \frac{5}{2}$, there exists a unique \emph{weak} global attractor $\mathcal{A}_w \subset \mathcal{B}_1$ for the equation (\ref{Fractional-NS}), while, when $\frac{5}{2} \leq \alpha$, there exists a unique \emph{strong} global attractor $\mathcal{A}_s \subset \mathcal{B}_1$ for the equation (\ref{Fractional-NS}).

\medskip

Remark that the weak/strobg global attractor of the equation (\ref{Fractional-NS}) is a subset of the space $L^{2}(\Rt)$, while, when comparing with the general framework of the  equation (\ref{Equation}), the $\mathcal{A}_{\bullet}$ global attractor is a subset of the space $H^{\frac{\beta}{2}}(\Rt)$. It it this interesting to  remark that the regularity properties (in the framework of the Sobolev spaces) of the global attractor are only given by the parameter $0 \leq \beta$. 

\medskip

Thereafter, as a direct consequence of Theorem \ref{Th3} we known that for any (divergence-free) external force $\fe \in L^2(\Rt)$ the equation (\ref{Fractional-NS}) has at least a stationary solution $\U \in H^{\frac{\alpha}{2}}(\Rt)$, which verifies (in the weak sense) the elliptic problem:
\[ \nu(-\Delta)^{\frac{\alpha}{2}} \U + \P \, \text{div}(\U \otimes \U) = \fe - \gamma \U, \qquad \text{div}(\U)=0,\]
as well as the energy estimate $\| \U \|_{H^{\frac{\alpha}{2}}} \leq \frac{1}{\c} \| \fe \|_{L^2}$. Moreover, by Theorem \ref{Th4} we get that all these finite-energy stationary solutions are \emph{orbital stable}, provided that $C \| \fe \|_{L^2} \leq 2 c^{\frac{3}{2}}$. Moreover, when $C \| \fe \|_{L^2} \leq  c^{\frac{3}{2}}$ this stationary solution in the unique one in the energy space $H^{\frac{\alpha}{2}}(\Rt)$ and it exponentially attracts all Leray-type solutions.

\medskip

Finally, in  equation (\ref{Fractional-NS}) is it worth focusing on the particular case when $\alpha=2$, which deals with the classical damped Navier-Stokes equations (\ref{Damped-NS}).  We known that  there exists a unique \emph{weak} global attractor $\mathcal{A}_w \subset \mathcal{B}_1$ for this equation and this result can be observed as the  counterpart, in the setting of the whole space $\Rt$,  of one of the main results proven in \cite{CheskidovFoias} for the classical Navier-Stokes equations with space-periodic conditions. Indeed,  in the space-periodic setting of the torus $\mathbb{T}^3=[0,L]^3$, the absorbing set in given by 
\[ \Vert \vu_0 \Vert^2_{L^{2}(\mathbb{T^3})} \lesssim  \frac{L^2}{\nu^2} \Vert \fe \Vert^{2}_{L^{2}(\mathbb{T}^3)}, \]
see the expression (\ref{Control-Per}), and getting back to the  expression (\ref{Aborbing-set-frac-NS}) we may  observe the term $\frac{1}{\gamma}$ plays the same role of  the ratio $\frac{L}{\nu}$.  Moreover,  Theorem \ref{Th4}  could also be adapted to the space-periodic setting (with $\frac{\nu}{L}$ instead of $\gamma$)  to obtain a new result on the stability of  stationary solutions. 

 \section{Preliminaries}\label{Sec:Preliminaries} Remark that equation (\ref{Equation}) can be rewritten as 
\begin{equation}\label{Equation-P}
\partial_t \vu + \left( \gamma I_d +  \nu (-\Delta)^{\frac{\alpha}{2}} \right) \,  \vu +    (I_d -\delta^2\Delta)^{- \frac{\beta}{2}}\,  \P  \left(  \text{div} (\vu \otimes \vu) \right) + \vec{\nabla} p = \fe.
\end{equation}
For the sake of simplicity, from now on  we shall denote 
\begin{equation}\label{Notation-Operators}
J^{\alpha}_{\gamma}=  \gamma I_d +  \nu (-\Delta)^{\frac{\alpha}{2}} ,\quad J^{-\beta}_{\delta}=  (I_d -\delta^2\Delta)^{- \frac{\beta}{2}} \quad \mbox{and} \quad J^{\beta}_{\delta}= (I_d -\delta^2\Delta)^{\frac{\beta}{2}}. 
\end{equation}
All the energy estimates that we shall perform depend on the operators $J^{\alpha}_{\gamma}$ and $J^{\beta}_{\delta}$. To handle these operators we shall use the following identities 
\begin{equation}\label{Identities}
 J^{\beta}_{\delta} = D(m_1) \left(I_d -\Delta\right)^{\frac{\beta}{2}} \quad \mbox{and} \quad
J^{\alpha}_{\gamma} = D(m_2) \left(I_d -\Delta\right)^{\frac{\alpha}{2}},
\end{equation}
where the pseudo-differential operators of zero $D(m_1)$ and $D(m_2)$ are defined in the Fourier level by the symbols
\begin{equation}\label{Symbols}
m_1(\xi)= \frac{ (1+\delta^2\vert \xi \vert^2)^{\frac{\beta}{2}}}{ (1+\vert \xi \vert^2)^{\frac{\beta}{2}}}  \quad \mbox{and} \quad   m_2(\xi)= \frac{\gamma+ \nu\vert \xi \vert^\alpha}{(1+\vert \xi \vert^2)^{\frac{\alpha}{2}}}. 
\end{equation}
We see that $m_1(\xi)$ and $m_2(\xi)$ are bounded functions in $\Rt$. Precisely,  for the quantities $\a, \b, \c$ given in (\ref{Notation}) and setting  $\d= M_\alpha \max(\gamma,\nu)$ with    $\ds{M_\alpha=\sup _{\xi \in \Rt} \frac{1+|\xi|^\alpha}{(1+|\xi|^2)^{\frac{\alpha}{2}}}}$, we have the following sharp  lower and upper bounds 
\begin{equation}\label{Bounds}
\a \leq m_1(\xi) \leq \b, \quad \c \leq m_2(\xi) \leq \d, \quad \text{for all} \ \ \xi \in \Rt.
\end{equation} 
\section{Leray-type solutions}
\subsection{Proof of Theorem \ref{Th1}}\label{Sec:Leray-Solutions}
The proof  follows the  Leray's method in the classical framework of the Navier-Stokes equations. See the Section $12.1$ of the book \cite{PGLR1} for more details. For the reader's convenience, we shall detail the main estimates.  

\medskip 

 For $0<\alpha$ we denote by $p_\alpha(t,x)$ the fundamental solution of the  linear   equation $\partial_t u + J^{\alpha}_{\gamma} u = 0$, where for all $0<t$ we have $p_\alpha(t,x)= \mathcal{F}^{-1}_{x}\left( e^{- t\, ( \gamma + \nu  \vert \xi \vert^\alpha) }\right)(x)$. Here $\mathcal{F}^{-1}_{x}$ stands for the inverse Fourier transform in the spatial variable. 
 
 \medskip
 
On the other hand, let  $\theta \in \mathcal{C}^{\infty}_{0}(\Rt)$ be  a positive and radial  function  such that $\int_{\Rt} \theta(x)dx=1$. For a parameter $0<\varepsilon$ we define $\theta_\varepsilon(x)=  \frac{1}{\varepsilon^3}\,  \theta \left(\frac{x}{\varepsilon}\right)$.  In the first step, for a time $0<T<1$ small enough    we consider the following (equivalent) regularized  integral problem: 
\begin{equation}\label{Integral}
\begin{split}
\vu(t,\cdot)=  &\,  p_\alpha(t,\cdot)\ast \vu_0 - \int_{0}^{t} p_\alpha(t-s,\cdot)\ast \fe(s,\cdot) \, ds \\
&\,  - \int_{0}^{t}  p_\alpha(t-s,\cdot)\ast J^{-\beta}_{\delta} \P\left(  \theta_\varepsilon \ast  \text{div}\big( (\theta_\varepsilon \ast \vu) \otimes (\theta_\varepsilon \ast  \vu)\big) \right) (s,\cdot)\, ds.
\end{split}
\end{equation} 
This regularized  problem  is solved  in the space $L^{\infty}\Big( [0,T], H^{\frac{\beta}{2}}(\Rt)\Big)\cap L^{2}\Big([0,T],H^{\frac{\alpha+\beta}{2}}(\Rt)\Big)$, with the natural  norm $\Vert \cdot \Vert_{T}= \Vert \cdot \Vert_{L^{\infty}_{t} H^{\frac{\beta}{2}}_{x}}+ \Vert \cdot \Vert_{L^{2}_{t} H^{\frac{\alpha+\beta}{2}}_{x}}$, to obtain a unique solution $\vu_\varepsilon$.  

\medskip

In the second step, we shall prove that this solution is global in time.  The function $\vu_\varepsilon $ also solves the regularized equation:
\begin{equation}\label{Eq-Reg}
\partial_t \vu_\varepsilon + \nu (-\Delta)^{\frac{\alpha}{2}}  \vu_\varepsilon + J^{-\beta}_{\delta}  \P  \left( \theta_\varepsilon \ast \text{div}((\theta_\varepsilon \ast \vu_\varepsilon)  \otimes ( \theta_\varepsilon \ast \vu_\varepsilon) )\right)  = \fe - \gamma \vu_\varepsilon,
\end{equation}
hence, applying the operator $J^{\beta}_{\delta}$  we have
\begin{equation*}
\partial_t \, J^{\beta}_{\delta} \vu_\varepsilon + \nu (-\Delta)^{\frac{\alpha}{2}} J^{\beta}_{\delta}  \vu_\varepsilon +  \P  (\theta_\varepsilon\ast \text{div}((\theta_\varepsilon \ast \vu_\varepsilon )  \otimes  (\theta_{\varepsilon} \ast \vu_\varepsilon )))  =  J^{\beta}_{\delta} \fe- \gamma J^\beta_\delta \vu_\varepsilon.
\end{equation*}
\begin{Remarque}\label{Rmk-Energy}
From this identity we directly  obtain that  $\vu_\varepsilon$ verifies the energy estimate stated in the second point of Definition \ref{Def-Leray}.
\end{Remarque}

\medskip

Moreover, by the operator $J^{\alpha}_{\gamma}$ we can write 
\begin{equation}\label{Eq-Energ-Control} 
\partial_t \, J^{\beta}_{\delta} \vu_\varepsilon + J^{\alpha}_{\gamma} J^{\beta}_{\delta}  \vu_\varepsilon +  \P  (\theta_\varepsilon\ast \text{div}((\theta_\varepsilon \ast \vu_\varepsilon )  \otimes  (\theta_{\varepsilon} \ast \vu_\varepsilon )))  =  J^{\beta}_{\delta} \fe,
\end{equation}
and  using the identities (\ref{Identities}) we obtain 
\begin{equation*}
\partial_t \, D(m_{1}) (I_d -\Delta)^{\frac{\beta}{2}}  \vu_\varepsilon + D(m_{1}) D(m_2) (I_d-\Delta)^{\frac{\alpha+\beta}{2}}  \vu_\varepsilon + \P  (\theta_\varepsilon\ast \text{div}((\theta_\varepsilon \ast \vu_\varepsilon )  \otimes  (\theta_{\varepsilon} \ast \vu_\varepsilon )))=  D(m_{1}) (I_d -\Delta)^{\frac{\beta}{2}}  \fe.
\end{equation*}
So  we can write
\begin{equation}\label{Iden}
\frac{1}{2}\, \frac{d}{dt} \Vert D(m^{1/2}_{1}) \vu_\varepsilon(t,\cdot) \Vert^{2}_{H^{\frac{\beta}{2}}}+  \Vert  D(m^{1/2}_{1}m^{1/2}_{2}) \vu(t,\cdot) \Vert^{2}_{H^{\frac{\alpha+\beta}{2}}} = \left\langle  D(m_{1}) (I_d -\Delta)^{\frac{\beta}{2}} \fe , \vu_\varepsilon(t,\cdot) \right\rangle_{L^2 \times L^2}, 
\end{equation}
Using the lower  bounds in (\ref{Bounds}) we have 
\begin{equation*}
\frac{{\bf a}}{2 }\, \frac{d}{dt} \Vert \vu_\varepsilon(t,\cdot) \Vert^{2}_{H^{\frac{\beta}{2}}}+ {\bf a}{\bf c}\,  \Vert  \vu(t,\cdot) \Vert^{2}_{H^{\frac{\alpha+\beta}{2}}} \leq  \, \left\langle D(m^{1/2}_{1}) (I_d -\Delta)^{\frac{\beta}{4}}  \fe , D(m^{1/2}_{1}) (I_d -\Delta)^{\frac{\beta}{4}} \vu_\varepsilon(t,\cdot) \right\rangle_{L^2 \times L^2}=I_3. 
\end{equation*}
Moreover, using the upper bounds in (\ref{Bounds}), the last term can be estimated as: 
\begin{equation*}
I_3 \leq   {\bf b} \, \Vert \fe (t,\cdot)  \Vert_{H^\frac{\beta}{2}} \Vert \vu_\varepsilon(t,\cdot) \Vert_{H^{\frac{\beta}{2}}}  \leq \,  \frac{{\bf b}^2}{2 {\bf a}{\bf c}} \Vert \fe(t,\cdot) \Vert^{2}_{H^\frac{\beta}{2}} + \frac{{\bf a}{\bf c}}{2} \Vert \vu_\varepsilon(t,\cdot) \Vert^{2}_{H^{\frac{\beta}{2}}} \leq  \frac{{\bf b}^2}{2 {\bf a}{\bf c}}\Vert \fe(t,\cdot) \Vert^{2}_{H^\frac{\beta}{2}} + \frac{{\bf a}{\bf c}}{2} \Vert \vu_\varepsilon(t,\cdot) \Vert^{2}_{H^{\frac{\alpha+\beta}{2}}}. 
\end{equation*}
We thus obtain the inequality:
\begin{equation}\label{Energy-estimate-1} 
\frac{{\bf a}}{2}\, \frac{d}{dt} \Vert \vu_\varepsilon(t,\cdot) \Vert^{2}_{H^{\frac{\beta}{2}}}+ \frac{{\bf a}{\bf c}}{2} \,   \Vert  \vu(t,\cdot) \Vert^{2}_{H^{\frac{\alpha+\beta}{2}}}  \leq  \frac{{\bf b}^2}{2 {\bf a}{\bf c}}\Vert \fe(t,\cdot) \Vert^{2}_{H^\frac{\beta}{2}},
\end{equation}
and we integrate on the interval of time $[0,t]$ to get the following control:
\begin{equation}\label{Energy-estimate-2}
\Vert \vu_\varepsilon(t,\cdot) \Vert^{2}_{H^{\frac{\beta}{2}}} + {\bf c}\, \int_{0}^{t} \Vert   \vu_\varepsilon(s,\cdot) \Vert^{2}_{H^{\frac{\alpha+\beta}{2}}}\, ds \leq \, \Vert \vu_0 \Vert^{2}_{H^{\frac{\beta}{2}}} + \frac{{\bf b}^2}{{\bf a}^2{\bf c}}\, \int_{0}^{t} \Vert \fe(s,\cdot) \Vert^{2}_{H^{\frac{\beta}{2}}} ds, 
\end{equation}
which allows us to extend the local solution $\vu_\varepsilon$ to the whole interval of time $[0,+\infty[$.  

\medskip

In the third step,  we study  the convergence to a weak solution of  equation (\ref{Equation-P}). By the Rellich-Lions lemma (see \cite{PGLR1}, Theorem $12.1$) there exists a sequence of positive numbers $(\varepsilon_n)_{n \in \mathbb{N}}$ and a function $\vu \in L^{2}_{loc}([0,+\infty[\times \Rt)$ such that the sequence $\vu_{\varepsilon_n}$ converges to $\vu$ in the strong topology of the space  $(L^{2}_{t,x})_{loc}$. Moreover, this sequence also converges to $\vu$ in the weak$-*$ topology of the spaces $L^{\infty}([0,T], H^{\frac{\beta}{2}}(\Rt))$ and $L^{2}([0,T], H^{\frac{\alpha+\beta}{2}}(\Rt))$ for all $0<T$.  We must study the convergence to the nonlinear term $\P\,\text{div}(\vu \otimes \vu)$;  and for this we shall need the following: 
\begin{Lemme}\label{Lemma-uniq} Let $\vu \in  (L^{\infty}_{t})_{loc} H^{\frac{\beta}{2}}_{x} \cap (L^{2}_{t})_{loc} H^{\frac{\alpha+\beta}{2}}_{x}$. If  $ 0 <\alpha + \beta <3$  then $ \text{div}( \vu \otimes \vu )\in (L^{2}_{t})_{loc} H^{\frac{-5+(\alpha +\beta)}{2}}_{x}$. 
\end{Lemme}	
\pv    As we have $\frac{\alpha + \beta}{2} < \frac{3}{2}$ then by the Hardy-Littlewood-Sobolev (HLS) inequalities we obtain  $\vu  \in (L^{2}_{t})_{loc} H^{\frac{\alpha+\beta}{2}}_{x}  \subset (L^{2}_{t})_{loc} \dot{H}^{\frac{\alpha+\beta}{2}}_{x} \subset  (L^{2}_{t})_{loc} L^{p}_{x}$, with $p=\frac{6}{3-(\alpha+\beta)}$. On the other hand, as $0 \leq \beta$ we have $\vu \in (L^{\infty}_{t})_{loc} H^{\frac{\beta}{2}}_{x} \subset (L^{\infty}_{t})_{loc} L^{2}_{x}$.  Thus, we can use the H\"older inequalities with $\frac{1}{q}= \frac{1}{2}+\frac{1}{p}$, hence we have $q=\frac{6}{6-(\alpha+\beta)}$, and we are able to write $\Vert \vu \otimes \vu \Vert_{(L^{2}_{t})_{loc} L^{q}_{x}} \leq C\,  \Vert \vw \Vert_{(L^{\infty}_{t})_{loc} L^{2}_{x}}\, \Vert \vu_1 \Vert_{(L^{2}_{t})_{loc} L^{p}_{x}}$.  Finally, by making use again of the (HLS) inequalities we have the embedding $L^{q}(\Rt) \subset \dot{H}^{\frac{-3+(\alpha+\beta)}{2}}(\Rt)$, hence we obtain $\vu \otimes \vu \in (L^{2}_{t})_{loc} \dot{H}^{\frac{-3+(\alpha+\beta)}{2}}_{x}$ and consequently we have $\text{div}(\vu \otimes \vu) \in  (L^{2}_{t})_{loc} \dot{H}^{\frac{-5+(\alpha+\beta)}{2}}_{x} \subset  (L^{2}_{t})_{loc} H^{\frac{-5+(\alpha+\beta)}{2}}_{x}$.  \finpv

\bigskip 

In the case when $0<\alpha + \beta <3$, by this lemma and by (\ref{Energy-estimate-2})  the family $(\theta_{\varepsilon} \ast \text{div} ((\theta_{\varepsilon} \ast \vu_\varepsilon) \otimes (\theta_\varepsilon \ast \vu_\varepsilon))$  is uniformly bounded respect to  the parameter $\varepsilon$ in the space $(L^{2}_{t})_{loc} H^{\frac{-5+ (\alpha+\beta)}{2} }_{x}$. On the other hand, in the case  when  $3 \leq \alpha +\beta$  we    set $0<\alpha'<\alpha$ and $0\leq \beta' <\beta$ such that $0<\alpha' + \beta'<3$  to obtain that the  family $(\theta_{\varepsilon} \ast \text{div} ((\theta_{\varepsilon} \ast \vu_\varepsilon) \otimes (\theta_\varepsilon \ast \vu_\varepsilon))$  is uniformly bounded  in  $(L^{2}_{t})_{loc} H^{\frac{-5+ (\alpha'+\beta')}{2} }_{x}$.  Consequently, in both cases  we   obtain the uniformly boundness of the family above in the larger space $(L^{2}_{t})_{loc} H^{-\frac{5}{2} }_{x}$.  From this fact and the convergences above we can deduce that the sequence $\left( \P \text{div} \left( \theta_{\varepsilon_n}\ast \vu_{\varepsilon_n} \otimes \vu_{\varepsilon_n}\right) \right)_{n\in \mathbb{N}}$ converges to $\P \, \text{div}\left( \vu \otimes \vu \right)$ in the weak$-*$ topology of the space $(L^{2}_{t})_{loc}H^{-\frac{5}{2}}_{x}$. 

\medskip

In the fourth step, to obtain the energy inequality given at the second point of Definition \ref{Def-Leray},  we get back to Remark \ref{Rmk-Energy} and applying classical tools (see the page $354$ of the  book\cite{PGLR1}) this inequality also holds  for the limit $\vu$.   

\subsection*{Uniqueness in the case $\frac{5}{2}\leq \alpha+\beta$.}  
 For simplicity, we shall omit the constants as well as the operators $D(m_1)$ and $D(m_2)$ given in (\ref{Identities}).  

\medskip

Let $\vu_1, \vu_2 \in (L^{\infty}_{t})_{loc} H^{\frac{\beta}{2}}_{x} \cap (L^{2}_{t})_{loc} H^{\frac{\alpha+\beta}{2}}_{x}$ be two Leray-type  solutions of  equation (\ref{Equation}) with the same external force $\fe \in (L^{2}_{t})_{loc} H^{\frac{\beta}{2}}_{x}$ and arising from the initial data $\vu_{0,1}$ and $\vu_{0,2}$ respectively. We define $\vw=\vu_1-\vu_2$, which,  up to the operators $D(m_1)$ and $D(m_2)$, it essentially solves the following problem: 
\begin{equation}\label{Eq-w}
\begin{cases}
\partial_t (I_d -\Delta)^{\frac{\beta}{2}}  \vw + (I_d -\Delta)^{\frac{\alpha+\beta}{2}}  \vw +  \P \left( (\vw \cdot \vec{\nabla}) \vu_1+ (\vu_2 \cdot \vec{\nabla}) \vw\right)=0,\\
\text{div}(\vw)=0,\\ 
\vw(0,\cdot)=\vw_0 = \vu_{0,1}-\vu_{0,2}. 
\end{cases}
\end{equation}
We shall perform an energy estimate on the solution $\vw$. First, we consider the case when $0<\alpha+\beta <3$. In this case, by Lemma \ref{Lemma-uniq} we have $\P \left( (\vw \cdot \vec{\nabla}) \vu_1 \right) \in (L^{2}_{t})_{loc} H^{\frac{5-(\alpha+\beta)}{2}}_{x}$ and $\P \left( (\vu_2 \cdot \vec{\nabla}) \vw \right) \in  (L^{2}_{t})_{loc} H^{\frac{5-(\alpha+\beta)}{2}}_{x}$.  Moreover, by assuming $\frac{5}{2}\leq \alpha+ \beta<3$ then we get $\frac{5-(\alpha+\beta)}{2} \leq \frac{\alpha + \beta }{2}$; so we have  $\vw \in (L^{2}_{t})_{loc} H^{\frac{5-(\alpha+\beta)}{2}}_{x}$. Then, we can write 
\begin{equation}\label{Identity}
\begin{split}
\frac{1}{2} \frac{d}{dt} \Vert \vw(t,\cdot) \Vert^{2}_{H^{\frac{\beta}{2}}} + & \Vert \vw(t,\cdot) \Vert^{2}_{H^{\frac{\alpha+\beta}{2}}}+ \left\langle \P \left( (\vw \cdot \vec{\nabla}) \vu_1 \right) , \vw \right\rangle_{H^{\frac{-5+(\alpha+\beta)}{2}} \times H^{\frac{5-(\alpha+\beta)}{2}}} \\
+ & \left\langle  \P \left( (\vu_2 \cdot \vec{\nabla}) \vw \right) , \vw \right\rangle_{H^{\frac{-5+(\alpha+\beta)}{2}} \times H^{\frac{5-(\alpha+\beta)}{2}}} =0. 
\end{split}
\end{equation} 
As $\text{div}(\vw)=0$, we have $ \left\langle  \P \left( (\vu_2 \cdot \vec{\nabla}) \vw \right) , \vw \right\rangle_{H^{\frac{-5+(\alpha+\beta)}{2}} \times H^{\frac{5-(\alpha+\beta)}{2}}}=0$, so it remains to estimate the term $\left\langle \P \left( (\vw \cdot \vec{\nabla}) \vu_1 \right) , \vw \right\rangle_{H^{\frac{-5+(\alpha+\beta)}{2}} \times H^{\frac{5-(\alpha+\beta)}{2}}} $. More precisely, the following estimate holds:
\begin{equation*}
\left\vert  \left\langle \P \left( (\vw \cdot \vec{\nabla}) \vu_1 \right) , \vw \right\rangle_{H^{\frac{-5+(\alpha+\beta)}{2}} \times H^{\frac{5-(\alpha+\beta)}{2}}} \right\vert  \lesssim  \,  \Vert \vw(t,\cdot) \Vert_{H^{\frac{\beta}{2}}}\, \Vert \vu_1(t,\cdot) \Vert_{H^{\frac{\alpha+\beta}{2}}}\, \Vert \vw(t,\cdot) \Vert_{H^{\frac{\alpha+\beta}{2}}}.
\end{equation*}

With this estimate,  we get back to the last identity to write  
 \begin{equation*} 
\frac{1}{2} \frac{d}{dt} \Vert \vw(t,\cdot)\Vert^{2}_{H^{\frac{\beta}{2}}} + \Vert \vw(t,\cdot) \Vert^{2}_{H^{\frac{\alpha+\beta}{2}}} \lesssim \,  \Vert \vw(t,\cdot) \Vert^{2}_{H^{\frac{\beta}{2}}}\, \Vert \vu_1(t,\cdot) \Vert^{2}_{H^{\frac{\alpha+\beta}{2}}} +  \Vert \vw(t,\cdot) \Vert^{2}_{H^{\frac{\alpha+\beta}{2}}}, 
\end{equation*}
hence we obtain:
\begin{equation*}
\frac{d}{dt} \Vert \vw(t,\cdot)\Vert^{2}_{H^{\frac{\beta}{2}}} \lesssim   \,  \Vert \vw(t,\cdot) \Vert^{2}_{H^{\frac{\beta}{2}}}\, \Vert \vu_1(t,\cdot) \Vert^{2}_{H^{\frac{\alpha+\beta}{2}}}.
\end{equation*}
From this estimate and the Gr\"onwall inequalities we have 
\begin{equation*}
\Vert \vw(t,\cdot) \Vert^{2}_{H^{\frac{\beta}{2}}} \lesssim \Vert \vw_0 \Vert^{2}_{H^{\frac{\beta}{2}}} \, \text{exp} \left( \ds{\int_{0}^{t} \Vert \vu_1(s,\cdot) \Vert^{2}_{H^{\frac{\alpha+\beta}{2}}} ds} \right). 
\end{equation*}
Moreover, as $\vu_1$ verifies the energy inequality stated in Definition \ref{Def-Leray}  then  we can  write: 
\begin{equation*}
\Vert \vu_1(t,\cdot)\Vert^{2}_{H^{\frac{\beta}{2}}} +  \int_{0}^{t} \Vert \vu_1(s,\cdot)\Vert^{2}_{H^{\frac{\alpha+\beta}{2}}}  ds \lesssim  \Vert \vu_{0,1} \Vert^{2}_{H^{\frac{\beta}{2}}} + \frac{1}{2} \int_{0}^{t} \Vert \fe(s,\cdot) \Vert^{2}_{H^{\frac{\beta}{2}}} ds +  \frac{1}{2} \int_{0}^{t} \Vert \vu_1(s,\cdot) \Vert^{2}_{H^{\frac{\alpha+ \beta}{2}}} ds,
\end{equation*}
and  we get 
\begin{equation*}
\Vert \vu_1(t,\cdot)\Vert^{2}_{H^{\frac{\beta}{2}}} + \int_{0}^{t} \Vert \vu_1(s,\cdot)\Vert^{2}_{H^{\frac{\alpha+\beta}{2}}}  ds \lesssim  \Vert \vu_{0,1} \Vert^{2}_{H^{\frac{\beta}{2}}} + \int_{0}^{t} \Vert \fe(s,\cdot) \Vert^{2}_{H^{\frac{\beta}{2}}} ds. 
\end{equation*}
Consequently we have  
\begin{equation}\label{estim-aux}
\int_{0}^{t} \Vert \vu(s,\cdot)\Vert^{2}_{H^{\frac{\alpha+\beta}{2}}}  ds \lesssim \Vert \vu_{0,1} \Vert^{2}_{H^{\frac{\beta}{2}}} + \int_{0}^{t} \Vert \fe(s,\cdot) \Vert^{2}_{H^{\frac{\beta}{2}}} ds.
\end{equation}
Then, we can write:
\begin{equation}\label{Uniqueness1}
\Vert \vw(t,\cdot) \Vert^{2}_{H^{\frac{\beta}{2}}} \lesssim \Vert \vw_0 \Vert^{2}_{H^{\frac{\beta}{2}}} \, \text{exp} \left(  \Vert \vu_{0,1} \Vert^{2}_{H^{\frac{\beta}{2}}}+\int_{0}^{t} \Vert \fe(s,\cdot) \Vert^{2}_{H^{\frac{\beta}{2}}} ds \right).
\end{equation}

\medskip 

On the other hand,  when $ 3 \leq \alpha +\beta$, we can always  set $0 \leq \beta' < \beta$ and $0< \alpha' < \alpha$ such that $\frac{5}{2} \leq \alpha' + \beta' <3$. Summarizing, for $\frac{5}{2}\leq \alpha +\beta$ we have the inequality  (\ref{Uniqueness1}) from which the uniqueness of Leray-type solutions directly follows. Theorem \ref{Th1} is  proven.  \finpv

\subsection{Energy estimates: proof of  Proposition \ref{Prop1}}

We consider the functions $\vu_{\varepsilon_n}$, which are solutions of  the regularized equation (\ref{Eq-Reg}).  Then we have the identity
\begin{equation}\label{iden-reg}
\begin{split}
& \frac{1}{2} \frac{d}{dt} \Vert D(m^{1/2}_{1}) \vu_{\varepsilon_n}(t,\cdot)\Vert^{2}_{H^{\frac{\beta}{2}}}+ \gamma\Vert D(m^{1/2}_{1}) \vu_{\varepsilon_n}(t,\cdot)\Vert^{2}_{H^{\frac{\beta}{2}}}+  \nu\left\Vert D(m^{1/2}_{1})) (-\Delta)^{\frac{\alpha}{4}} (I_d-\Delta)^{\frac{\beta}{4}} \vu_{\varepsilon_n}(t,\cdot) \right\Vert^{2}_{L^2} \\
=& \,  \left\langle D(m^{1/2}_{1}) (I_d-\Delta)^{\frac{\beta}{4}}\fe(t,\cdot), D(m^{1/2}_{1})(I_d-\Delta)^{\frac{\beta}{4}}\vu_{\varepsilon_n}(t,\cdot) \right\rangle_{L^2\times L^2},
\end{split}
\end{equation}
hence, using the lower and upper bounds in (\ref{Bounds}) we can write
\begin{equation*}
\frac{\a}{2} \frac{d}{dt} \Vert \vu_{\varepsilon_n}(t,\cdot)\Vert^{2}_{H^{\frac{\beta}{2}}} +  \gamma \a\Vert \vu_{\varepsilon_n}(t,\cdot)\Vert^{2}_{H^{\frac{\beta}{2}}} \leq  \, \b \Vert \fe \Vert_{H^{\frac{\beta}{2}}}\, \Vert \vu_{\varepsilon_n}(t,\cdot)\Vert_{H^{\frac{\beta}{2}}} \leq \frac{\b^2}{2\a\gamma }\Vert \fe(t,\cdot) \Vert^{2}_{H^{\frac{\beta}{2}}}+ \frac{\gamma\a}{2} \Vert \vu_{\varepsilon_n}(t,\cdot)\Vert^{2}_{H^\frac{\beta}{2}}, 
\end{equation*} 
and then we obtain 
\[ \frac{d}{dt}\Vert \vu_{\varepsilon_n}(t,\cdot)\Vert^{2}_{H^{\frac{\beta}{2}}} + \gamma  \Vert \vu_{\varepsilon_n}(t,\cdot)\Vert^{2}_{H^{\frac{\beta}{2}}} \leq \frac{\b^2}{\a^2\gamma } \Vert \fe(t,\cdot) \Vert^{2}_{H^{\frac{\beta}{2}}}. \]
Thereafter, by applying the Gr\"onwall inequalities we have 
\begin{equation}\label{Ineq}
\Vert \vu_{\varepsilon_n}(t,\cdot) \Vert^{2}_{H^{\frac{\beta}{2}}} \leq e^{- \gamma \, t} \, \Vert \vu_0\Vert^{2}_{H^{\frac{\beta}{2}}}  +  \frac{\b^2 e^{-\gamma t} }{\a^2\gamma}\, \int_{0}^{t} e^{\gamma \,s}\, \Vert \fe(s,\cdot)\Vert^{2}_{H^{\frac{\beta}{2}}} ds.  
\end{equation}
We will recover this control in time for the limit function $\vu$: we regularize in the time variable the quantity $\Vert \vu_{\varepsilon_n}(t,\cdot) \Vert^{2}_{H^{\frac{\beta}{2}}}$ by a convolution product with a positive function $w \in \mathcal{C}^{\infty}_{0}([-\eta, \eta])$ (with $0<\eta$) such that $\ds{\int_{\R} w(t) dt =1}$. Thus, in the previous estimate we have
\[  \Vert  w \ast \vu_{\varepsilon_n}(t,\cdot) \Vert^{2}_{H^{\frac{\beta}{2}}}  \leq w \ast \Vert \vu_{\varepsilon_n}(t,\cdot) \Vert^{2}_{H^{\frac{\beta}{2}}} \leq w \ast \left( e^{- \gamma \, t} \, \Vert \vu_0\Vert^{2}_{H^{\frac{\beta}{2}}}  + \frac{\b^2 e^{-\gamma t} }{\a^2\gamma}\, \int_{0}^{t} e^{\gamma \,s}\, \Vert \fe(s,\cdot)\Vert^{2}_{H^{\frac{\beta}{2}}} ds \right).\] 
As $(\vu_{\varepsilon_n})_{n\in \mathbb{N}}$ converges weakly$-*$ to $\vu$ in the space $(L^{\infty}_{t})_{loc}H^{\frac{\beta}{2}}_{x}$ then $ w \ast \vu_{\varepsilon_n}(t,\cdot)$ converges weakly$-*$ to $\vu(t,\cdot)$ in the space $H^{\frac{\beta}{2}}(\Rt)$ and we are able to write:
\[  \Vert  w \ast \vu (t,\cdot) \Vert^{2}_{H^{\frac{\beta}{2}}}  \leq \liminf_{n\to + \infty} \Vert  w \ast \vu_{\varepsilon_n}(t,\cdot) \Vert^{2}_{H^{\frac{\beta}{2}}}  \leq  w\ast \left(  e^{- \gamma \, t} \, \Vert \vu_0\Vert^{2}_{H^{\frac{\beta}{2}}}  + \frac{\b^2 e^{-\gamma t} }{\a^2\gamma} \, \int_{0}^{t} e^{\gamma \,s}\, \Vert \fe(s,\cdot)\Vert^{2}_{H^{\frac{\beta}{2}}} ds\right).\]
In this fashion, for $0\leq t$ a Lebesgue point of the function $t \mapsto \Vert \vu(t,\cdot)\Vert^{2}_{H^{\frac{\beta}{2}}}$ we have the energy control stated in the first point of Proposition \ref{Prop1}. Moreover, this energy control is extended to all time $0 \leq t $ by the weak continuity of the function  $t \mapsto \Vert \vu(t,\cdot)\Vert^{2}_{H^{\frac{\beta}{2}}}$. 

\medskip 

In order to prove the second point of  Proposition \ref{Prop1} we get back to the inequality (\ref{Energy-estimate-1}).  Integrating in the interval of time $[t, t+T]$ we have: 
\[ \Vert \vu_{\varepsilon_n}(t+T,\cdot) \Vert^{2}_{H^{\frac{\beta}{2}}} +\c \,  \int_{t}^{t+T} \Vert \vu_{\varepsilon_n}(s,\cdot) \Vert^{2}_{H^{\frac{\alpha+\beta}{2}}} ds \leq \Vert \vu_{\varepsilon_n}(t,\cdot) \Vert^{2}_{H^{\frac{\beta}{2}}} + \frac{\b2}{\a^2 \c }\,\int_{t}^{t+T}\Vert \fe(s,\cdot) \Vert^{2}_{H^{\frac{\beta}{2}}} ds. \]
Hence, by the estimate (\ref{Ineq}) we can write:
\[ \c  \int_{t}^{t+T} \Vert \vu_{\varepsilon_n}(s,\cdot) \Vert^{2}_{H^{\frac{\alpha+\beta}{2}}} ds \leq  e^{- \gamma\, t} \, \Vert \vu_0\Vert^{2}_{H^{\frac{\beta}{2}}}  + \frac{\b^2 e^{-\gamma t}}{\a^2 \gamma} \, \int_{0}^{t} e^{\gamma\,s}\, \Vert \fe(s,\cdot)\Vert^{2}_{H^{\frac{\beta}{2}}} ds + \frac{\b^2}{\a^2 \c }\,\int_{t}^{t+T}\Vert \fe(s,\cdot) \Vert^{2}_{H^{\frac{\beta}{2}}} ds.\]
By recalling that  $(\vu_{\varepsilon_n})_{n\in \mathbb{N}}$ converges weakly$-*$ to $\vu$ in the space $(L^{2}_{t})_{loc}H^{\frac{\alpha+\beta}{2}}_{x}$ we obtain the desired estimate. 

\medskip

  Proposition \ref{Prop1} is proven. \finpv

\section{Long time behavior of Leray-type  solutions}\label{Sec:Long-time-behavior} 
As explained in Section \ref{Sec:Results}, the notion of absorbing set is of key importance when studying  the existence of  global attractors either  the weak and the strong case.  Our starting point is then to verify that the set $\mathcal{B}$ given in (\ref{Aborbing-set}) is an absorbing set for the  equation (\ref{Equation}).
\subsection{The absorbing set: proof of Proposition \ref{Prop2}}
Let $\vu_0 \in  H^{\frac{\beta}{2}}(\Rt)$ be an initial datum, and let $\vu(t,\cdot)$ be a Leray-type solution of the equation (\ref{Equation}) arising from $\vu_0$.  By the first point of Proposition \ref{Prop1} and as $\fe \in H^{\frac{\beta}{2}}(\Rt)$  is a time-independent function  we can write:
\begin{equation}\label{Estim-absorbing}
 \Vert \vu(t,\cdot) \Vert^{2}_{H^{\frac{\beta}{2}}} \leq e^{- \gamma\, t} \Vert \vu_0 \Vert^{2}_{H^{\frac{\beta}{2}}} + \frac{\b^2}{\a^2\gamma^2}\,  \Vert \fe \Vert^{2}_{H^{\frac{\beta}{2}}}  (1-e^{- \gamma\,t}) \leq  e^{- \gamma\, t} \Vert \vu_0 \Vert^{2}_{H^{\frac{\beta}{2}}} + \frac{\b^2}{\a^2\gamma^2}\,  \Vert \fe \Vert^{2}_{H^{\frac{\beta}{2}}}.
\end{equation} 
Hence, we can set a time $0<T=T\left(\a,\b,\gamma,\Vert \vu_0 \Vert^2_{H^{\frac{\beta}{2}}}, \Vert \fe \Vert^{2}_{H^{\frac{\beta}{2}}}\right)$ large enough such that for all $T<t$ we have the inequality  $ e^{- \gamma\, t} \Vert \vu_0 \Vert^{2}_{H^{\frac{\beta}{2}}} \leq  \frac{\b^2}{\a^2\gamma^2}\,  \Vert \fe \Vert^{2}_{H^{\frac{\beta}{2}}}$. Consequently, for all $T<t$ we have  $\vu(t,\cdot) \in \mathcal{B}$. Proposition \ref{Prop2} is proven. \finpv

\subsection{Global attractor: proof of Theorem \ref{Th3}}
\subsubsection{The weak  global attractor in the case $0<\alpha+\beta<\frac{5}{2}$.}
The existence and uniqueness of a weak global attractor $\mathcal{A}_w$ for the equation (\ref{Equation}) bases on the following previous results that we summarize as follows.   First, for the absorbing set  $\mathcal{B}$ given in (\ref{Aborbing-set}) and for a time $0\leq t$ we define the set 
\begin{equation}\label{family-R}
R(t)\mathcal{B}= \big\{  \vu(t,\cdot) : \,\, \text{$\vu$ is a Leray-type solution of (\ref{Equation}) arising from}\,\, \vu_0 \in \mathcal{B}\big\} \subset H^{\frac{\beta}{2}}(\Rt).
\end{equation}
As uniqueness of Leray-type solutions is unknown for this range of values of the parameter $\alpha+\beta$,  the family  $(R(t))_{t \geq 0}$ does not define a semigroup on the space $H^{\frac{\beta}{2}}(\Rt)$. However, this family enjoys the following property: $R(t_1+t_2)\mathcal{B} \subset R(t_1) R(t_2)\mathcal{B}$, for all $0\leq t_1,t_2$. We introduce now the following:
\begin{Definition}[Weakly  uniformly compact family]\label{Def-Uniform-weak-compact} The family $(R(t))_{t\geq 0}$ is uniformly weakly compact if there exists a time $0<T$ such that set $\ds{\bigcup_{T\leq t} R(t)\mathcal{B}}$ is relatively compact in $(\mathcal{B}, d_w )$, where the distance $d_w$ is given in (\ref{weak-distance}).
\end{Definition}

Now we can state the  following result on the existence of a weak global attractor. For a proof see  \cite{CheskidovFoias}, Theorem $2.11$ and Corollary $2.5$. 
\begin{Theoreme}[Existence of a weak global attractor]\label{Th-existence-weak}  If the  family 	$(R(t))_{t\geq 0}$ given in (\ref{family-R}) is  uniformly weak compact in the sense of Definition \ref{Def-Uniform-weak-compact}, then  there exists a unique weak global attractor $\mathcal{A}_w$ in the sense of  Definition \ref{Def-global-attractor}.  
\end{Theoreme}	

\noindent
{\bf Proof of the first point in Theorem \ref{Th3}}.  By Theorem \ref{Th-existence-weak}, we shall prove that the family $(R(t))_{t\geq 0}$  is  uniformly weak compact.   Let $0\leq t$ and let $\vu(t,\cdot) \in R(t)\mathcal{B}$. By definition of the set $R(t)\mathcal{B}$ given in (\ref{family-R}) we known that $\vu(t,\cdot)$ is a Leray-type solution of the equation (\ref{Equation}) arising from an initial datum $\vu_0 \in \mathcal{B}$. Then,  by  definition of the absorbing set $\mathcal{B}$ given in (\ref{Aborbing-set}), and moreover, by the estimate (\ref{Estim-absorbing}) we have 
\begin{equation}\label{Estim-Arb}
 \Vert \vu(t,\cdot) \Vert^{2}_{H^{\frac{\beta}{2}}} \leq e^{- \gamma\, t} \frac{2\b^2}{\a^2\gamma^2} \| \fe \|^{2}_{H^{\frac{\beta}{2}}} + \frac{\b^2}{\a^2\gamma^2}\,  \Vert \fe \Vert^{2}_{H^{\frac{\beta}{2}}}.
\end{equation}
We can set a time $0<T$, which does not depend on $\vu_0 \in \mathcal{B}$, such that for all $T<t$ we have $2\, e^{- \gamma\, t}  \leq  1$.  Thus, for all $T<t$ we have $\Vert \vu(t,\cdot)\Vert^{2}_{H^{\frac{\beta}{1}}} \leq \frac{2\b^2}{\a^2\gamma^2}\Vert \fe \Vert^{2}_{H^{\frac{\beta}{2}}}$, hence, always by definition of the set $\mathcal{B}$, we obtain  $\ds{\bigcup_{T\leq t} R(t)\mathcal{B} \subset \mathcal{B}}$. Finally, as $(\mathcal{B}, d_w)$ is a compact metric space the family $(R(t))_{t \geq 0}$ is then uniformly weak compact in the sense of Definition \ref{Def-Uniform-weak-compact}.  Thus, by Theorem \ref{Th-existence-weak} there exists a unique weak global attractor $\mathcal{A}_w$.  

\medskip

We prove now the characterization of  $\mathcal{A}_w$ given in (\ref{Charac-weak}). For this,  we  shall need to introduce some  notation.  We denote by  $\mathcal{L}_{\mathcal{B}}([0,+\infty])$ the set of all the Leray-type solutions of the equation (\ref{Equation}) arising from initial data in $\mathcal{B}$.  Moreover, we denote  by $H^{\frac{\beta}{2}}_{w} (\Rt)$ the space $H^{\frac{\beta}{2}}(\Rt)$  endowed with its weak topology. Then, we  denote by  $\mathcal{C}_{w}\Big([0,+\infty[, H^{\frac{\beta}{2}}_{w}(\Rt)\Big)$ the  space  of weak$-*$ continuous ${H^{\frac{\beta}{2}}}-$ valued functions on  the interval $[0,+\infty[$.  We thus have the embedding $\mathcal{L}_{\mathcal{B}}([0,+\infty]) \subset \mathcal{C}_{w}\Big([0,+\infty[, H^{\frac{\beta}{2}}_{w}(\Rt)\Big)$.

\medskip 

We shall prove that $\mathcal{L}_{\mathcal{B}}([0,+\infty])$ is compact in $\mathcal{C}_{w}\Big([0,+\infty[, H^{\frac{\beta}{2}}_{w}(\Rt)\Big)$. Let $(\vu_n)_{n\in \mathbb{N}}$ be a sequence in $\mathcal{L}_{\mathcal{B}}([0,+\infty[)$. Indeed, we just remark that by the estimate (\ref{Estim-Arb}) we have the following uniform bound: $ \Vert \vu_n(t,\cdot) \Vert^{2}_{H^{\frac{\beta}{2}}} \leq  \frac{5}{\gamma^2} + \frac{4}{\gamma^2}\,  \Vert \fe \Vert^{2}_{H^{\frac{\beta}{2}}}$, hence,  the sequence $(\vu_n)_{n\in \mathbb{N}}$ has a subsequence which converges in the space $\mathcal{C}_{w}\Big([0,+\infty[, H^{\frac{\beta}{2}}_{w}(\Rt)\Big)$.  Then, by Theorem $2.14$ of \cite{CheskidovFoias} we have (\ref{Charac-weak}). The first point of Theorem \ref{Th3}  is now  proven. \finpv

\subsubsection{The strong  global  attractor in the case $\frac{5}{2} \leq \alpha+\beta$}\label{Sec:strong-global-attractor} 
Uniqueness of Leray-type  solutions allows us to define the semigroup $S(t):H^{\frac{\beta}{2}}(\Rt) \to H^{\frac{\beta}{2}}(\Rt)$ as: 
\begin{equation}\label{Seligroup}
S(t)\vu_0= \vu(t,\cdot), \quad 0 \leq t, \quad \vu_0 \in H^{\frac{\beta}{2}}(\Rt),
\end{equation}
where $\vu(t,\cdot)$ is the \emph{unique} Leray-type solution of the equation (\ref{Equation}) which arises from $\vu_0$. Due to the uniqueness of solutions, it is easy to verify that $(S(t))_{t\geq 0}$ is a strongly continuous semigroup on the Hilbert space $H^{\frac{\beta}{2}}(\Rt)$. 

\medskip 
We recall now the following:
\begin{Definition}[Strongly asymptotically compact semigroup]\label{def-asymptotically-compact}
	The semigroup $(S(t))_{t\geq 0}$ is  strongly  asymptotically compact 
	if for any bounded sequence $(\vu_{0,n})_{n\in \mathbb{N}}$ in $H^{\frac{\beta}{2}}(\Rt)$, and moreover,  for any sequence of times $(t_n)_{n \in \mathbb{N}}$ such that  $t_n \to \infty$ when $n\to \infty$, the sequence $(S(t_n)\vu_{0,n})_{n\in \mathbb{N}}$  is strongly precompact in $H^{\frac{\beta}{2}}(\Rt)$. 
\end{Definition}

Now we are able to state the following theorem on the existence of a strong global attractor. For a proof of this result see \cite{Raugel} and \cite{Temam}. 
\begin{Theoreme}[Existence of a strong global attractor]\label{Th:Existence-Attractor}
Assume that:
	\begin{enumerate}
		\item The semigroup $(S(t))_{t\geq 0}$ has a bounded and closed absorbing set $\mathcal{B} \subset H^{\frac{\beta}{2}}(\Rt)$. 
		\item The semigroup $(S(t))_{t\geq 0}$ is asymptotically compact in the sense of definition above. 
		\item  For every  $0 \leq t$ fixed,  the map $S(t) : \mathcal{B} \to H^{\frac{\beta}{2}}(\Rt)$ is continuous.
	\end{enumerate}
	Then, the semigroup $(S(t))_{t \geq 0}$ has a unique strong  global attractor $\mathcal{A}_{s} \subset  H^{\frac{\beta}{2}}(\Rt)$ given in Definition \ref{Def-global-attractor}.  Moreover,  the following statements hold true:
	\begin{enumerate}
    \item  The set 	$\mathcal{A}_s$ is invariant, \emph{i.e.}, for all $0 \leq t$ we have:
	\[ \mathcal{A}_s= \Big\{  \vu(t,\cdot) : \,\, \text{$\vu$ is a Leray-type solution of (\ref{Equation}) arising from}\,\, \vu_0 \in   \mathcal{A}_s  \Big\}. \]
	\item We have the  characterization of $\mathcal{A}_s$ given in (\ref{Charac-strong}). 
	\end{enumerate} 
\end{Theoreme} 
\noindent
{\bf Proof of the second point in Theorem \ref{Th3}}.  We shall prove that the semigroup $(S(t))_{t\geq 0}$  verify  all  the assumptions  in Theorem \ref{Th:Existence-Attractor}. For the reader's convenience, we will study  each of them separately. 

\medskip 

{\bf Point $1$}.  This point was already satisfied by Proposition \ref{Prop2}.

\medskip

{\bf Point $2$}.   Let $(\vu_{0,n})_{n\in \mathbb{N}}$ be a bounded sequence in $H^{\frac{\beta}{2}}(\Rt)$, and moreover, let $(t_n)_{n\in \mathbb{N}}$ be a sequence of positive times such that  $t_n \to +\infty$ when $n\to +\infty$. We must show that the sequence $(S(t_n)\vu_{0,n})_{n\in \mathbb{N}}$ is strongly precompact in the space $H^{\frac{\beta}{2}}(\Rt)$ and for this  we shall perform the following  energy method: for each $n \in \mathbb{N}$,  and for $J^{\alpha}_{\gamma}$, $J^{\beta}_{\delta}$ defined in (\ref{Notation-Operators}), we consider the following initial value  problem for the equation (\ref{Equation}): 
\begin{equation}\label{Bardina-n}
\left\{ \begin{array}{ll}\vspace{2mm}
\partial_t \vu_n + J^{\alpha}_{\gamma} \vu_n + J^{-\beta}_{\delta} \, \P\,  \text{div} \,(\vu_n \otimes \vu_n)  =\fe, \qquad \text{div}(\vu_n)=0, \\ 
\vu_n(-t_n,\cdot)= \vu_{0,n}.
\end{array} 
\right.
\end{equation}
By Theorem \ref{Th1}  there exists a unique  Leray-type solution  $\vu_n : [-t_n, +\infty[\times \Rt \to \Rt$.  	Moreover,   by definition of the semigroup $S(t)$ given in (\ref{Seligroup}),  for all $n\in \mathbb{N}$ we have the identity  $\ds{S(t_n)\vu_{0,n}= \vu_{n}(0, \cdot)}$. Therefore, we shall prove that  the sequence $(\vu_{n}(0,\cdot))_{n\in \mathbb{N}}$ is strongly precompact in $H^{\frac{\beta}{2}}(\Rt)$. Our general strategy is the following:  first, we  shall prove the existence of  an eternal  solution  associated to the equation (\ref{Equation}).  We recall that an eternal solution associated to  this equation is a function 
\begin{equation}\label{Eternal-sol-1}
\vu_e \in L^{\infty}_{loc}\Big(  ]-\infty, +\infty[, H^{\frac{\beta}{2}}(\Rt)\Big) \cap L^{2}_{loc}\Big(  ]-\infty, +\infty[, H^{\frac{\alpha+\beta}{2}}(\Rt)\Big), 
\end{equation}
which is a weak solution of  equation (\ref{Equation-Eternal}). Thus, we will show that the sequence $(\vu_{n}(0,\cdot))_{n\in \mathbb{N}}$ converges (via a sub-sequence) to $\vu_e(0,\cdot)$   in the strong topology of the  space $H^{\frac{\beta}{2}}(\Rt)$. 
 
 \medskip
 
 Our starting point is then to prove existence of an eternal solution:
\begin{Proposition}\label{Prop:eternal-sol} There exists a function  $\vu_e$ which verifies (\ref{Eternal-sol-1}) and (\ref{Equation-Eternal}).
\end{Proposition}
\pv  This function will be obtained as the limit when $n\to +\infty $ of the solutions $\vu_n : [-t_n, +\infty[\times \Rt \to \Rt$  to the initial value problems  (\ref{Bardina-n}).  By the first point in Proposition \ref{Prop1},   for all $n\in \mathbb{N}$  and for all  $ -t_n \leq t $  we have
\begin{equation*}
\Vert \vu_n(t,\cdot)\Vert^{2}_{H^{\frac{\beta}{2}}}  \leq  e^{- \gamma(t+t_n)}  \Vert \vu_{0,n} \Vert^{2}_{H^{\frac{\beta}{2}}}  + \frac{\b^2}{\a^2\gamma^2} \, \Vert \fe \Vert^{2}_{H^{\frac{\beta}{2}}}.
\end{equation*}
Moreover, as  the sequence $(\vu_{0,n})_{n\in \mathbb{N}}$ is bounded in $H^{\frac{\beta}{2}}(\Rt)$, there exists $0<R$ such that  we can write 
\begin{equation}\label{estim1-n}
\sup_{n \in \mathbb{N}} \,  \sup_{t \geq -t_n}  \, \Vert \vu_n(t,\cdot) \Vert^{2}_{H^{\frac{\beta}{2}}}  \leq    e^{-\gamma\, (t+t_n)}  R^2+ \frac{\b^2}{\a^2\gamma^2} \Vert \fe \Vert^{2}_{H^{\frac{\beta}{2}}} \leq R^2+ \frac{\b^2}{\a^2\gamma^2} \Vert \fe \Vert^{2}_{H^{\frac{\beta}{2}}}.
\end{equation}

On the other hand,  by the second point in Proposition  \ref{Prop1}, for all  $-t_n \leq t$ and $T=1$  we have 
\begin{equation*}
\c\, \int_{t}^{t+1} \Vert \vu_n(s,\cdot)\Vert^{2}_{H^{\frac{\alpha+\beta}{2}}} ds \leq  e^{- \gamma (t+t_n)}  \Vert \vu_{0,n} \Vert^{2}_{H^{\frac{\beta}{2}}}  + \left( \frac{\b^2}{\a^2\gamma^2}+ \frac{\b^2}{\a^2 \c} \right)  \Vert \fe \Vert^{2}_{H^{\frac{\beta}{2}}}, 
\end{equation*}
hence, for the constant $0<R$ above we get: 
\begin{equation}\label{estim2-n}
\sup_{n \in \mathbb{N}} \, 
\sup_{ t \geq -t_n} \left( \c \int_{t}^{t+1} \Vert \vu_n(s,\cdot)\Vert^{2}_{H^{\frac{\alpha+\beta}{2}}}\,  ds \right) \leq R^2+  \left( \frac{\b^2}{\a^2\gamma^2}+ \frac{\b^2}{\a^2 \c} \right)  \Vert \fe \Vert^{2}_{H^{\frac{\beta}{2}}}.  
\end{equation} 

In this fashion, by the estimates (\ref{estim1-n}) and (\ref{estim2-n}) and by the Banach-Alaoglu theorem, there exists $\vu_e \in L^{\infty}_{loc}(\R, H^{\frac{\beta}{2}}(\Rt)) \cap L^{2}_{loc}(\R, H^{\frac{\alpha+\beta}{2}}(\Rt))$ such that the sequence $(\vu_n)_{n\in \mathbb{N}}$ converges (via a sub-sequence) to $\vu_e$ in the weak $-*$ topology of the spaces $L^{\infty}([-\tau,\tau], H^{\frac{\beta}{2}}(\Rt))$   and $L^{2}([-\tau, \tau], H^{\frac{\alpha+\beta}{2}}(\Rt))$, for all $0<\tau$. Moreover,  as in the proof of Theorem \ref{Th1}, by using the Rellich-Lions lemma we obtain that  the limit $\vu_e$  is a weak solution of   equation (\ref{Equation-Eternal}).   \finpv

\medskip 

 We will prove now the convergence  (via a sub-sequence) of the  sequence $(\vu_n(0,\cdot))_{n\in \mathbb{N}}$   to $\vu_e(0,\cdot)$ in the strong topology of the space  $H^{\frac{\beta}{2}}(\Rt)$. For this, by following the same ideas of (\ref{Iden}) we have: 
\begin{Lemme}\label{Lemme:Identity} Let $\frac{5}{2} \leq \alpha +\beta$. Moreover,  let $D(m_1)$ and $D(m_2)$  be the pseudo-differential operators of order zero defined in (\ref{Symbols}). Then, for all $0 < t$, Leray-type solutions of  equation (\ref{Equation}) verify the identity: 
\begin{equation*}
\begin{split}
 &\,\frac{d}{dt} \Vert D(m^{1/2}_{1}) \vu(t,\cdot) \Vert^{2}_{H^{\frac{\beta}{2}}} + 2 \left\Vert D(m^{1/2}_{1})D(m^{1/2}_{2}) \vu(t,\cdot) \right\Vert^{2}_{H^{\frac{\alpha+\beta}{2}}} \\
 = &\,2 \left\langle D(m^{1/2}_{1}) (I_d -\Delta)^{\frac{\beta}{4}} \fe(t,\cdot), D(m^{1/2}_{1})(I_d -\Delta)^{\frac{\beta}{4}} \vu(t,\cdot)  \right\rangle_{L^2 \times L^2}.
 \end{split}
\end{equation*} 
Moreover, the eternal solution $\vu_e$ constructed in Proposition  \ref{Prop:eternal-sol}  also verifies this identity. 
\end{Lemme} 	

We multiply each term in this identity by $e^{2 \, t}$, and moreover,  we integrate  in the interval $[-t_n, 0]$ to get:
\begin{equation*}
\begin{split}
& \,  \Vert D(m^{1/2}_{1}) \vu_n(0,\cdot)\Vert^{2}_{H^{\frac{\beta}{2}}} -  e^{- 2  \,  t_n }\Vert D(m^{1/2}_{1}) \vu_{0,n} \Vert^{2}_{H^{\frac{\beta}{2}}}- 2  \int_{-t_n}^{0} e^{2 \,  t} \Vert D(m^{1/2}_{1}) \vu_n(t,\cdot)\Vert^{2}_{H^{\frac{\beta}{2}}}  dt  \\
&\, + 2\int_{-t_n}^{0} e^{2\,   t} \left\Vert  D(m^{1/2}_{1}) D(m^{1/2}_{2})\vu_n(t,\cdot) \right\Vert^{2}_{H^{\frac{\alpha+\beta}{2}}}  dt  \\
=  &\,  2 \int_{-t_n}^{0} e^{2\, t} \left\langle D(m^{1/2}_{1}) (I_d -\Delta)^{\frac{\beta}{4}} \fe,  D(m^{1/2}_{1}) (I_d -\Delta)^{\frac{\beta}{4}} \vu_n(t,\cdot)\right\rangle_{L^2\times L^2} \,  dt.
\end{split}
\end{equation*}
By applying the $\ds{\limsup}$ when $n\to +\infty$ in each term of this identity  we obtain:  
\begin{equation}\label{Estim-energ-limsup}
\begin{split}
&\limsup_{n\to +\infty} \Vert D(m^{1/2}_{1}) \vu_n(0,\cdot)\Vert^{2}_{H^{\frac{\beta}{2}}}\\
\leq & \,  \limsup_{n\to +\infty}   e^{- 2  t_n }\Vert D(m^{1/2}_{1}) \vu_{0,n} \Vert^{2}_{H^{\frac{\beta}{2}}}   +  \limsup_{n\to +\infty} \left( 2\,\int_{-t_n}^{0}  e^{2  t} \Vert D(m^{1/2}_{1}) \vu_n(t,\cdot) \Vert^{2}_{H^{\frac{\beta}{2}}}  dt \right) \\
&\, + \limsup_{n\to +\infty} \left( - 2 \,\int_{-t_n}^{0} e^{2  t} \left\Vert D(m^{1/2}_{1})D(m^{1/2}_{2})\,\vu_n(t,\cdot) \right\Vert^{2}_{\frac{\alpha+\beta}{2}}  dt \right)\\
&\, +\limsup_{n\to +\infty} \left(2 \int_{-t_n}^{0} e^{2  t} \left\langle D(m^{1/2}_{1})(I_d -\Delta)^{\frac{\beta}{4}} \fe, D(m^{1/2}_{1}) (I_d -\Delta)^{\frac{\beta}{4}}\vu_n(t,\cdot)\right\rangle_{L^2\times L^2}  dt\right),
\end{split}
\end{equation}
where we must study each term on the right side. For the first term, always by the fact that  the sequence $(\vu_{0,n})_{n\in \mathbb{N}}$ is bounded in $H^{\frac{\beta}{2}}(\Rt)$, we have
\begin{equation}\label{E1}
\limsup_{n\to +\infty}   e^{- 2  t_n }\Vert D(m^{1/2}_{1}) \vu_{0,n} \Vert^{2}_{H^{\frac{\beta}{2}}} =0.    
\end{equation}
For the second  term,   by the  estimate (\ref{estim1-n}) the sequence $(\vu_n)_{n \in \mathbb{N}}$ converges to $\vu_e$ in the weak$-*$ topology of the space $L^{\infty}_{loc}(\R, H^{\frac{\beta}{2}}(\Rt))$; and then  it converges in the weak$-*$ topology of the space  $L^{2}_{loc}(\R, H^{\frac{\beta}{2}}(\Rt))$. We thus have:  
\begin{equation*}
\liminf_{n\to +\infty} \left( 2\,\int_{-t_n}^{0}  e^{2  t} \Vert D(m^{1/2}_{1}) \vu_n(t,\cdot) \Vert^{2}_{H^{\frac{\beta}{2}}}  dt \right) \geq   2 \int_{-\infty}^{0} e^{2t}\Vert  D(m^{1/2}_{1}) \vu_e(t,\cdot) \Vert^{2}_{H^{\frac{\beta}{2}}} dt, 
\end{equation*}
hence we can write 
\begin{equation*}
\liminf_{n\to +\infty} \left( 2\,\int_{-t_n}^{0}  e^{2  t} \Vert D(m^{1/2}_{1}) \vu_n(t,\cdot) \Vert^{2}_{H^{\frac{\beta}{2}}}  dt \right) \leq   - 2 \int_{-\infty}^{0} e^{2t}\Vert D(m^{1/2}_{1}) \vu_e(t,\cdot) \Vert^{2}_{H^{\frac{\beta}{2}}} dt.
\end{equation*}
Similarly, for the third term,   by the  estimate (\ref{estim2-n})  the sequence $(\vu_n)_{n \in \mathbb{N}}$ converges to $\vu_e$ in the weak$-*$ topology of the space $L^{2}_{loc}(\R, H^{\frac{\alpha+\beta}{2}}(\Rt))$, then we have
\begin{equation}\label{E2}
\begin{split}
&\, \limsup_{n\to +\infty} \left( -2\,  \int_{-t_n}^{0} e^{2  t} \Vert D(m^{1/2}_{1})D(m^{1/2}_{2}) \vu_n(t,\cdot) \Vert^{2}_{H^{\frac{\alpha+\beta}{2}}}  dt \right) \\
\leq &\,  -2\,  \int_{-\infty}^{0} e^{2\, t} \Vert D(m^{1/2}_{1})D(m^{1/2}_{2}) \vu_e(t,\cdot)\Vert^{2}_{H^{\frac{\alpha+\beta}{2}}} dt. 
\end{split}
\end{equation}
Moreover, for the fourth term we obtain
\begin{equation}\label{E4}
\begin{split}
&\, \limsup_{n\to +\infty} \left(2 \int_{-t_n}^{0} e^{2  t} \left\langle D(m^{1/2}_{1})(I_d -\Delta)^{\frac{\beta}{4}} \fe,  D(m^{1/2}_{1}) (I_d -\Delta)^{\frac{\beta}{4}}\vu_n(t,\cdot)\right\rangle_{L^2\times L^2}  dt\right)\\
= &\,  2 \int_{-\infty}^{0}  e^{2  t} \left\langle D(m^{1/2}_{1}) (I_d -\Delta)^{\frac{\beta}{4}} \fe, D(m^{1/2}_{1})(I_d -\Delta)^{\frac{\beta}{4}}\vu_e(t,\cdot)\right\rangle_{L^2\times L^2} \, dt, 
\end{split}
\end{equation} 
Thus, with these estimates  we get back to (\ref{Estim-energ-limsup}) to write: 
\begin{equation*}
\begin{split}
\limsup_{n\to +\infty} \Vert D(m^{1/2}_{1}) \vu_n(0,\cdot)\Vert^{2}_{H^{\frac{\beta}{2}}}\leq & \,\,  2  \int_{-\infty}^{0} e^{2\beta t} \Vert D(m^{1/2}_{1}) \vu_e(t,\cdot)\Vert^{2}_{H^{\frac{\beta}{2}}} dt\\
&\,  -2\,  \int_{-\infty}^{0} e^{2  t} \Vert  D(m^{1/2}_{1})D(m^{1/2}_{2}) \vu_e(t,\cdot)\Vert^{2}_{H^{\frac{\alpha+\beta}{2}}}\, dt \\
&+ 2 \int_{-\infty}^{0}  e^{2 t} \left\langle  (I_d -\Delta)^{\frac{\beta}{4}}\fe, (I_d -\Delta)^{\frac{\beta}{4}}\vu_e(t,\cdot)\right\rangle_{L^2\times L^2} \, dt =(A). 
\end{split}
\end{equation*} 
We shall study now the term $(A)$. By Lemma \ref{Lemme:Identity} the  eternal solution $\vu_e$  of  equation (\ref{Equation-Eternal}) verifies the identity:
\begin{equation*}
\begin{split}
 \frac{d}{dt} \Vert  D(m^{1/2}_{1})\vu_e(t,\cdot)\Vert^{2}_{H^{\frac{\beta}{2}}} = &\,  - 2 \left\Vert D(m^{1/2}_{1}) D(m^{1/2}_{2})\vu_e(t,\cdot) \right\Vert^{2}_{H^{\frac{\alpha+\beta}{2}}} \\
 &\,+ 2 \left\langle  D(m^{1/2}_{1})(I_d -\Delta)^{\frac{\beta}{4}}\fe, D(m^{1/2}_{1}) (I_d -\Delta)^{\frac{\beta}{4}}\vu_e(t,\cdot)\right\rangle_{L^2\times L^2}.  
 \end{split}
\end{equation*}
We multiply each term by $e^{2  t}$, then  we integrate  in the interval $]-\infty, 0]$ to get: 
\begin{equation*}
\begin{split}
 \Vert D(m^{1/2}_{1}) \vu_e(0,\cdot)\Vert^{2}_{H^{\frac{\beta}{2}}}  = & \,   \int_{-\infty}^{0} e^{2  t} \Vert D(m^{1/2}_{1}) \vu_e(t,\cdot)\Vert^{2}_{H^{\frac{\beta}{2}}}  dt   - 2 \int_{-\infty}^{0} e^{2  t} \left\Vert D(m^{1/2}_{1}) D(m^{1/2}_{2})\, \vu_e(t,\cdot) \right\Vert^{2}_{H^{\frac{\alpha+\beta}{2}}} \,  dt  \\
 & + 2\, \int_{-\infty}^{0} e^{2  t} \left\langle  D(m^{1/2}_{1})(I_d -\Delta)^{\frac{\beta}{4}}\fe,  D(m^{1/2}_{1})(I_d -\Delta)^{\frac{\beta}{4}}\vu_e(t,\cdot)\right\rangle_{L^2\times L^2}  \,  dt= (A). 
\end{split}
\end{equation*}
 
In this fashion,  by the previous estimate  we get  $\ds{\limsup_{n\to +\infty} \Vert D(m^{1/2}_{1}) \vu_n(0,\cdot)\Vert^{2}_{H^{\frac{\beta}{2}}}\leq \Vert D(m^{1/2}_{1})\vu_e(0,\cdot)\Vert^{2}_{H^{\frac{\beta}{2}}}}$. Moreover,   as  sequence  $(\vu_n)_{n\in \mathbb{N}}$ converges (via a sub-sequence) to $\vu_e$ in the weak$-*$ topology of the space $L^{\infty}(\R, H^{\frac{\beta}{2}}(\Rt))$, we are able to write   $\ds{\Vert D(m^{1/2}_{1}) \vu_e(0,\cdot)\Vert^{2}_{H^{\frac{\beta}{2}}} \leq \liminf_{n\to +\infty} \Vert D(m^{1/2}_{1}) \vu_{n}(0,\cdot)\Vert^{2}_{H^{\frac{\beta}{2}}}}$. We thus  obtain the desired strong convergence: $\ds{\lim_{n\to +\infty} \Vert D(m^{1/2}_{1}) \vu_n(0,\cdot)\Vert^{2}_{H^{\frac{\beta}{2}}}= \Vert D(m^{1/2}_{1}) \vu_e(0,\cdot)\Vert^{2}_{H^{\frac{\beta}{2}}}}$. 
			
\medskip 	

{\bf Point $3$}. The continuity of the map $S(t) : \mathcal{B} \to H^{\frac{\beta}{2}}(\Rt)$  directly follows from the estimate (\ref{Uniqueness1}),  where we have $ \vw(t,\cdot)=\vu_1(t,\cdot)-\vu_2(t,\cdot)= S(t)\vu_{0,1}-S(t)\vu_{0,2}$, and $\vw(0,\cdot)=\vu_{0,1}-\vu_{0,2}$.  

\medskip 

At this point, we are able to apply   Theorem \ref{Th:Existence-Attractor} to deduce that  the semigroup $(S(t))_{t\geq 0}$ has a unique strong global attractor $\mathcal{A}_{s}$. The second point of  Theorem \ref{Th3} is  now proven. \finpv 

\section{Stationary solutions}\label{Sec:Stationary}
\subsection{Proof of Theorem \ref{Th3}} 
 We shall use the following approximated equation.  Let $\theta \in \mathcal{C}^{\infty}_{0}(\Rt)$ be such that $0\leq \theta \leq 1$,  $\theta(x)= 1$ when $\vert x \vert \leq 1$ and $\theta(x)=0$ when $2 \leq \vert x \vert$. For $0<R$ we define the cut-off function $\theta_R(x)= \theta\left( \frac{x}{R}\right)$. Then, for $0<\varepsilon$ we consider the approximated problem:
\begin{equation}\label{Stationary-Approx}
-\varepsilon \Delta \U + \nu(-\Delta)^{\frac{\alpha}{2}} \U + J^{-\beta}_{\delta} \P  \left( (\theta_R \U \cdot \vec{\nabla}) \theta_R \U \right)  = \fe - \gamma \U, \quad \text{div}(\U)=0, \quad  0 <\varepsilon, R,
\end{equation}
where the operator $J^{-\beta}_{\delta}$ is defined in (\ref{Notation-Operators}).  Remark that when  $R \to +\infty$ and $\varepsilon \to 0^{+}$, solutions to the approximated equation (\ref{Stationary-Approx}) formally converge to solutions of  equation (\ref{Stationary}). 

\medskip

For $0 <\varepsilon, R$ fixed, our starting point is to construct solutions of  equation (\ref{Stationary-Approx}). For this, we shall use the following theorem. For a proof of this result we refer to  Theorem $16.1$, page $529$ of \cite{PGLR1}. 
\begin{Theoreme}[Sheafer's fixed point]\label{Th-Sheafer} Let $E$ be a Banach space and let $T:E \to E$ such that: 
	\begin{enumerate}
		\item[1)] $T$ is a continuous operator and compact operator.
		\item[2)] There exists a constant $0<M$, such that for all $\lambda \in [0,1]$, if $e\in E$ verifies $e=\lambda T(e)$ then we have $\Vert e \Vert_E \leq M$.
	\end{enumerate}
	Then, there exists $e  \in E$ such that $e=T(e)$.  
\end{Theoreme}
 Within the framework of this theorem, we set the Banach space $E=\left\{ \U \in H^{1}(\Rt): \,\, \text{div}(\U)=0 \right\}$. Moreover, we rewrite equation (\ref{Stationary-Approx}) as 
\[ -\varepsilon \Delta \U  + \frac{\gamma}{2} \U+ \nu(-\Delta)^{\frac{\alpha}{2}} \U +  \frac{\gamma}{2} \U+ J^{-\beta}_{\delta} \P  \left( (\theta_R \U \cdot \vec{\nabla}) \theta_R \U \right)  = \fe, \]
hence, using the operator $J^{\alpha}_{\gamma/2}$ given in (\ref{Notation-Operators}) we get 
\[  \left(\frac{\gamma}{2} I_d - \varepsilon \Delta \right) \U +  J^{\alpha}_{\gamma/2} \U + J^{-\beta}_{\delta}\P  \left( (\theta_R \U \cdot \vec{\nabla}) \theta_R \U \right)  = \fe.\]
Then  we have the  following (equivalent) fixed point problem
\begin{equation}\label{U-fixed-point-Approximated}
\U=  \frac{J^{-\beta}_{\delta}}{\left( \frac{\gamma}{2} I_d - \varepsilon \Delta\right) + J^{\alpha}_{\gamma/2}}\left(  \P\left(  (\theta_R \U \cdot \vec{\nabla}) \theta_R \U \right)\right) + \frac{1}{\left(\frac{\gamma}{2} I_d - \varepsilon \Delta\right)+ J^{\alpha}_{\gamma/2}} \left( \fe \,\right)=T_{\varepsilon,R}(\U). 
\end{equation} 

In the following technical lemmas, we  verify each point of Theorem \ref{Th-Sheafer}.
\begin{Lemme} For $0<\varepsilon$ and $0<R$, the operator $T_{\varepsilon,R}(\cdot)$ defined in (\ref{U-fixed-point-Approximated}) is continuous and compact in the space $E$. 
\end{Lemme}	
\pv From the expression (\ref{U-fixed-point-Approximated})  we can write 
\begin{equation*}
\begin{split}
T_{\varepsilon,R}(\U)=  & \,  \frac{-\Delta \,  J^{-\beta}_{\delta}}{\left( \frac{\gamma}{2} I_d - \varepsilon \Delta\right) + J^{\alpha}_{\gamma/2}}\left( \frac{1}{-\Delta}  \P\left(  (\theta_R \U \cdot \vec{\nabla}) \theta_R \U \right)\right) + \frac{1}{\left(\frac{\gamma}{2} I_d - \varepsilon \Delta\right)+ J^{\alpha}_{\gamma/2}}  \left( \fe \,\right),
\end{split}
\end{equation*}
and denoting $\ds{T_R(\U)= \frac{1}{-\Delta}  \P\left(  (\theta_R \U \cdot \vec{\nabla}) \theta_R \U \right)}$ we get 
\begin{equation*}
\begin{split}
T_{\varepsilon,R}(\U)=  &\,   \frac{-\Delta \,  J^{-\beta}_{\delta}}{\left( \frac{\gamma}{2} I_d - \varepsilon \Delta\right) + J^{\alpha}_{\gamma/2}}\left(  T_{R} (\U)\right) + \frac{1}{\left(\frac{\gamma}{2} I_d - \varepsilon \Delta\right)+ J^{\alpha}_{\gamma/2}}   \left( \fe \,\right).
\end{split}
\end{equation*}
Hence, for any $\U_1, \U_2\in E$ we have
\[ \|  T_{\varepsilon,R}(\U_1) -   T_{\varepsilon,R}(\U_2) \|_{H^1}= \left\|   \frac{-\Delta \,  J^{-\beta}_{\delta}}{\left( \frac{\gamma}{2} I_d - \varepsilon \Delta\right) + J^{\alpha}_{\gamma/2}}\left(  T_{R} (\U_2)- T_{R} (\U_1)\right)  \right\|_{H^1}.\]
Moreover, since operator $\ds{\frac{-\Delta \,  J^{-\beta}_{\delta}}{\left( \frac{\gamma}{2} I_d - \varepsilon \Delta\right) + J^{\alpha}_{\gamma/2}}(\cdot)}$ has a bounded symbol in the Fourier variable we directly have
\[ \|  T_{\varepsilon,R}(\U_1) -   T_{\varepsilon,R}(\U_2) \|_{H^1} \lesssim  \left\|    T_{R} (\U_2)- T_{R} (\U_1)\right\|_{H^1}.\]
From the proof of Theorem $16.2$, page $530$ of \cite{PGLR1} we known that $T_R$ is continuous and compact in the space $E$, which yields the wished result. \finpv 
\begin{Lemme} Let $\lambda \in [0,1]$, and let $\U \in E$ be such that $\U= \lambda T_{\varepsilon,R}(\U)$. There exists a constant $0<C_\gamma$, which only depends on $\gamma$, such that the following estimate holds:
	\begin{equation}\label{Energy-Stationary}
	\varepsilon\, C_\gamma\,    \Vert \U \Vert^{2}_{H^{1+ \frac{\beta}{2}}} + \frac{\c }{2} \,\Vert \U \Vert^{2}_{H^{\frac{\alpha+\beta}{2}}}   \leq \frac{\b^2}{2\a^2\c}  \Vert \fe \Vert^{2}_{H^{\frac{\beta}{2}}}, 
	\end{equation}
	with  the constants $\a,\b$ and $\c$ defined in (\ref{Notation}). 
\end{Lemme}	
\pv  Let $\U \in E$ be such that $\U = \lambda T_{\varepsilon,R} (\U)$, with $0 \leq \lambda \leq 1$.  Then $\U$ solves the equation
\begin{equation*}
-\varepsilon \Delta \U + \nu(-\Delta)^{\frac{\alpha}{2}} \U + \lambda\, (I_d-\delta^2\Delta)^{-\frac{\beta}{2}} \P  \left( (\theta_R \U \cdot \vec{\nabla}) \theta_R \U \right)  =  \lambda\, \fe - \gamma \U.
\end{equation*}
Recalling the operators $J^{\beta}_{\delta}$ and $J^{-\beta}{\delta}$ defined in (\ref{Notation-Operators}) we can write  
\begin{equation*}
-\varepsilon \Delta  J^{\beta}_{\delta} \U + \nu(-\Delta)^{\frac{\alpha}{2}} J^{\beta}_{\delta} \U + \lambda\, \P  \left( (\theta_R \U \cdot \vec{\nabla}) \theta_R \U \right)  =  \lambda\, J^{\beta}_{\delta}  \fe - \gamma\, J^{\beta}_{\delta} \U.
\end{equation*}
Then,  we multiply this equation by $\U$, integrating by parts and using the lower and upper bounds in (\ref{Bounds}) we obtain (\ref{Energy-Stationary}).
\finpv 

\medskip

Now, remark that from the estimate (\ref{Energy-Stationary}) and by  the continuous embedding $H^{1+\frac{\beta}{2}}(\Rt) \subset H^1(\Rt)$,  we are able to write $\ds{\Vert \U \Vert^{2}_{H^1} \leq \frac{1}{\varepsilon\, C_\gamma}  \frac{\b^2}{2\a^2\c}  \Vert \fe \Vert^{2}_{H^{\frac{\beta}{2}}}}$. We set  the constant $\ds{M=\frac{1}{\varepsilon\, C_\gamma} \frac{\b^2}{2\a^2\c}  \Vert \fe \Vert^{2}_{H^{\frac{\beta}{2}}}}$ to verify  the second point of Theorem \ref{Th-Sheafer}.   

\medskip

By theorem \ref{Th-Sheafer} there exists $\U= \U_{\varepsilon,R} \in E$ a solution to the fixed point problem $\U_{\varepsilon, R}= T_{\varepsilon,R}(\U_{\varepsilon,R})$. Moreover, this solutions also solve the elliptic problem (\ref{Stationary-Approx}). Consequently, by the estimate (\ref{Energy-Stationary}) we also have the uniform  control (in $0<\varepsilon$ and $0<R$): \[\Vert \U_{\varepsilon, R} \Vert^{2}_{H^{\frac{\alpha+\beta}{2}}} \leq  \frac{\b^2}{\a^2 \c^2} \Vert \fe \Vert^{2}_{H^{\frac{\beta}{2}}}.\]  
Thus, the end of the proof follows standard arguments: first, we set $0<\varepsilon$, by the uniform control above and by the Rellich-Lions lemma there exists $\U_\varepsilon \in H^{\frac{\alpha+\beta}{2}}(\Rt)$ and a sequence $R_n \to +\infty$ such that $\U_{\varepsilon, R_n}$ converges to $\U_\varepsilon$ in the weak topology of the space $H^{\frac{\alpha+\beta}{2}}(\Rt)$ and in the strong topology of the space $L^{2}_{loc}(\Rt)$. Consequently, the limit $\U_\varepsilon$ solves the equation:
\[ -\varepsilon \Delta \U_\varepsilon + \nu(-\Delta)^{\frac{\alpha}{2}} \U_\varepsilon + J^{-\beta}_{\delta}\P  \left( (\U_\varepsilon \cdot \vec{\nabla})  \U_\varepsilon \right)  = \fe - \gamma \U_\varepsilon, \quad \text{div}(\U)=0, \quad  0 <\varepsilon.  \] 
Similarly, by the uniform  control $\ds{\Vert \U_{\varepsilon} \Vert^{2}_{H^{\frac{\alpha+\beta}{2}}} \leq \frac{\b^2}{\a^2 \c^2} \Vert \fe \Vert^{2}_{H^{\frac{\beta}{2}}}}$ and using again the Rellich-Lions lemma  the family $\U_\varepsilon$ converges to a solution $\U \in H^{\frac{\alpha+\beta}{2}}(\Rt)$ of the equation (\ref{Stationary}). Moreover, this solution also verifies the energy estimate $\ds{\Vert \U \Vert^{2}_{H^{\frac{\alpha+\beta}{2}}} \leq  \frac{\b^2}{\a^2 \c^2} \Vert \fe \Vert^{2}_{H^{\frac{\beta}{2}}}}$. Theorem \ref{Th3}  is now proven. \finpv

\subsection{Proof of Theorem \ref{Th4}}
Let $\vu\in L^{\infty}_{t}H^{\frac{\beta}{2}}_{x}\cap (L^{2}_{t})_{loc}H^{\frac{\alpha+\beta}{2}}_{x}$ be a Leray-type solution to equation (\ref{Equation})  (obtained in Theorem \ref{Th1}) and let  $\U \in H^{\frac{\alpha+\beta}{2}}(\Rt)$ be a stationary solution of equation (\ref{Stationary})  (obtained in Theorem \ref{Th3}).  We define $\vw(t,\cdot)=\vu(t,\cdot)-\U$. This function solves the equation:
\begin{equation}\label{Eq-Difference}
\partial_{t} \vw + \nu (-\Delta)^{\frac{\alpha}{2}} \vw + J^{-\beta}_{\delta} \P \left( \text{div}(\vu \otimes \vu) - \text{div}(\U \otimes \U)  \right) + \gamma \vw=0, \quad \text{div}(\vw)=0.
\end{equation}

The key estimate to prove Theorem \ref{Th4} is the following new energy control: 
\begin{Proposition}\label{Prop:Key-Estimate} Let $\vw \in  L^{\infty}_{t}H^{\frac{\beta}{2}}_{x}\cap (L^{2}_{t})_{loc}H^{\frac{\alpha+\beta}{2}}_{x}$ be a weak solution to  equation (\ref{Eq-Difference}). Then, $\vw$ verifies: 
	\begin{equation}\label{Key-Estim}
	\begin{split}
	\Vert J^{\frac{\beta}{2}}_{\delta}  \vw(t,\cdot) \Vert^{2}_{H^{L^2}} \leq & \, \,  \Vert J^{\frac{\beta}{2}}_{\delta} ( \vu_0 - \U) \Vert^{2}_{L^2}   -2\nu \, \int_{0}^{t}\Vert (-\Delta)^{\frac{\alpha}{4}} J^{\frac{\beta}{2}}_{\delta}  \vw(s,\cdot) \Vert^{2}_{L^2}\, ds \\
	&- 2\gamma \int_{0}^{t} \Vert  J^{\frac{\beta}{2}}_{\delta} \vw(s,\cdot) \Vert^{2}_{L^2}\, ds + C \frac{\b}{\a \c} \| \fe \|_{H^{\frac{\beta}{2}}} \int_{0}^{t} \Vert \vw(s,\cdot) \Vert^{2}_{H^{\frac{\alpha+\beta}{2}}} \, ds,
	\end{split}
	\end{equation}
where the operator $J^{\beta}_{\delta}$ is defined (\ref{Notation-Operators}) and the quantities $\a,\b,\c$ are given in (\ref{Notation}).	
\end{Proposition}	
\pv  We write
\[\text{div}(\vu \otimes \vu) - \text{div}(\U \otimes \U) =  \text{div}(\vw \otimes \vw)+ \text{div}(\vw \otimes \U)+ \text{div}(\U\otimes \vw),\]
hence the function $\vw$ verifies the equation:
\begin{equation}\label{eq-w}
\partial_{t} \vw + \nu (-\Delta)^{\frac{\alpha}{2}} \vw + J^{-\beta}_{\delta} \P \left(  \text{div}(\vw \otimes \vw)+ \text{div}(\vw \otimes \U)+ \text{div}(\U\otimes \vw)  \right) + \gamma \vw=0, \quad \text{div}(\vw)=0.
\end{equation}

To prove (\ref{Key-Estim}) our starting point is  the  direct identity 
\begin{equation}\label{iden1}
\Vert  J^{\frac{\beta}{2}}_{\delta}\vw(t,\cdot) \Vert^{2}_{L^2}= \Vert  J^{\frac{\beta}{2}}_{\delta} \vu(t,\cdot) \Vert^{2}_{L^2}-2 \left(   J^{\frac{\beta}{2}}_{\delta}\vu(t,\cdot) ,   J^{\frac{\beta}{2}}_{\delta}\U \right)_{L^2} + \Vert  J^{\frac{\beta}{2}}_{\delta} \U \Vert^{2}_{L^2},
\end{equation}
Here, we must  study  the second term on the right-hand  side, and for this we shall proof the following identity:
\begin{equation}\label{Identity-1}
- 2 \left(  J^{\frac{\beta}{2}}_{\delta} \vu(t,\cdot) ,  J^{\frac{\beta}{2}}_{\delta}\U \right)_{L^2} = -2 \left( J^{\frac{\beta}{2}}_{\delta} \vu_0,  J^{\frac{\beta}{2}}_{\delta}\U\right)_{L^2} - 2 \int_{0}^{t}\left\langle \partial_t \, J^\beta_\delta \vw(s,\cdot) , \U \right\rangle_{H^{-\frac{3}{2}} \times H^{\frac{3}{2}}} \, ds. 
\end{equation}
Indeed, let us start by  verifying that the last  term in the right-hand side is well-defined. This fact follows from the next technical lemmas.   
\begin{Lemme}\label{LemThec1} Let $2\leq \alpha +\beta$  and let $\U \in H^{\frac{\alpha+\beta}{2}}(\Rt)$ be the solution of equation (\ref{Stationary}). Then  we have the following gain of regularity $\U \in H^{\alpha+\frac{\beta}{2}}(\Rt)$. 
\end{Lemme} 
\pv We write $\U$ as the solution of the fixed point problem
\[ \U= \frac{J^{-\beta}_{\delta}}{ J^{\alpha}_{\gamma}} \left(  \P \,  \text{div}(\U \otimes \U)\right) + \frac{1}{ J^{\alpha}_{\gamma}} \left( \fe \,\right),\]
where the operators $J^{\alpha}_{\gamma}$ and $J^{-\beta}_{\delta}$ are given in (\ref{Notation-Operators}). Moreover, by the identities (\ref{Identities}) we obtain
\[ \U=  D(m^{-1}_{1}) D(m^{-1}_{2}) (I_d - \Delta)^{-\frac{\alpha+\beta}{2}} \left(  \P \,  \text{div}(\U \otimes \U)\right) + D(m^{-1}_{2}) (I_d -\Delta)^{-\frac{\alpha}{2}} \fe.\]

\medskip

In this last identity, since $\fe \in H^{\frac{\beta}{2}}(\Rt)$ then we have $D(m^{-1}_{2})(I_d-\Delta)^{\frac{\alpha}{2}} \fe \in H^{\alpha+\frac{\beta}{2}}(\Rt)$. Moreover, as $\U \in H^{\frac{\alpha+\beta}{2}}(\Rt)$ (with $1\leq \frac{\alpha+\beta}{2}$) by the product laws in  Sobolev spaces and by an iterative argument we obtain $D(m^{-1}_{1})D(m^{-1}_{2})(I_d -\Delta)^{-\frac{\alpha+\beta}{2}} \P \, \text{div} (\U \otimes \U)  \in H^{\alpha+\frac{\beta}{2}}(\Rt)$. We thus have $\U \in H^{\alpha+\frac{\beta}{2}}(\Rt)$.  It is worth emphasizing this gain of regularity of $\U$ is sharp in the sense that the term $D(m^{-1}_{2})(I_d-\Delta)^{\frac{\alpha}{2}} \fe$ only belongs to the space $H^{\alpha+ \frac{\beta}{2}}(\Rt)$, provided that $\fe \in H^{\frac{\beta}{2}}(\Rt)$.  \finpv 

\begin{Lemme}\label{LemTech2} Let $ 2 \leq \alpha + \beta$ and let $\vu$ be a Leray-type solution of  equation (\ref{Equation}). Then  we have  $\partial_t \, J^\beta_\delta \vu \in ( L^{2}_{t})_{loc} H^{-\frac{3}{2}}_{x}$. 
\end{Lemme}	
\pv  Recall that $\vu$ solves the equation: 
\[ \partial_t J^{\beta}_\delta \vu= -\nu (-\Delta)^{\frac{\alpha}{2}} J^{\beta}_\delta\, \vu - \gamma \, J^\beta_\delta \, \vu-\P ( \text{div}(\vu \otimes \vu) ) + J^{\beta}_\delta\fe, \]
where we must verify that  each term on the right-hand side belong to the space $(L^{2}_{t})_{loc} H^{-\frac{3}{2}}_{x}$.  Without loss of generality, we shall assume that $2 \leq \alpha + \beta < \frac{5}{2}$. The case $\frac{5}{2}\leq \alpha +\beta$ is directly treated as the previous one by setting $\alpha' \leq \alpha$ and $\beta' \leq \beta$ such that  $2 \leq \alpha' + \beta' < \frac{5}{2}$. 

\medskip

As $\vu \in (L^{2}_{t})_{loc} H^{\frac{\alpha+\beta}{2}}_{x}$, and as $\alpha+\beta <\frac{5}{2}$,  we have $-\nu (-\Delta)^{\frac{\alpha}{2}} J^{\beta}_\delta\, \vu  \in (L^{2}_{t})_{loc} H^{-\frac{\alpha+\beta}{2}}_{x} \subset (L^{2}_{t})_{loc} H^{-\frac{3}{2}}_{x}$. We also have $-\gamma J^\beta_\delta \, \vu \in (L^{2}_{t})_{loc} H^{\frac{\alpha-\beta}{2}}_{x} \subset (L^{2}_{t})_{loc} H^{-\frac{3}{2}}_{x}$. Thereafter, as $2 \leq \alpha+\beta<\frac{5}{2}$ by Lemma \ref{Lemma-uniq}  we have $-\P ( \text{div}(\vu \otimes \vu) )  \in (L^{2}_{t})_{loc} H^{\frac{-5 +(\alpha+\beta)}{2}}_{x} \subset (L^{2}_{t})_{loc} H^{-\frac{3}{2}}_{x}$. Finally, as $\fe \in H^{\frac{\beta}{2}}(\Rt)$, and  as $\beta \leq 3$ (since we have $\beta<\alpha+\beta<\frac{5}{2}$), we get  $J^\beta_\delta \fe \in (L^{2}_{t})_{loc} H^{-\frac{\beta}{2}}_{x}  \subset (L^{2}_{t})_{loc} H^{-\frac{3}{2}}_{x}$. \finpv 

\medskip

By Lemma \ref{LemThec1} and by the assumption $\frac{3}{2} < \alpha+\frac{\beta}{2}$ we have $\U \in H^{\frac{3}{2}}(\Rt)$. Moreover, by Lemma \ref{LemTech2} and the fact that  $\U$ is a time independent function we have 
$\partial_t \, J^\beta \vw = \partial _t \, J^\beta \vu \in (L^{2}_{t})_{loc} H^{-\frac{3}{2}}_{x}$. In this fashion, the term $\ds{\int_{0}^{t}\left\langle \partial_t\, J^\beta_\delta \vw(s,\cdot) , \U \right\rangle_{H^{-\frac{3}{2}} \times H^{\frac{3}{2}}} \, ds}$ is well defined for all $0\leq t$.

\medskip

Now,  for \emph{a.e.} $0 \leq s$ we can write 
\begin{equation*}
\begin{split}
\partial_t \left(  J^{\frac{\beta}{2}}_{\delta} \vu(s,\cdot) ,   J^{\frac{\beta}{2}}_{\delta}\U \right)_{L^2} = &  \partial_t  \left\langle   J^\beta_\delta \vu(s,\cdot) , \U \right\rangle_{H^{-\frac{\beta}{2}} \times H^{\frac{\beta}{2}}} =\left\langle \partial_t \, J^\beta_\delta \vu(s,\cdot) , \U \right\rangle_{H^{-\frac{3}{2}} \times H^{\frac{3}{2}}} \\
= & \, \left\langle \partial_t\, J^\beta_\delta \vw(s,\cdot) , \U \right\rangle_{H^{-\frac{3}{2}} \times H^{\frac{3}{2}}},
\end{split}
\end{equation*}
and integrating on the interval of time $[0,t]$ we obtain the desired identity  (\ref{Identity-1}). 

\medskip

We substitute the identity (\ref{Identity-1}) in the second term of the identity (\ref{iden1}) to get 
\[ \Vert  J^{\frac{\beta}{2}}_{\delta} \vw(t,\cdot) \Vert^{2}_{L^2}= \Vert  J^{\frac{\beta}{2}}_{\delta}\vu(t,\cdot) \Vert^{2}_{L^2} -2 \left(  J^{\frac{\beta}{2}}_{\delta} \vu_0,   J^{\frac{\beta}{2}}_{\delta}\U \right)_{L^2} - 2 \int_{0}^{t}\left\langle \partial_t \, J^\beta_\delta \vw(s,\cdot) , \U \right\rangle_{H^{-\frac{3}{2}} \times H^{\frac{3}{2}}} \, ds  + \Vert  J^{\frac{\beta}{2}}_{\delta} \U \Vert_{L^2}. \]

Here, we substitute the term $\Vert  J^{\frac{\beta}{2}}_{\delta}\vu(t,\cdot) \Vert^{2}_{L^2}$  with the right-hand side of the energy  estimate given in Definition \ref{Def-Leray} to obtain: 
\begin{equation*}
\begin{split}
\Vert  J^{\frac{\beta}{2}}_{\delta} \vw(t,\cdot) \Vert^{2}_{L^2} \leq & \,  \Vert  J^{\frac{\beta}{2}}_{\delta}\vu_0 \Vert^{2}_{L^2} -  \, 2 \nu \, \int_{0}^{t}  \left\Vert  (-\Delta)^{\frac{\alpha}{4}} J^{\frac{\beta}{2}}_{\delta}\, \vu(s,\cdot)\right\Vert^{2}_{L^2}   ds \\
& + 2 \int_{0}^{t} \left( J^{\frac{\beta}{2}}_{\delta}\fe, J^{\frac{\beta}{2}}_{\delta} \vu(s,\cdot) \right)_{L^2} \, ds -2 \gamma \int_{0}^{t}\Vert J^{\frac{\beta}{2}}_{\delta} \vu(s,\cdot)\Vert^{2}_{L^2} ds \\
& -2 \left(  J^{\frac{\beta}{2}}_{\delta} \vu_0,   J^{\frac{\beta}{2}}_{\delta}\U \right)_{L^2} - 2 \int_{0}^{t}\left\langle \partial_t \, J^\beta_\delta \vw(s,\cdot) , \U \right\rangle_{H^{-\frac{3}{2}} \times H^{\frac{3}{2}}} \, ds  + \Vert  J^{\frac{\beta}{2}}_{\delta} \U \Vert_{L^2}. 
\end{split}
\end{equation*} 
Rearranging  terms we get 
\begin{equation}\label{estim2}
\begin{split}
\Vert  J^{\frac{\beta}{2}}_{\delta} \vw(t,\cdot) \Vert^{2}_{L^2} \leq & \, \,  \underbrace{\Vert  J^{\frac{\beta}{2}}_{\delta}\vu_0 \Vert^{2}_{L^2}  -2 \left(  J^{\frac{\beta}{2}}_{\delta} \vu_0,  J^{\frac{\beta}{2}}_{\delta} \U \right)_{L^2} +  \Vert  J^{\frac{\beta}{2}}_{\delta} \U \Vert_{L^2}}_{(A)} \\
&   +   \underbrace{2\, \int_{0}^{t}  \left(  J^{\frac{\beta}{2}}_{\delta} \fe,  J^{\frac{\beta}{2}}_{\delta} \vu(s,\cdot) \right)_{L^2}\, ds}_{(B)} \\
&-   \underbrace{2\, \int_{0}^{t}\left\langle \partial_t \, J^\beta_\delta \vw(s,\cdot) , \U \right\rangle_{H^{-\frac{3}{2}} \times H^{\frac{3}{2}}} \, ds}_{(C)} \\ 
& - 2 \nu \, \int_{0}^{t}  \left\Vert  (-\Delta)^{\frac{\alpha}{4}} J^{\frac{\beta}{2}}_\delta\,  \vu(s,\cdot)\right\Vert^{2}_{L^2}   ds  -  2\gamma \int_{0}^{t}\Vert  J^{\frac{\beta}{2}}_{\delta} \vu(s,\cdot)\Vert^{2}_{L^2} ds,
\end{split}
\end{equation}
where we must study the expressions $(A)$, $(B)$ and $(C)$ separately. The term $(A)$ is easy to handle  and we have 
\begin{equation}\label{(A)}
(A)= \Vert  J^{\frac{\beta}{2}}_{\delta} (\vu_0 - \U) \Vert^{2}_{L^2}.
\end{equation}

Next, in order to study   term $(B)$, we remark  that by  equation (\ref{Stationary}) we can write 
\[ J^\beta_\delta \fe  = \nu (- \Delta)^{\frac{\alpha}{2}} \, J^\beta_\delta \U + \P\, \text{div}(\U \otimes \U) + \gamma \, J^\beta_\delta \U,\]
and we have 
\begin{equation}\label{(B)}
\begin{split}
(B) = &\,\,  2\, \int_{0}^{t}  \left(  J^\beta_{\delta}  \fe, \vu(s,\cdot) \right)_{L^2}\, ds  \\
= &\,\, 2\, \int_{0}^{t}  \left(  \nu (- \Delta)^{\frac{\alpha}{2}} \, J^\beta_\delta \U + \P\, \text{div}(\U \otimes \U) + \gamma \, J^\beta_\delta \U , \vu(s,\cdot) \right)_{L^2}\, ds \\
= & \,\, 2\nu \int_{0}^{2} \left( (-\Delta)^{\frac{\alpha}{4}} J^{\frac{\beta}{2}}_\delta \U, (-\Delta)^{\frac{\alpha}{4}} J^{\frac{\beta}{2}}_\delta \vu(s,\cdot) \right)_{L^2}\, ds + 2 \gamma \int_{0}^{t} \left( J^{\frac{\beta}{2}}_{\delta} \U, J^{\frac{\beta}{2}}_{\delta}  \vu(s,\cdot) \right)_{L^2}\, ds. \\
& + 2 \int_{0}^{t} \left(  \text{div}(\U \otimes \U)  , \vu(s,\cdot) \right)_{L^2}\, ds. 
\end{split}
\end{equation}

Finally, to study the term  $(C)$, we remark now that by the equation  (\ref{eq-w}) we can write 
\[ \partial_t \, J^{\beta}_\delta \vw= - \nu (-\Delta)^{\frac{\alpha}{2}} J^\beta_\delta \vw - \P \left( \text{div}(\vw \otimes \vw)+ \text{div}(\vw \otimes \U)+ \text{div}(\U\otimes \vw) \right) - \gamma J^\beta_\delta \, \vw.\] 
Therefore, we obtain
\begin{equation*}
\begin{split}
(C) =& \,\, -2 \nu \int_{0}^{t} \left\langle  (-\Delta)^{\frac{\alpha}{2}} J^\beta \vw(s,\cdot) , \U \right\rangle_{H^{-\frac{3}{2}} \times H^{\frac{3}{2}}}\, ds  - 2  \int_{0}^{t} \left\langle  \text{div}(\vw \otimes \vw) , \U \right\rangle_{H^{-\frac{3}{2}} \times H^{\frac{3}{2}}}\, ds \\
& - 2  \int_{0}^{t} \left\langle  \text{div}(\vw\otimes \U) , \U \right\rangle_{H^{-\frac{3}{2}} \times H^{\frac{3}{2}}}\, ds  - 2  \int_{0}^{t} \left\langle  \text{div}(\U\otimes \vw) , \U \right\rangle_{H^{-\frac{3}{2}} \times H^{\frac{3}{2}}}\, ds\\
&-2 \gamma \int_{0}^{t}\left( J^{\frac{\beta}{2}}_{\delta}  \vw(s,\cdot),  J^{\frac{\beta}{2}}_{\delta} \U \right)_{L^2}\, ds.
\end{split}
\end{equation*}

In this identity, for the first term on the right-hand side we directly  write 
\[ 2 \nu \int_{0}^{t} \left\langle  (-\Delta)^{\frac{\alpha}{2}} J^\beta_\delta \vw , \U \right\rangle_{H^{-\frac{3}{2}} \times H^{\frac{3}{2}}}\, ds= -2 \nu \int_{0}^{t} \left( (-\Delta)^{\frac{\alpha}{4}} J^{\frac{\beta}{2}}_\delta \vw(s,\cdot) , (-\Delta)^{\frac{\alpha}{4}} J^{\frac{\beta}{2}}_\delta\,\U \right)_{L^2} \, ds. \]
Thereafter, to study the second and the third term on the right-hand side, we have the following remarks.  On the one hand, by Lemma \ref{LemThec1} we have $\U \in H^{\alpha+\frac{\beta}{2}}(\Rt)$, and as we have $\frac{3}{2}<\alpha +\frac{\beta}{2}$  by the Sobolev embeddings we get $\U \in L^{\infty}(\Rt)$. Then, as $\vw \in (L^{2}_{t})_{loc} L^{2}_{x}$ (recall that $\vw$ belongs to the energy space) we have $\vw \otimes \U \in (L^{2}_{t})L^{2}_{x}$ and consequently $\text{div}(\vw \otimes \U) \in (L^{2}_{t})_{loc}H^{-1}_{x}$.  On the other hand, we recall that $ \vw  \in  (L^{2}_{t})_{loc}H^{\frac{\alpha+\beta}{2}}_{x}$, and moreover,  as we have $2 \leq \alpha +\beta$  we get $ \vw  \in  (L^{2}_{t})_{loc}H^{1}_{x}$.  With these remarks and integrating by parts we are able to write:
\[ - 2  \int_{0}^{t} \left\langle  \text{div}(\vw \otimes \vw) , \U \right\rangle_{H^{-\frac{3}{2}} \times H^{\frac{3}{2}}}\, ds =  2  \int_{0}^{t} \left\langle   (\vw \cdot \vec{\nabla})  \U, \vw   \right\rangle_{H^{-1} \times H^{1}}\, ds, \]
and 
\[ - 2  \int_{0}^{t} \left\langle  \text{div}(\vw \otimes \U) , \U \right\rangle_{H^{-\frac{3}{2}} \times H^{\frac{3}{2}}}\, ds =  2  \int_{0}^{t} \left\langle   (\U \cdot \vec{\nabla})  \U, \vw   \right\rangle_{H^{-1} \times H^{1}}\, ds. \] 
Finally, as $\text{div}(\U)=\text{div}(\vw)=0$, in the fourth term on the right-hand side we have
\[ - 2  \int_{0}^{t} \left\langle  \text{div}(\U\otimes \vw) , \U \right\rangle_{H^{-\frac{3}{2}} \times H^{\frac{3}{2}}}\, ds=0. \]

In this fashion,  term $(C)$ writes down as follows:
\begin{equation}\label{(C)}
\begin{split}
(C) = & \,\, -2 \nu \int_{0}^{t} \left( (-\Delta)^{\frac{\alpha}{4}} J^{\frac{\beta}{2}}_{\delta} \vw(s,\cdot) , (-\Delta)^{\frac{\alpha}{4}} J^{\frac{\beta}{2}}\,\U \right)_{L^2} \, ds  - 2 \gamma \int_{0}^{t}\left( J^{\frac{\beta}{2}}_{\delta}  \vw(s,\cdot), J^{\frac{\beta}{2}}_{\delta} \U \right)_{L^2}\, ds\\
&+  2  \int_{0}^{t} \left\langle   (\vw \cdot \vec{\nabla})  \U, \vw   \right\rangle_{H^{-1} \times H^{1}}\, ds + 2  \int_{0}^{t} \left\langle   (\U \cdot \vec{\nabla})  \U, \vw   \right\rangle_{H^{-1} \times H^{1}}\, ds.
\end{split}
\end{equation}

With identities (\ref{(A)}), (\ref{(B)}) and (\ref{(C)}) at hand,  we get back to the inequality (\ref{estim2}) to write 
\begin{equation}\label{estim4}
\begin{split}
\Vert J^{\frac{\beta}{2}}_{\delta}  \vw(t,\cdot) \Vert^{2}_{H^{L^2}} \leq & \, \,  \Vert J^{\frac{\beta}{2}}_{\delta} ( \vu_0 - \U) \Vert^{2}_{L^2} \\
 &+  \underbrace{2\nu \int_{0}^{2} \left( (-\Delta)^{\frac{\alpha}{4}} J^{\frac{\beta}{2}}_\delta \U, (-\Delta)^{\frac{\alpha}{4}} J^{\frac{\beta}{2}}_\delta \vu(s,\cdot) \right)_{L^2}\, ds}_{(A_1)} \, + \, \underbrace{2 \gamma \int_{0}^{t} \left( J^{\frac{\beta}{2}}_{\delta} \U, J^{\frac{\beta}{2}}_{\delta}  \vu(s,\cdot) \right)_{L^2}\, ds}_{(B_1)} \\
& + \underbrace{2 \int_{0}^{t} \left(  \text{div}(\U \otimes \U)  , \vu(s,\cdot) \right)_{L^2}\, ds}_{(C_1)}\\
&\underbrace{2 \nu \int_{0}^{t} \left( (-\Delta)^{\frac{\alpha}{4}} J^{\frac{\beta}{2}}_{\delta} \vw(s,\cdot) , (-\Delta)^{\frac{\alpha}{4}} J^{\frac{\beta}{2}}\,\U \right)_{L^2} \, ds}_{(A_2)} \, + \, \underbrace{2 \gamma \int_{0}^{t}\left( J^{\frac{\beta}{2}}_{\delta}  \vw(s,\cdot), J^{\frac{\beta}{2}}_{\delta} \U \right)_{L^2}\, ds}_{(B_2)}\\
&- 2  \int_{0}^{t} \left\langle   (\vw \cdot \vec{\nabla})  \U, \vw   \right\rangle_{H^{-1} \times H^{1}}\, ds  - \underbrace{2  \int_{0}^{t} \left\langle   (\U \cdot \vec{\nabla})  \U, \vw   \right\rangle_{H^{-1} \times H^{1}}\, ds}_{(C_2)} \\
& - \underbrace{2 \nu \, \int_{0}^{t}  \left\Vert  (-\Delta)^{\frac{\alpha}{4}} J^{\frac{\beta}{2}}_\delta\,  \vu(s,\cdot)\right\Vert^{2}_{L^2}ds}_{(A_3)}  \,  - \, \underbrace{2  \gamma \int_{0}^{t}\Vert  J^{\frac{\beta}{2}}_{\delta} \vu(s,\cdot)\Vert^{2}_{L^2} ds}_{(B_3)}.
\end{split}
\end{equation}

Here, we have the following  remarks. First, we can prove that 
\begin{equation}\label{Id1}
\sum_{i=1}^{3} (A_i)= -2\nu \, \int_{0}^{t}\Vert (-\Delta)^{\frac{\alpha}{4}} J^{\frac{\beta}{2}}_{\delta}  \vw(s,\cdot) \Vert^{2}_{L^2}\, ds.
\end{equation}
Indeed,  as $\vw=\vu - \U$ we  write:
\begin{equation*}
\begin{split}
&-2 \nu \left\Vert  (-\Delta)^{\frac{\alpha}{4}} J^{\frac{\beta}{2}}_{\delta} \, \vu \right\Vert^{2}_{L^2} + 2 \nu  \left( (-\Delta)^{\frac{\alpha}{4}} J^{\frac{\beta}{2}}_{\delta}  \U, (-\Delta)^{\frac{\alpha}{4}} J^{\frac{\beta}{2}}_{\delta}  \vu  \right)_{L^2} + 2\nu \left( (-\Delta)^{\frac{\alpha}{4}} J^{\frac{\beta}{2}}_{\delta} \vw , (-\Delta)^{\frac{\alpha}{4}} J^{\frac{\beta}{2}}_{\delta} \,\U \right)_{L^2}\\
=& -2 \nu \left\Vert  (-\Delta)^{\frac{\alpha}{4}} J^{\frac{\beta}{2}}_{\delta} \, \vu \right\Vert^{2}_{L^2}  + 4 \nu \left( (-\Delta)^{\frac{\alpha}{4}} J^{\frac{\beta}{2}}_{\delta}  \U, (-\Delta)^{\frac{\alpha}{4}} J^{\frac{\beta}{2}}_{\delta}  \vu  \right)_{L^2} + 2 \nu \left\Vert  (-\Delta)^{\frac{\alpha}{4}}J^{\frac{\beta}{2}}_{\delta} \, \U \right\Vert^{2}_{L^2}= -2\nu \left\Vert  (-\Delta)^{\frac{\alpha}{4}} J^{\frac{\beta}{2}}_{\delta} \, \vw \right\Vert^{2}_{L^2}.  
\end{split}
\end{equation*} 
The same arguments yield
\begin{equation}\label{Id2}
\sum_{i=1}^{3}(B_i) = - 2\gamma \int_{0}^{t} \Vert  J^{\frac{\beta}{2}}_{\delta} \vw(s,\cdot) \Vert^{2}_{L^2}\, ds.
\end{equation}
Moreover,  always by the identity   $\vw=\vu - \U$ and the divergence-free of $\vu$ and $\U$   we obtain
\begin{equation}\label{Id3}
\sum_{i=1}^{2}(C_i)=0.
\end{equation}
Finally, for a generic $0<C$ and for the quantities $\a,\b,\c$ defined in (\ref{Notation}),  the following estimate holds:
\begin{equation}\label{Id4}
\left|  2  \int_{0}^{t} \left\langle   (\vw \cdot \vec{\nabla})  \U, \vw   \right\rangle_{H^{-1} \times H^{1}}\, ds \right| \leq  C\, \frac{\b}{\a \c} \|\fe \|_{H^\frac{\beta}{2}} \int_{0}^{t} \Vert \vw(s,\cdot) \Vert^{2}_{H^{\frac{\alpha+\beta}{2}}} \, ds. 
\end{equation}
Indeed, by the  H\"older inequalities, the Hardy-Littlewood-Sobolev inequalities,  moreover, as $2 \leq \alpha+\beta$  we can write: 
\begin{equation*}
\begin{split}
& \left\vert  2  \int_{0}^{t} \left\langle   (\vw \cdot \vec{\nabla})  \U, \vw   \right\rangle_{H^{-1} \times H^{1}}\, ds  \right\vert \leq \,\, 2 \int_{0}^{t} \left\vert  \left\langle   (\vw \cdot \vec{\nabla})  \U, \vw   \right\rangle_{H^{-1} \times H^{1}}  \right\vert\, ds \\
\leq & \,\, 2 \int_{0}^{t} \left\Vert  (\vw(s,\cdot) \cdot \vec{\nabla})  \U  \right\Vert_{H^{-1}} \, \Vert \vw(s,\cdot) \Vert_{H^1}\, ds \leq   \,\, 2 \int_{0}^{t} \left\Vert \text{div} (\U \otimes \vw(s,\cdot))  \right\Vert_{\dot{H}^{-1}}\,   \Vert \vw(s,\cdot) \Vert_{H^1}\, ds\\
\leq & \,\, 2 \int_{0}^{t} \Vert \U \otimes \vw(s,\cdot) \Vert_{L^2}\, \Vert \vw(s,\cdot) \Vert_{H^1}\, ds \leq   \,\, 2 \int_{0}^{t} \Vert \U \Vert_{L^3} \Vert \vw(s,\cdot) \Vert_{L^6} \Vert \vw(s,\cdot) \Vert_{H^{1}}\, ds \\
\leq & \,\, C\, \Vert \U \Vert_{L^3} \, \int_{0}^{t} \Vert \vw(s,\cdot) \Vert_{\dot{H}^{1}}\, \Vert \vw(s,\cdot) \Vert_{H^{1}}\, ds \leq  \,\, C\, \Vert \U \Vert_{L^3} \int_{0}^{t} \Vert \vw(s,\cdot) \Vert^{2}_{H^{\frac{\alpha+\beta}{2}}} \, ds.
\end{split}
\end{equation*}
Then, by the interpolation inequalities,  the Hardy-Littlewood-Sobolev inequalities, the fact that $2\leq \alpha+\beta$ and the energy estimate obtained in Theorem \ref{Th3},  we are able to write: 
\[ C\, \Vert \U \Vert_{L^3} \leq C\, \Vert \U \Vert^{1/2}_{L^2}\, \Vert \U \Vert^{1/2}_{L^6}\leq C\, \Vert \U \Vert^{1/2}_{L^2}\, \Vert \U \Vert^{1/2}_{\dot{H}^{1}} \leq C \, \Vert \U \Vert_{H^1}\leq C \, \Vert \U \Vert_{H^{\frac{\alpha+\beta}{2}}} \leq C \frac{\b}{\a \c} \| \fe \|_{H^{\frac{\beta}{2}}}.\]

We substitute (\ref{Id1}), (\ref{Id2}), (\ref{Id3}) and (\ref{Id4}) in estimate (\ref{estim4}) to obtain (\ref{Key-Estim}). Proposition \ref{Prop:Key-Estimate} is proven.  \finpv 

\medskip 

\medskip

With this energy estimate at hand, we are able to conclude each point stated in Theorem \ref{Th4}. 

\begin{enumerate}
	\item   Recalling our notation $J^{\frac{\alpha}{2}}_{\gamma}=\gamma I_d + \nu(-\Delta)^{\frac{\alpha}{4}}$, we can write 
	\begin{equation}\label{Estim1}
	-2\nu \, \int_{0}^{t}\Vert (-\Delta)^{\frac{\alpha}{2}} J^{\frac{\beta}{2}}_{\delta}  \vw(s,\cdot) \Vert^{2}_{L^2}\, ds  - 2\gamma \int_{0}^{t} \Vert  J^{\frac{\beta}{2}}_{\delta} \vw(s,\cdot) \Vert^{2}_{L^2}\, ds= -2 \int_{0}^{t} \| J^{\frac{\alpha}{2}}_{\gamma} J^{\frac{\beta}{2}}_{\delta} \vw(s,\cdot) \|^{2}_{L^2} ds.
	\end{equation}
Moreover, 	by the identities (\ref{Identities}) and by the lower bounds in (\ref{Bounds}) we have
\begin{equation}\label{Estim2}
\begin{split}
-2 \int_{0}^{t} \| J^{\frac{\alpha}{2}}_{\gamma} J^{\frac{\beta}{2}}_{\delta} \vw(s,\cdot) \|^{2}_{L^2} ds = &\, -2 \int_{0}^{t} \| D(m^{1/2}_{1})D(m^{1/2}_{2}) (I_d-\Delta)^{\frac{\alpha+\beta}{4}} \vw(s,\cdot) \|^{2}_{L^2} ds\\
= &\,  -2 \int_{0}^{t} \| D(m^{1/2}_1)D(m^{1/2}_2) \vw(s,\cdot) \|^{2}_{H^{\frac{\alpha+\beta}{2}}} ds\\
\leq & \, -2 \sqrt{\a \c} \int_{0}^{t} \|\vw(s,\cdot) \|^{2}_{H^{\frac{\alpha+\beta}{2}}} ds. 
\end{split}
\end{equation}

Getting back to the estimate (\ref{Key-Estim}) we obtain
\begin{equation*}
\Vert J^{\frac{\beta}{2}}_{\delta}  \vw(t,\cdot) \Vert^{2}_{H^{L^2}} \leq   \Vert J^{\frac{\beta}{2}}_{\delta} ( \vu_0 - \U) \Vert^{2}_{L^2}  + \left(-2\sqrt{\a \c} + C \frac{\b}{\a \c} \| \fe \|_{H^{\frac{\beta}{2}}} \right)  \int_{0}^{t} \Vert \vw(s,\cdot) \Vert^{2}_{H^{\frac{\alpha+\beta}{2}}} \, ds.
\end{equation*}
We assume (\ref{Assumption1}) hence we have $ -2\sqrt{\a \c} + C \frac{\b}{\a \c} \| \fe \|_{H^{\frac{\beta}{2}}} \leq 0$. We thus get the estimate
\[ \Vert J^{\frac{\beta}{2}}_{\delta}  \vw(t,\cdot) \Vert^{2}_{H^{L^2}} \leq   \Vert J^{\frac{\beta}{2}}_{\delta} ( \vu_0 - \U) \Vert^{2}_{L^2},\]
from which directly follows the orbital stability of stationary solutions. 

\item   We get back to the estimate (\ref{Key-Estim}), hence we write
 \begin{equation}\label{Key-Estim-2}
 \begin{split}
 \Vert J^{\frac{\beta}{2}}_{\delta}  \vw(t,\cdot) \Vert^{2}_{H^{L^2}} \leq & \, \,  \Vert J^{\frac{\beta}{2}}_{\delta} ( \vu_0 - \U) \Vert^{2}_{L^2}   -\nu \, \int_{0}^{t}\Vert (-\Delta)^{\frac{\alpha}{4}} J^{\frac{\beta}{2}}_{\delta}  \vw(s,\cdot) \Vert^{2}_{L^2}\, ds - \gamma \int_{0}^{t} \Vert  J^{\frac{\beta}{2}}_{\delta} \vw(s,\cdot) \Vert^{2}_{L^2}\, ds \\
 &+ C \frac{\b}{\a \c} \| \fe \|_{H^{\frac{\beta}{2}}} \int_{0}^{t} \Vert \vw(s,\cdot) \Vert^{2}_{H^{\frac{\alpha+\beta}{2}}} \, ds - \gamma \int_{0}^{t} \Vert  J^{\frac{\beta}{2}}_{\delta} \vw(s,\cdot) \Vert^{2}_{L^2}\, ds.
 \end{split}
 \end{equation}
 By (\ref{Estim1}) and (\ref{Estim2}) we have 
 \[ -\nu \, \int_{0}^{t}\Vert (-\Delta)^{\frac{\alpha}{4}} J^{\frac{\beta}{2}}_{\delta}  \vw(s,\cdot) \Vert^{2}_{L^2}\, ds - \gamma \int_{0}^{t} \Vert  J^{\frac{\beta}{2}}_{\delta} \vw(s,\cdot) \Vert^{2}_{L^2}\, ds  \leq - \sqrt{\a\c}  \int_{0}^{t} \|\vw(s,\cdot) \|^{2}_{H^{\frac{\alpha+\beta}{2}}} ds,\]
 and we get
\begin{equation*}
\begin{split}
\Vert J^{\frac{\beta}{2}}_{\delta}  \vw(t,\cdot) \Vert^{2}_{L^2} \leq &\,    \Vert J^{\frac{\beta}{2}}_{\delta} ( \vu_0 - \U) \Vert^{2}_{L^2}  + \left(-\sqrt{\a \b} + C \frac{\b}{\a \c} \| \fe \|_{H^{\frac{\beta}{2}}} \right)  \int_{0}^{t} \Vert \vw(s,\cdot) \Vert^{2}_{H^{\frac{\alpha+\beta}{2}}} \, ds \\
&\, - \gamma \int_{0}^{t} \Vert  J^{\frac{\beta}{2}}_{\delta} \vw(s,\cdot) \Vert^{2}_{L^2}\, ds. 
\end{split}
\end{equation*} 
Here we assume (\ref{Assumption2}) to obtain that $\sqrt{\a \b} + C \frac{\b}{\a \c} \| \fe \|_{H^{\frac{\beta}{2}}} \leq 0$. Then, by the Gr\"owall inequalities we obtain
\begin{equation}\label{Inequality-01}
\Vert J^{\frac{\beta}{2}}_{\delta}  \vw(t,\cdot) \Vert^{2}_{L^2} \leq e^{-\gamma t} \Vert J^{\frac{\beta}{2}}_{\delta} ( \vu_0 - \U) \Vert^{2}_{L^2} ,
\end{equation}
which yields the  estimate (\ref{Inequality}).  Remark that this estimate also yields the uniqueness of the stationary solution $\U$. Indeed, if  $\U_1, \U_2$ are  two solutions to equation (\ref{Stationary}) then we can set 
 $\vw=\U_1 - \U_2$. Observe that since  $\partial_t \U_1=0$ then $\U_1$ is also a solution to the evolution equation (\ref{Equation}) with initial datum $\vu_0 = \U_1$. Then,  for a time $0<t$ such that $e^{-\gamma t} < \frac{1}{2}$  by (\ref{Inequality-01}) we have $\Vert J^{\frac{\beta}{2}}_{\delta} (\U_1 - \U_2)\Vert^{2}_{L^2} \leq \frac{1}{2} \Vert J^{\frac{\beta}{2}}_{\delta} ( \U_1 - \U_2) \Vert^{2}_{L^2}$, hence $\U_1=\U_2$. 
\end{enumerate}	 
Theorem \ref{Th4} is proven. \finpv

\subsection{Proof of Corollary \ref{Corollary}}
We shall prove that $\{ \U\}$ is a strong global attractor in the sense of Definition \ref{Def-global-attractor} with $\bullet=s$. The first point of Definition \ref{Def-global-attractor}  is evident, while   by estimate   (\ref{Inequality})  we directly obtain that $\{ \U \}$ is an strong attracting  set for equation (\ref{Equation}) in the sense of Definition \ref{Attracting-set}. Then, by   uniqueness of the global attractor  we have $\mathcal{A}_w= \mathcal{A}_s = \{ \U \}$.  \finpv

\begin{appendices}
\section{Appendix}\label{Appendix-Fractal-dim}
In all this section, for the parameter $0<\delta$ fixed  we shall consider the space $H^{\frac{\beta}{2}}(\Rt)$ with the \emph{equivalent} inner product $\ds{(f, g)_{H^{\frac{\beta}{2}}_{\delta}}= \int_{\Rt} (1+\delta^2 |\xi|^2)^{\frac{\beta}{2}} \widehat{f}(\xi) \bar{\widehat{g}}(\xi)d\xi}$. Moreover, we shall denote by  $\| f \|_{H^{\frac{\beta}{2}}_{\delta}}$ its corresponding norm, where for the quantities $\a$ and $\b$ defined in (\ref{Notation}) we have 
\begin{equation}\label{Equivalence}
\a \| f \|_{H^{\frac{\beta}{2}}} \leq \| f \|_{H^{\frac{\beta}{2}}_{\delta}} \leq \b \| f \|_{H^{\frac{\beta}{2}}}.
\end{equation}

Let us recall that the fractal dimension $\text{dim}\left(\mathcal{A}_s\right)$ is commonly estimated by the so-called box-counting method. By the Hausdorff criterion, for every $0<\varepsilon$ the compact set $\mathcal{A}_s$ can be covered by a finite number of $\varepsilon-$balls in the space $H^{\frac{\beta}{2}}(\Rt)$. We denote by $N_{\varepsilon}(\mathcal{A}_s)$ the minimal number of such   $\varepsilon-$balls. Then, we have the following: 

\begin{Definition}\label{def-fractal-dim-attractor} The fractal (box-counting) dimension of the strong global attractor $\mathcal{A}_s$ is given by the quantity:
	\begin{equation}
	\text{dim}\left(\mathcal{A}_s\right)= \limsup_{\varepsilon \to 0^{+}} \frac{ \ln \left(N_{\varepsilon} (\mathcal{A}_s)\right)}{\ln\left( \frac{1}{\varepsilon}\right)}. 
	\end{equation}	
\end{Definition}	

Our next result reads as follows: 

\begin{Theoreme}\label{Prop-fractal-dim} Let $1\leq \alpha$ and $2 \leq \beta$. There exists a constant $ 0<\mathfrak{C}$ given in (\ref{C}),  depending on the quantities $\a,\b,\c$ defined in (\ref{Notation}) and the damping parameter $\gamma$, such that the following estimate holds: 
	\begin{equation}\label{Fractal-dimension}
	\text{dim}\left(\mathcal{A}_s\right) \leq \frac{2 \mathfrak{C}}{\gamma} \, \max\left( \Vert \fe \Vert^{2}_{H^{\frac{\beta}{2}}_{\delta}},  \Vert \fe  \Vert^{4}_{H^{\frac{\beta}{2}}_{\delta}} \right).
	\end{equation}	
\end{Theoreme} 	

Some comments are in order. We observe that the fractal dimension of $\mathcal{A}_s$ is essentially   controlled by the damping parameter $\gamma$ and  the size  of the external force in the $H^{\frac{\beta}{2}}_{\delta}-$ norm.  Precisely, for $\| \fe \|_{H^{\frac{\beta}{2}}_{\delta}}$ fixed, large values of $\gamma$ yield small values of 	$\text{dim}\left(\mathcal{A}_s\right)$. This type of control was also pointed out in \cite{Cao} and \cite{Ilyn} for some related models.

\medskip 

The assumptions  $1\leq \alpha$ and $2 \leq \beta$ are essentially required to  adapt the method used in previous works \cite{Constantin1,Ilyn,Temam} to the more general equation (\ref{Equation}), but we think that in further investigations this estimate could be improved to the less restrictive assumption $\frac{5}{2} \leq \alpha + \beta$, where  uniqueness of Leray-type solutions is known  and it is one  the key ideas to  derive the estimate (\ref{Fractal-dimension}). On the other hand,  this method cannot be applied to the case of  weak global attractors  (when $0<\alpha +\beta<\frac{5}{2}$) where uniqueness of Leray-type solutions is unknown. To the best of our knowledge, upper estimates of the fractal dimension for  weak global attractors  is matter of deeper  and far from obvious research.  

\medskip

We start by introducing some definition and notation  that we shall use in the sequel. The first definition concerns the following quasi-differential operator. Let $0\leq  t$ be a fixed time and let $\vu_0 \in \mathcal{A}_{s}$ be an initial datum. Moreover, let $u(t,\cdot)$ be the unique solution  of  equation (\ref{Equation}) arising from $\vu_0$ and given by Theorem \ref{Th1}. Thus, for  $u(t,\cdot)$ fixed, let $\vv \in (L^{\infty}_{loc}([0,+\infty[, H^{\frac{\beta}{2}}(\Rt)) \cap L^{2}_{loc}([0, +\infty[, H^{\frac{\alpha+\beta}{2}}(\Rt))$ be the solution of the following linearized problem: 

\begin{equation}\label{Equation-Lin}
\left\{ \begin{array}{ll}\vspace{2mm}
\partial_t \vv + \nu (-\Delta)^{\frac{\alpha}{2}} \vv + (I_d -\delta^2\Delta)^{-\frac{\beta}{2}}\,   \P \left((\vv \cdot \vec{\nabla}) \vu + (\vu \cdot \vec{\nabla})\vv \right)   = -\gamma \vv, \quad \text{div}(\vv)=0,  \\ 
\vv(0,\cdot)= \vv_0 \in H^{\frac{\beta}{2}}(\Rt),
\end{array} 
\right.
\end{equation}
where $\vv_0$ denotes an initial datum. As we assume $\frac{5}{2}\leq \alpha +\beta$, the existence and uniqueness of a solution $\vv \in (L^{\infty}_{t})_{loc} H^{\frac{\beta}{2}}_{x} \cap (L^{2}_{t})_{loc} H^{\frac{\alpha+\beta}{2}}_{x}$ essentially follows the ideas in the proof 
of Theorem \ref{Th1}, so we will omit this proof.

\begin{Definition}[Quasi-differential operator]\label{def-quasi-diff-op} 
	The quasi-differential operator $DS(t,\vu_0)$, depending on the time $0\leq t$ and the datum $\vu_0 \in \mathcal{A}_{s}$, is the linear and bounded operator $DS(t,\vu_0): H^{\frac{\beta}{2}}(\Rt) \to H^{\frac{\beta}{2}}(\Rt)$ defined as 
	\begin{equation*}
	DS(t,\vu_0) \vv_0= \vv(t,\cdot),    
	\end{equation*}
	where  $\vv(t,\cdot)$ is the unique solution of  equation (\ref{Equation-Lin}). 
\end{Definition}

Once we have defined this operator, our second definition is devoted to the notion of a semigroup uniformly quasi-differentiable. 

\begin{Definition}[Semigroup uniformly quasi-differentiable]\label{def-unif-quasi-diff} Let $0\leq t$ fixed and let $S(t)$ be the semi-group associated to equation (\ref{Equation}) and defined in (\ref{Seligroup}). We say that this semigroup is uniformly quasi-differentiable on the global attractor $\mathcal{A}_{s}\subset H^{\frac{\beta}{2}}(\Rt)$,  if for all $\vu_{0,1}, \vu_{0,2}\in \mathcal{A}_{s}$ we have 
	\begin{equation*}
	\left\Vert S(t)\vu_{0,2}-S(t)\vu_{0,1}-DS(t,\vu_{0,1})(\vu_{0,2}-\vu_{0,1}) \right\Vert_{H^{\frac{\beta}{2}}_\delta}\leq \mathfrak{o} \left( \Vert \vu_{0,2}-\vu_{0,1}\Vert_{H^{\frac{\beta}{2}}_\delta}\right),    
	\end{equation*}
	where  $DS(t,\vu_{0,1})$ is given in Definition \ref{def-quasi-diff-op} and the quantity $\mathfrak{o}(\cdot)$ verifies: $\ds{\lim_{h \to 0^{+}} \frac{\mathfrak{o}(h)}{h}=0}$. 
\end{Definition}

Finally, in our last definition, we  introduce the notion of the $n-$ global Lyapunov exponent, with  $n \in \mathbb{N}$. For this,  we shall need to precise   some notation. On the one hand,  we denote by $\mathcal{O}_n$  the set of all the   families  $(\vw_i)_{1\leq i \leq n}$ in the space $H^{\frac{\alpha+\beta}{2}}(\Rt)$, with $\text{div}(\vw_i)=0$ and which are orthonormal in the  space $H^{\frac{\beta}{2}}(\Rt)$ with the inner product $(\cdot , \cdot)_{H^{\frac{\beta}{2}}_{\delta}}$. On the other hand, we get back to  equation (\ref{Equation-Lin}) to write
\begin{equation*}
\partial_t \vv=  - \nu (-\Delta)^{\frac{\alpha}{2}}\vv - (I_d -\delta^2\Delta)^{-\frac{\beta}{2}}\, \P \left((\vv \cdot \vec{\nabla}) \vu + (\vu \cdot \vec{\nabla})\vv \right)  -\gamma \vv. 
\end{equation*} and then, from the right-hand  side of this identity, and for all $\vw \in H^{\frac{\beta}{2}}(\Rt)$,  we define now the linear operator
\begin{equation}\label{defi-op-L}
\mathcal{L}(t,\vu_0)\vw= - \nu (-\Delta)^{\frac{\alpha}{2}} \vw - (I_d -\delta^2\Delta)^{-\frac{\beta}{2}}\, \P \left((\vw \cdot \vec{\nabla}) \vu + (\vu \cdot \vec{\nabla})\vw \right)-\gamma \vw.    
\end{equation}

Once we have introduced the set $\mathcal{O}_n$ and the linear operator $\mathcal{L}(t,\vu_0)(\cdot)$ above, we have the following:
\begin{Definition}[$n-$ global Lyapunov exponent]\label{def-global-lyapunov-expo} Let $n\in \mathbb{N}$ fixed. We define the $n-$ global Lyapunov exponent $\ell (n)$ by the quantity: 
	\begin{equation*}
	\ell (n)= \limsup_{T\to+\infty} \left(  \sup_{\vu_0 \in \mathcal{A}_{s}}\, \,  \sup_{(\vw_i)_{1\leq i \leq n} \in \mathcal{O}_n} \left( \frac{1}{T}\int_{0}^{T}\sum_{i=1}^{n} \left( \mathcal{L}(t,\vu_0)\vw_i, \vw_i \right)_{H^{\frac{\beta}{2}}_\delta} dt \right)\right).   
	\end{equation*}
\end{Definition} 

We have now all the tools to state the following technical result, which allows us  to find  an upper bound of the fractal dimension for the attractor $\mathcal{A}_{s}$. For a proof of this result see \cite{Chepyzhov}.
\begin{Theoreme}[Upper bound of the fractal dimension]\label{Th-tec-upper-bound} Let $\text{dim}(\mathcal{A}_s)$ be the fractal box counting dimension of the global  attractor $\mathcal{A}_{s}$ given in Definition \ref{def-fractal-dim-attractor}. If the following statements hold: 
	\begin{enumerate}
		\item The semigroup $S(t)$   is uniformly quasi-differentiable on  $\mathcal{A}_{s}$ in the sense of Definition \ref{def-unif-quasi-diff}.
		\item The quasi-differential operator $DS(t,\vu_0)(\cdot)$, given in Definition \ref{def-quasi-diff-op}, depends continuously  on the initial datum $\vu_0 \in \mathcal{A}_{s}$. 
		\item There exists $1 \leq \kappa $, and there exist two constants $c_1,c_2>0$ such that  for all $n\in \mathbb{N}$  the $n-$ global Lyapunov exponent $\ell(n)$ given in Definition \ref{def-global-lyapunov-expo} verifies:
		\begin{equation}\label{estim-lyapunov-exp}
		\ell(n)\leq -c_1 \, n^{\kappa}+c_2.  \vspace{-5mm}  
		\end{equation} 
	\end{enumerate}
	Then, we have the following upper bound:  $\ds{dim\left(\mathcal{A}_{s}\right) \leq \left(\frac{c_2}{c_1}\right)^{1/\kappa}}$. 
\end{Theoreme}

\subsection{Proof of Theorem \ref{Prop-fractal-dim}}
We must verify all the points stated in Theorem \ref{Th-tec-upper-bound}. The first and the second point are classical to verify and they  are essentially proven in \cite{Babin}. So,  we will focus on the third point which is more delicate. 

\medskip

In order to estimate  the $n-$ global Lyapunov exponent $\ell(n)$ according to  (\ref{estim-lyapunov-exp}), we shall prove  the following technical estimates. First,  in the expression of the quantity $\ell(n)$ given in Definition \ref{def-global-lyapunov-expo}, we shall estimate the  term $\ds{\sum_{i=1}^{n} \left(\mathcal{L}(t,\vu_0) \vw_{i}, \vw_i \right)_{H^{\frac{\beta}{2}}_{\delta}}}$ as follows:
\begin{Proposition}\label{Prop:estim-lyapunov} Let $n\in \mathbb{N}$ fixed and let $(\vw_i)_{1\leq i \leq n} \in \mathcal{O}_n$. Moreover,  let $\mathcal{L}(t,\vu_0)(\cdot)$ be the linear operator given in (\ref{defi-op-L}). Then,   we have: 
	\begin{equation}\label{estim-tech}
	\sum_{i=1}^{n} \left( \mathcal{L}(t,\vu_0)\vw_i, \vw_i \right)_{H^{\frac{\beta}{2}}_{\delta}} \leq   - \frac{\gamma \a}{2} n +  \frac{2}{5}\frac{C_{LT}^{5/2}}{(\a\c)^{3/2}}\, \Vert \vec{\nabla} \otimes \vu(t,\cdot) \Vert^{5/2}_{L^{5/2}}, 
	\end{equation} 
	where $0<C_{LT}$ is a numerical constant given in (\ref{LT}),  and the quantities $\a$ and $\b$ given in (\ref{Notation}). 
\end{Proposition}	
\pv Recall that $\text{div}(\vu)=0$ and $\| \vw_{i}\|_{H^{\frac{\beta}{2}}_\delta}=1$,  then  we write 
\begin{equation}\label{Iden-Lyapunov}
\begin{split}
&\sum_{i=1}^{n}\left( \mathcal{L}(t, \vu_0) \vw_i, \vw_i \right)_{H^{\frac{\beta}{2}}_{\delta}}  \\
=&\, \sum_{i=1}^{n} \left( -\nu (-\Delta)^{\frac{\alpha}{2}} \vw_i - \gamma \vw_i , \vw_i\right)_{H^{\frac{\beta}{2}}_{\delta}}  - \sum_{i=1}^{n}  \left(  (I_d-\delta^2\Delta)^{-\frac{\beta}{2}}( (\vw_i \cdot \vec{\nabla}) \vu), \vw_i \right)_{H^{\frac{\beta}{2}}_{\delta}} \\
&\,  -  \sum_{i=1}^{n}  \left( (I_d-\delta^2\Delta)^{-\frac{\beta}{2}}( (\vu \cdot \vec{\nabla}) \vw_i), \vw_i \right)_{H^{\frac{\beta}{2}}_{\delta}} \\
=&\, \sum_{i=1}^{n}  \left( - \nu (-\Delta)^{\frac{\alpha}{2}} \vw_i  - \frac{\gamma}{2}\vw_i, \vw_i  \right)_{H^{\frac{\beta}{2}}_{\delta}} - \frac{\gamma}{2} \sum_{i=1}^{n}  \left( \vw_i , \vw_i  \right)_{H^{\frac{\beta}{2}}_{\delta}} -  \sum_{i=1}^{n} \left((\vw_i \cdot \vec{\nabla}) \vu), \vw_i \right)_{L^2}\\
&\, - \sum_{i=1}^{n}  \left((\vu \cdot \vec{\nabla}) \vw_i), \vw_i \right)_{L^2}\\
=&\, \underbrace{ - \sum_{i=1}^{n}   \left( \left( \frac{\gamma}{2} I_d  + \nu (-\Delta)^{\frac{\alpha}{2}} \right)  \vw_i, \vw_i   \right)_{H^{\frac{\beta}{2}}_{\delta}}}_{(A)}  -  \frac{\gamma}{2} n  + \underbrace{- \sum_{i=1}^{n}  \left(  (\vw_i \cdot \vec{\nabla}) \vu), \vw_i \right)_{L^2}}_{(B)}.
\end{split} 
\end{equation}
To estimate term $(A)$, recall that by (\ref{Notation-Operators}) and (\ref{Identities}) (with $\frac{\gamma}{2}$ instead of $\gamma$) we have  $ \frac{\gamma}{2} I_d  + \nu (-\Delta)^{\frac{\alpha}{2}}= J^{\alpha}_{\frac{\gamma}{2}}= D(m_2) (I_d - \Delta)^{\frac{\alpha}{2}}$.  Recall also that $J^{\frac{\beta}{2}}_{\delta} = (I_d - \delta^2 \Delta)^{\frac{\beta}{2}}=D(m_1)(I_d - \Delta)^{\frac{\beta}{2}}$.  Moreover,  by  the lower  bounds in (\ref{Bounds}) we have  
\begin{equation*}
\begin{split}
 (A)  = & \,  -  \sum_{i=1}^{n} \left(  D(m_2) (I_d - \Delta)^{\frac{\alpha}{2}}    \vw_i, \vw_i   \right)_{H^{\frac{\beta}{2}}_{\delta}}  -  \sum_{i=1}^{n} \left( D(m_1) D(m_2) (I_d - \Delta)^{\frac{\alpha+\beta}{2}}    \vw_i, \vw_i   \right)_{L^2} \\
= & \, - \sum_{i=1}^{n}   \Vert D(m^{1/2}_1) D(m^{1/2}_{2}) \vw_i \Vert^{2}_{H^{\frac{\alpha+\beta}{2}}} \leq - \a \c \sum_{i=1}^{n}   \Vert  \vw_i \Vert^{2}_{H^{\frac{\alpha+\beta}{2}}}.
\end{split}
\end{equation*}

To estimate term $(B)$,   following the same computations performed in \cite{Ilyn} (see the estimate $(3.5)$ in the page $16$) we write  
\begin{equation*}
(B) =  \, - \sum_{i=1}^{n} \int_{\Rt} \sum_{j,k=1}^{3} w_{i,k} (\partial_k u_j) w_{i,j}\, dx  \leq \, \int_{\Rt} \vert \vec{\nabla} \otimes \vu \vert \, \sum_{i=1}^{n} \vert \vw_i \vert^2\, dx. 
\end{equation*}
Then, by  H\'older inequalities, and  by the Lieb-Thirring inequality  \cite{Lieb} (see the estimate $(6)$, page $2$),  for the constant  
\begin{equation}\label{LT}
C_{LT}= \frac{3}{5^{5/3}} \left( 16 \pi^{3/2} \frac{\Gamma(7/2)}{\Gamma(5)}\right)^{2/3},
\end{equation}
we have
\begin{equation*}
\begin{split}
&\, \int_{\Rt} \vert \vec{\nabla} \otimes \vu \vert \, \sum_{i=1}^{n} \vert \vw_i \vert^2\, dx \leq  \, \Vert \vec{\nabla} \otimes \vu \Vert_{L^{5/2}}\, \left\Vert \sum_{i=1}^{n} \vert \vw_i \vert^2 \right\Vert_{L^{5/3}} \leq C_{LT}  \Vert \vec{\nabla} \otimes \vu \Vert_{L^{5/2}} \left( \sum_{i=1}^{n} \Vert \vw_i \Vert^{2}_{\dot{H}^{1}}\right)^{3/5}\\
\leq &\, \frac{C_{LT}}{(\a\c)^{3/5}} \Vert \vec{\nabla} \otimes \vu \Vert_{L^{5/2}}\, \left(\a\c\,  \sum_{i=1}^{n}  \Vert \vw_i \Vert^{2}_{\dot{H}^{1}}\right)^{3/5} \leq \frac{2}{5}\frac{C_{LT}^{5/2}}{(\a\c)^{3/2}}\, \Vert \vec{\nabla} \otimes \vu \Vert^{5/2}_{L^{5/2}}+ \frac{3}{5} \a\c\, \sum_{i=1}^{n} \Vert \vw_i \Vert^{2}_{\dot{H}^{1}} \\
\leq & \, \frac{2}{5}\frac{C_{LT}^{5/2}}{(\a\c)^{3/2}}\, \Vert \vec{\nabla} \otimes \vu \Vert^{5/2}_{L^{5/2}}+ \frac{3}{5} \a\c\, \sum_{i=1}^{n} \Vert \vw_i \Vert^{2}_{H^{\frac{\alpha+\beta}{2}}}. 
\end{split} 
\end{equation*}
With these estimates at hand, we get back to the identity (\ref{Iden-Lyapunov}) where we can write
\begin{equation*}
\begin{split}
\sum_{i=1}^{n} \left( \mathcal{L}(t, \vu_0) \vw_i, \vw_i  \right)_{H^{\frac{\beta}{2}}_{\delta}} \leq &\, - \a\c \, \sum_{i=1}^{n} \Vert \vw_i \Vert^{2}_{H^{\frac{\alpha+\beta}{2}}} - \frac{\gamma }{2}n + \frac{2}{5}\frac{C_{LT}^{5/2}}{A^{3/2}}\, \Vert \vec{\nabla} \otimes \vu \Vert^{5/2}_{L^{5/2}}+ \frac{3}{5} \a\c\, \sum_{i=1}^{n} \Vert \vw_i \Vert^{2}_{H^{\frac{\alpha+\beta}{2}}}\\
\leq &\, - \frac{\gamma}{2} n +  \frac{2}{5}\frac{C_{LT}^{5/2}}{(\a\c)^{3/2}}\, \Vert \vec{\nabla} \otimes \vu \Vert^{5/2}_{L^{5/2}}. 
\end{split}
\end{equation*}
Proposition \ref{Prop:estim-lyapunov} is proven.  \finpv

\medskip

For $0<T$, we take the time-average $\frac{1}{T} \int_{0}^{T}(\cdot) dt$ in each term of the inequality (\ref{estim-tech}) to obtain 
\begin{equation*}
\frac{1}{T} \int_{0}^{T} \sum_{i=1}^{n} \left( \mathcal{L}(t, \vu_0) \vw_i, \vw_i\right)_{H^{\frac{\beta}{2}}_{\delta}} dt \leq - \frac{\gamma}{2} n +   \frac{2}{5}\frac{C_{LT}^{5/2}}{(\a\c)^{3/2}} \, \frac{1}{T} \int_{0}^{T} \Vert \vec{\nabla} \otimes \vu (t,\cdot) \Vert^{5/2}_{L^{5/2}}\, dt,  
\end{equation*}
where we must estimate the last expression on the right-hand side. By  interpolation inequalities, the Hardy-Littlewood-Sobolev inequalities, and recalling that we have $1\leq \alpha$ and $2\leq \beta$ (hence $1\leq \frac{\beta}{2}$ and $\frac{3}{2} \leq \frac{\alpha+\beta}{2}$)  we can write
\begin{equation}\label{Estim-0} 
\begin{split}
&\,\frac{1}{T}\int_{0}^{T}  \Vert \vec{\nabla} \otimes \vu (t,\cdot) \Vert^{5/2}_{L^{5/2}}\, dt \leq \,  \frac{1}{T}\int_{0}^{T}  \Vert \vec{\nabla} \otimes \vu (t,\cdot) \Vert_{L^{2}}\, \Vert \vec{\nabla} \otimes \vu (t,\cdot) \Vert^{3/2}_{L^{3}} dt \\
\leq & \,  \frac{1}{T}\int_{0}^{T}  \Vert  \vu (t,\cdot) \Vert_{\dot{H}^{1}}\, \Vert \vec{\nabla} \otimes \vu (t,\cdot) \Vert^{3/2}_{\dot{H}^{1/2}} dt
\leq \,  \frac{1}{T}\int_{0}^{T}  \Vert  \vu (t,\cdot) \Vert_{\dot{H}^{1}}\, \Vert \vu (t,\cdot) \Vert^{3/2}_{\dot{H}^{3/2}} dt \\
\leq &\, \frac{1}{4} \frac{1}{T}\int_{0}^{T} \Vert \vu(t,\cdot) \Vert^{4}_{\dot{H}^1}\,  dt+ \frac{3}{4} \,  \frac{1}{T}\int_{0}^{T} \Vert \vu(t,\cdot) \Vert^{2}_{\dot{H}^{3/2}}\,  dt \\
\leq &\, \frac{1}{4} \frac{1}{T}\int_{0}^{T} \Vert \vu(t,\cdot) \Vert^{4}_{H^{\frac{\beta}{2}}}\,  dt+ \frac{3}{4} \,  \frac{1}{T}\int_{0}^{T} \Vert \vu(t,\cdot) \Vert^{2}_{H^{\frac{\alpha+\beta}{2}}}\,  dt.
\end{split}
\end{equation}
In order to estimate the first term on the right-hand side, by (\ref{Estim-absorbing})  and (\ref{Equivalence}) we get  
\begin{equation}\label{Estim-1}
\begin{split}
\frac{1}{4} \frac{1}{T}\int_{0}^{T} \Vert \vu(t,\cdot) \Vert^{4}_{H^{\frac{\beta}{2}}}\,  dt \leq &\, \, \frac{1}{8 \gamma T} \Vert \vu_0 \Vert^{4}_{H^{\frac{\beta}{2}}} \left(1-e^{-2\gamma\,T} \right) + \frac{\b^4}{4\a^4 \gamma^4} \, \Vert \fe \Vert^{4}_{H^{\frac{\beta}{2}}} \leq  \frac{1}{8 \gamma T} \Vert \vu_0 \Vert^{4}_{H^{\frac{\beta}{2}}} + \frac{\b^4}{4 \a^4\gamma^4} \, \Vert \fe \Vert^{4}_{H^{\frac{\beta}{2}}} \\
\leq & \, \frac{1}{8 \gamma T} \Vert \vu_0 \Vert^{4}_{H^{\frac{\beta}{2}}} + \frac{\b^4}{4 \a^8\gamma^4} \, \Vert \fe \Vert^{4}_{H^{\frac{\beta}{2}}_\delta}.  
\end{split}
\end{equation} 
For the second term on the right-hand side,  by the second point of Proposition \ref{Prop1}  with $t=0$, recalling that $\fe$ does not depend on the time variable, and using again (\ref{Equivalence}) we have
\begin{equation}\label{Estim-2}
\begin{split}
\frac{3}{4} \frac{1}{T}\int_{0}^{T} \Vert \vu(t,\cdot) \Vert^{2}_{H^{\frac{\alpha+\beta}{2}}} dt \leq \frac{3 }{4\ \c \, T} \Vert \vu_0\Vert^{2}_{H^{\frac{\beta}{2}}}+ \frac{3 \b^2}{4\, \a^2\c^2} \Vert \fe \Vert^{2}_{H^{\frac{\beta}{2}}} \leq \frac{3 }{4\ \c \, T} \Vert \vu_0\Vert^{2}_{H^{\frac{\beta}{2}}}+ \frac{3 \b^2}{4\, \a^4\c^2} \Vert \fe \Vert^{2}_{H^{\frac{\beta}{2}}_\delta}. 
\end{split}
\end{equation}
With the estimates (\ref{Estim-1}) and (\ref{Estim-2}) at hand, we get back to the estimate (\ref{Estim-0})  to obtain 
\begin{equation*}
\frac{1}{T}\int_{0}^{T}  \Vert \vec{\nabla} \otimes \vu (t,\cdot) \Vert^{5/2}_{L^{5/2}}\, dt \leq  \frac{1}{T} \left( \frac{1}{8\gamma} \Vert \vu_0\Vert^{4}_{H^{\frac{\beta}{2}}}+\frac{3}{4 \c} \Vert \vu_0 \Vert_{H^{\frac{\beta}{2}}}\right)+ \left( \frac{\b^4}{4 \a^8\gamma^4}+\frac{3\b^2}{4\a^4\c^2}\right)\max\left( \Vert \fe \Vert^{2}_{H^{\frac{\beta}{2}}_\delta},  \Vert \fe \Vert^{4}_{H^{\frac{\beta}{2}}_\delta} \right). 
\end{equation*}

Now, we define the constant  
\begin{equation}\label{C}
\mathfrak{C}= \frac{2}{5}\frac{C_{LT}^{5/2}}{(\a\c)^{3/2}} \left( \frac{\b^4}{4 \a^8\gamma^4}+\frac{3\b^2}{4\a^4\c^2}\right). 
\end{equation}
Then, we get the following upper bound on the quantity $\ell(n)$ given in  Definition \ref{def-global-lyapunov-expo}:  $\ell(n) \leq - \frac{\gamma}{2} n + \mathfrak{C}\max\left( \Vert \fe \Vert^{2}_{H^{\frac{\beta}{2}}_\delta},  \Vert \fe \Vert^{4}_{H^{\frac{\beta}{2}}_\delta} \right)$. Finally,  by Theorem \ref{Th-tec-upper-bound} we obtain the desired estimate (\ref{Fractal-dimension}). Theorem \ref{Prop-fractal-dim} is proven. \finpv 
\end{appendices}

\section*{Author's declarations}

Data sharing not applicable to this article as no datasets were generated or analyzed during the current
study. This work has not received any financial support. In addition, the author declares that he has no
conflicts of interest.

\end{document}